\documentclass{pcms-l}
\DeclareMathSymbol{\rtimes}{\mathbin}{AMSb}{"6F}

\newif\iffigs\figstrue

\iffigs
  \input epsf.tex
\else
\fi

\makeatletter
\def\serieslogo@{\vtop to\z@{%
  \parindent\z@ \fontsize{8}{9\p@}\bfseries
        IAS/Park City Mathematics Series\newline
        Volume \currentvolume, \currentyear\par}}
\makeatother

\makeatletter
\def\copyrightbox@{2}
\makeatother

\def\currentvolume{3}
\def\currentyear{1993}

\makeatletter
\def\lecture{\global\Monographfalse\global\Lecturetrue
  \global\let\sectionmark\@gobble 
        \addtocounter{lecturenum}1\relax
        \refstepcounter{chapter}\relax 
        {\Large\bfseries
	\raggedleft
		LECTURE {\LARGE\thelecturenum}\\
		\vspace*{3pt}%
		\thelecturename
  \endgraf}%
  \chaptermark{\thelectureseries}%
  \lecturemark{\thelecturename}%
  \addcontentsline{toc}{chapter}{Lecture \thelecturenum.\ \thelecturename}%
	\vspace{34\p@}\noindent}
\makeatother

\makeatletter
\def\thebibliography#1{\ifLecture
		\else
			\chapter*\bibname
			  \markboth{\bibname}{\bibname}%
		\fi
  \normalsize\labelsep .5em\relax
  \list{\@arabic\c@enumi.}{\settowidth\labelwidth{\@biblabel{#1}}%
  \leftmargin\labelwidth
  \advance\leftmargin\labelsep
	\bibsetup\relax
	\usecounter{enumi}}\sloppy
  \clubpenalty9999 \widowpenalty\clubpenalty  \sfcode`\.\@m}
\makeatother

\makeatletter
\def\lecturestar#1{
\global\Monographfalse\global\Lecturetrue
  \global\let\sectionmark\@gobble 
        {\Large\bfseries
	\raggedleft
		#1 \\
  \endgraf}%
  \chaptermark{\thelectureseries}%
  \lecturemark{\thelecturename}%
  \addcontentsline{toc}{chapter}{\thelecturename}%
	\vspace{34\p@}\noindent}

\def\lectureoptionstar#1#2{
\global\Monographfalse\global\Lecturetrue
  \global\let\sectionmark\@gobble 
        {\Large\bfseries
	\raggedleft
		#1 \\
                \vspace*{3pt}%
                #2
  \endgraf}%
  \chaptermark{\thelectureseries}%
  \lecturemark{\thelecturename}%
  \addcontentsline{toc}{chapter}{\thelecturename}%
	\vspace{34\p@}\noindent}
\makeatother

\numberwithin{section}{chapter}
\numberwithin{equation}{chapter}

\theoremstyle{plain}
\newtheorem{theorem}[equation]{Theorem}
\newtheorem{lemma}[equation]{Lemma}
\newtheorem*{construction}{Construction}
\newtheorem*{ConeConjecture}{The Cone Conjecture}
\newtheorem*{MonodromyTheorem}{The Monodromy Theorem}
\newtheorem*{NilpotentOrbitTheorem}{The Nilpotent Orbit Theorem}
\newtheorem*{HTmirrorconjecture}{The Hodge-Theoretic Mirror Symmetry
Conjecture}
\newtheorem*{converse}{Converse Conjecture}

\theoremstyle{definition}
\newtheorem{definition}[equation]{Definition}
\newtheorem{example}[equation]{Example}
\newtheorem*{exampleonebis}{Example \ref{exampleone} (bis)}
\newtheorem*{examples}{Examples}
\newtheorem*{remark}{Remark}
\newtheorem*{exercise}{Exercise}

\def\lhk{\mathbin{
\hbox{\vrule height1.4pt width4pt depth-1pt
\vrule height4pt width0.4pt depth-1pt}}}

\newcommand{\one}{\mbox{\rm 1\kern-2.7pt l}}

\def\mmu{
{\mathchoice
{\raise-.65ex\hbox{$l$}\!\mu}
{\raise-.65ex\hbox{$l$}\!\mu}
{\raise-.4ex\hbox{$\scriptstyle l$}\!\mu}
{\raise-.25ex\hbox{$\scriptscriptstyle l$}\mskip-4mu\mu}
}}

\newcommand{\cA}{\mathcal{A}}
\newcommand{\cB}{\mathcal{B}}
\newcommand{\C}{\mathbb{C}}

\newcommand{\cD}{\mathcal{D}}
\newcommand{\cE}{\mathcal{E}}
\newcommand{\cF}{\mathcal{F}}
\newcommand{\BH}{\mathbb{H}}
\newcommand{\cH}{\mathcal{H}}
\newcommand{\cJ}{\mathcal{J}}
\newcommand{\cK}{\mathcal{K}}
\newcommand{\MM}{\mathcal{M}}
\newcommand{\cN}{\mathcal{N}}
\renewcommand{\O}{\mathcal{O}}
\renewcommand{\P}{\mathbb{P}}
\newcommand{\Q}{\mathbb{Q}}
\newcommand{\R}{\mathbb{R}}
\newcommand{\cR}{\mathcal{R}}
\newcommand{\BS}{\mathbb{S}}
\newcommand{\cS}{\mathcal{S}}
\newcommand{\cU}{\mathcal{U}}
\newcommand{\BW}{\mathbb{W}}
\newcommand{\WW}{\mathcal{W}}
\newcommand{\Z}{\mathbb{Z}}
\newcommand{\cZ}{\mathcal{Z}}

\newcommand{\ad}{\operatorname{ad}}
\newcommand{\Ann}{\operatorname{Ann}}
\newcommand{\Aut}{\operatorname{Aut}}
\newcommand{\Flags}{\operatorname{Flags}}
\newcommand{\GL}{\operatorname{GL}}
\newcommand{\Gr}{\operatorname{Gr}}
\newcommand{\Hom}{\operatorname{Hom}}
\renewcommand{\Im}{\operatorname{Im}}
\newcommand{\Image}{\operatorname{Image}}
\newcommand{\Maps}{\operatorname{Maps}}
\newcommand{\Mov}{\operatorname{Mov}}
\newcommand{\NE}{\mathop{\overline{\text{NE}}}\nolimits}
\newcommand{\PGL}{\operatorname{PGL}}
\newcommand{\rank}{\operatorname{rank}}
\newcommand{\Res}{\operatorname{Res}}
\newcommand{\Spec}{\operatorname{Spec}}
\newcommand{\Spf}{\operatorname{Spf}}
\newcommand{\Stab}{\operatorname{Stab}}
\newcommand{\SU}{\operatorname{SU}}
\newcommand{\Sym}{\operatorname{Sym}}

\newcommand{\U}{\operatorname{U}}

\newcommand{\MMhol}[1]{\widetilde{\MM}^*_{(#1,J)}}
\newcommand{\Mcx}{\MM_{\text{cx}}}
\newcommand{\cue}[1]{\frac{q^{#1}}{1-q^{#1}}}
\newcommand{\q}{\cue{}}
\newcommand{\qq}{\cue2}
\newcommand{\qqq}{\cue3}
\newcommand{\qqqq}{\cue4}

\newcommand{\cu}[2]{\frac{{#1}^{#2} q^{#1}}{1-q^{#1}}}
\newcommand{\cuo}[1]{\frac{{#1} \cdot q^{#1}}{1-q^{#1}}}
\newcommand{\cuone}[1]{\frac{1^{#1} q}{1-q}}

\newcommand{\CP}{\C\P}
\newcommand{\dbar}{\overline{\partial}}
\newcommand{\dlog}{d\mskip0.5mu\log}
\newcommand{\DR}{\text{DR}}
\newcommand{\dt}{{\scriptscriptstyle\bullet}}
\newcommand{\suchthat}{\ |\ }
\newcommand{\Thol}{T^{(1,0)}}

\begin{document}
\frontmatter
\title{Mathematical Aspects of Mirror~Symmetry}
\author{David R. Morrison}

\address{Department of Mathematics, Box 90320, Duke University,
Durham, NC 27708-0320}
\email{drm@math.duke.edu}

\date{September 26, 1996}

\mainmatter
\chapter*{Mathematical Aspects of Mirror~Symmetry}
\auth{David R. Morrison}

\lecturename{Introduction}

\markboth{D. R. Morrison, Mathematical Aspects of Mirror
Symmetry}{Mathematical Aspects of Mirror Symmetry}

\addcontentsline{toc}{chapter}{Introduction}

\addvspace\linespacing
\noindent{\large\bfseries Introduction}\par
\addvspace{.5\linespacing}

\noindent
{\em Mirror symmetry}\/ is the remarkable discovery in string theory that
certain ``mirror pairs'' of Calabi--Yau manifolds
apparently
produce isomorphic physical theories---related by an isomorphism
which reverses the sign of a certain
quantum number---when used as backgrounds for string
propagation.
The sign reversal in the isomorphism has profound effects
on the geometric interpretation of the pair of physical theories.
It leads
to startling predictions that certain geometric invariants of
one Calabi--Yau
manifold (essentially the numbers of rational curves of
various degrees) should be
related to a completely different set of geometric invariants of
the mirror partner (period integrals of holomorphic forms).
The period integrals are much easier to calculate than the numbers of
rational curves, so this idea has been used to make very specific
predictions about numbers of curves on certain Calabi--Yau manifolds;
hundreds of these predictions have now been explicitly verified.
Why either the pair of manifolds, or these different invariants, should
have anything to do with each other is a great mathematical mystery.

The focus in these lectures will be on giving a precise mathematical
description of two string-theoretic quantities which play a primary r\^ole
in mirror symmetry: the so-called
$A$-model and $B$-model correlation functions on a Calabi--Yau manifold.
The first of these is related to the problem of counting rational curves
while the second is related to period integrals and
variations of Hodge structure.
A natural mathematical consequence of mirror symmetry is the assertion that
Calabi--Yau manifolds often come in pairs with
the property that the $A$-model correlation function of the first
manifold coincides
with the $B$-model correlation function of the second, and {\em vice
versa}.  Our goal will be to formulate this statement as a precise
mathematical conjecture.  There are other recent mathematical
expositions of mirror
symmetry, by Voisin \cite{voisin} and by Cox and Katz \cite{coxkatz}, which
concentrate on other aspects of the subject;
the reader may wish to consult those as well in order to
obtain a complete picture.

I have only briefly touched on the physics which inspired mirror symmetry
(in lectures one and eight), since there are a number of good places to read
about some of the physics background:  I recommend
Witten's address at the International
Congress in Berkeley \cite{witten:physgeom},
a book on ``Differential Topology and Quantum Field Theory'' by
Nash \cite{nash}, and the first chapter of
H\"ubsch's ``Calabi--Yau Manifolds: A Bestiary for
Physicists'' \cite{hubsch}.  There are, in addition, three collections of
papers related to string theory and mirror symmetry which contain some very
accessible expository material:
``Mathematical
Aspects of String Theory'' (from a 1986 conference at U.C. San Diego)
\cite{stringbook},
``Essays on Mirror Manifolds'' (from a 1991 conference at MSRI)
\cite{mirrorbook}, and its successor volume ``Mirror Symmetry II''
\cite{MSII}.
I particularly recommend the paper by Greene and Plesser ``An introduction
to mirror manifolds'' \cite{GP:intro}
and the paper by Witten ``Mirror manifolds and
topological field theory'' \cite{witten:mirror},
both in the MSRI volume.

This is a revised version of the lecture notes which I prepared in conjunction
with my July, 1993 Park City lectures, and which I supplemented
when delivering a
similar lecture series in Trento during June, 1994.  The field of mirror
symmetry is a rapidly developing one, and in finalizing these notes for
publication I have elected to let them remain as a ``snapshot'' of the
field as it was in 1993 or 1994, making only minor modifications to the
main text to accommodate subsequent developments.  I have, however, added a
postscript that sketches the progress which has been made in a number of
different directions since then.

\addvspace\linespacing
\noindent{\large\bfseries Acknowledgments}\par
\addvspace{.5\linespacing}

\noindent
The ideas presented here concerning
the mathematical aspects of
mirror symmetry were largely shaped
through conversations and collaborations I have had
with Paul Aspinwall, Brian Greene, Sheldon Katz,
Ronen Plesser, and Edward Witten.  It is a pleasure to thank them all for
their contributions.

I am grateful to Antonella Grassi and Yiannis Vlassopoulos for providing me
with copies of the notes they took during the lectures.  I am also grateful
to Grassi, Katz, Plesser,
and Vlassopoulos, as well as to Michael Johnson, Lisa Traynor,
and the referee of \cite{compact}, for pointing out errors in the first
drafts of these notes.

This research was supported in part by the National Science
Foundation under grant number DMS-9103827.

\tableofcontents
\chapter*{}

\lecturename{Some Ideas From String Theory}
\lecture

\markboth{D. R. Morrison, Mathematical Aspects of Mirror Symmetry}{Lecture 1.
Some Ideas From String Theory}

\label{stringtheory}

\section{String theory and quantum field theory}

The origins of the startling calculations which have led to tremendous
interest among mathematicians in the
phenomenon of ``mirror
symmetry'' lie in string theory.  String theory is a proposed model
of the physical world which idealizes its fundamental constituent
particles as one-dimensional mathematical
objects (``strings'') rather than zero-dimensional objects (``points'').
In theories such as general
relativity, one has traditionally imagined a
point as tracing out what is known as
a ``worldline'' in spacetime; the
corresponding notion in string theory is of a ``worldsheet'' which
will describe the trace of a {\em string}\/ in spacetime.  We will consider
here
only ``closed string theory'' in which the string is a closed loop; the
worldsheets are then (locally) closed surfaces: if we look at a portion
of a worldsheet which represents the history of several interacting
particles over a finite time interval, we will see a closed surface which
has a boundary consisting of a finite number of closed loops.

An early version of string theory
was proposed as a model for nuclear processes in the 1960's.
Those early investigations revealed a somewhat disturbing property:
in order to get a sensible physical theory, the
spacetime $M$ in which the string is propagating must have dimension
twenty-six.
Obviously, when we look around us, we do not see twenty-six dimensions.
A later variant which incorporates supersymmetry\footnote{I shall not
attempt to explain supersymmetry in these lectures.}
is sensible exactly when the spacetime has dimension ten---again a bit larger
than the four-dimensional spacetime which we observe.   Partly for this
reason, but primarily because a better model for nuclear processes was
found, the original research activity in string theory largely
died out in the early
1970's.

String theory was subsequently revived in the 1980's when it was
shown \cite{quantum:gravity}
that if the ten-dimensional
string theory were used to model things at much smaller
distance
scales, an apparently consistent quantum theory
 of gravity could be
produced.  (In fact, gravity is predicted as an essential ingredient of
this theory.)
This ``anomaly cancellation'' result explained how certain potential
inconsistencies in the quantum theory are avoided through an interaction
between gravity and the other forces present.  Tremendous optimism and
excitement pervaded this period, particularly
since the new model contains a rich spectrum of elementary particles at low
energies and exhibits many features one would expect of
a ``grand unified field theory'' which could describe in a
single theory all of the forces observed in nature.
The ``problem'' of ten dimensions
in this context can be
resolved by assuming that
the ten-dimensional spacetime is locally
a product $M=M^{1,3}\times M^6$ of a macroscopic four-dimensional
spacetime and a compact six-dimensional space whose size is on
the order of the Planck length ($10^{-33}$ cm).  Because this is so small
compared to macroscopic lengths,
one wouldn't expect to observe the compact
space
directly, but its effect on four-dimensional physics could be detectable
in various indirect ways.

The next step was even more remarkable for mathematicians---a group of string
theorists \cite{CHSW}
calculated that the compact six-dimensional space must
have a Ricci-flat metric on it.  (The physically relevant metric
is actually a perturbation of this Ricci-flat one.)  This is a very
restrictive property---it implies, for example, that the six-manifold
is a complex K\"ahler manifold of complex dimension three which has trivial
canonical bundle;  conversely, such K\"ahler manifolds always admit
Ricci-flat metrics.  (This had been conjectured by Calabi \cite{calabi}
in the late 1950's and proved by Yau \cite{yau} in the mid 1970's.)
These manifolds have since been named ``Calabi--Yau
manifolds;'' finding and studying them become problems in algebraic
geometry, thanks to Yau's theorem.

The model being described
here of a string propagating in a spacetime (with a specified metric)
is generally regarded as a woefully inadequate
description of the ``true'' string theory, a good formulation
of which is as yet unknown.  Indeed, if string theory is truly
 a theory of gravity as we observe it, then the theory should approximate
general relativity when the distance scale approaches macroscopic
levels.  Since the metric on spacetime is part of general relativity,
it should be a part of that ``approximation'' which is somehow to
be deduced from a solution to the ultimate ``string
equations,'' rather than being something which is put in by hand in advance.
Even the {\em topology}\/ of spacetime should be dictated by the string theory.
However,
neither these ``string equations'' nor their exact solutions are known at
present.

The Calabi--Yau manifolds and their connections with string theory have been
studied intensively for more than a decade.  In the earliest period, these
manifolds were analyzed using standard mathematical techniques, and the
results were applied in a string-theoretic context.  However, at the same
time, other advances were being made in string theory which suggested other
ways of looking at certain aspects of the theory of Calabi--Yau manifolds.
This eventually led to the discovery of a surprising new phenomenon known as
``mirror symmetry,'' in which
it was observed that different Calabi--Yau manifolds could lead to
identical physical theories in a way that implied surprising connections
between certain geometric features of the manifolds.

To explain this mirror symmetry observation in more detail, we must first
describe a few aspects of
quantum field theory and its relationship to string theory.
In classical mechanics, the worldline 
representing a
particle is required to minimize the ``action'' (which is the energy
integrated with respect to time), or more precisely, to be a stationary
path for the action functional.  Due to this ``stationary action
principle,'' the location of the path in
spacetime is completely determined by a knowledge of boundary conditions.
Other physically measurable quantities associated to the particle (which
are often represented as some kind of ``internal variables'') will also
evolve from their boundary states in a completely predictable manner, again
minimizing the action.

In quantum field theory, however, this changes.  Only the probability of
various possible outcomes can be predicted with certainty, and {\em all}\/
trajectories---not just the action-minimizing ones---contribute to the
measurement of this probability.  
The probability is calculated from
an integral over the space of
all possible paths with these initial
and final states,\footnote{There are enormous
mathematical difficulties in dealing with these ``path integrals'' or
``functional integrals,''
and they do not in general have a rigorous mathematical formulation.
Nevertheless, in the hands of skilled practitioners they can be
used to make predictions which agree with laboratory experiments to a
remarkable degree of precision.}
and the classical trajectory is recovered as the leading term in a
stationary phase approximation to the path integral.

Relativistic quantum field theories are frequently studied by treating the
theory as a small perturbation of a simple type of theory---called a free
field theory---whose functional integrals are well-understood.  For
example, the path integral describing the interaction of two charged
particles can be expanded in a perturbative series whose terms are
described by ``Feynman diagrams.'' The zeroth order term is the diagram
\iffigs
$$\vbox{\centerline{\epsfysize=2cm\epsfbox{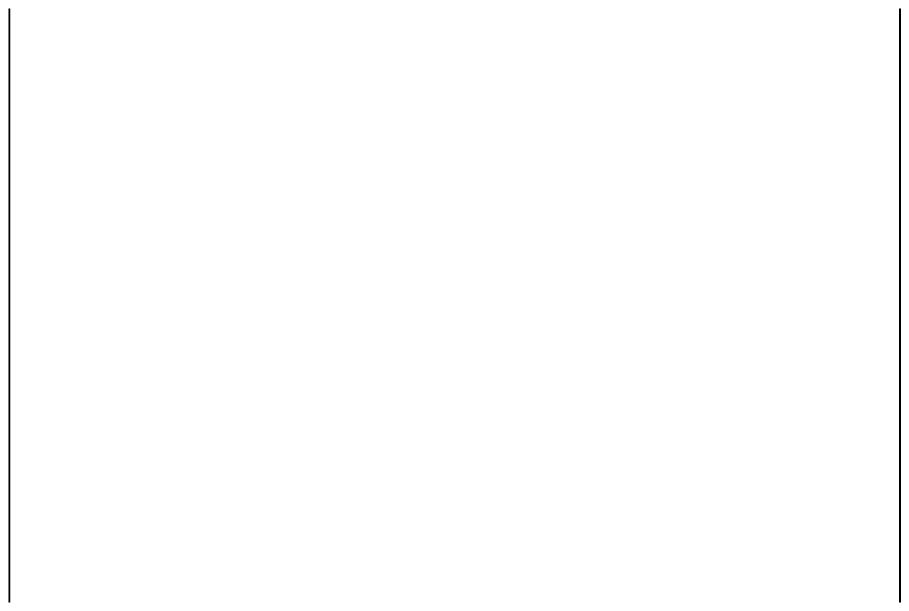}}
}$$
\else
\vskip1in
\noindent
\fi
which represents two particles which do not interact at all, the leading
perturbative correction is described by the diagram
\iffigs
$$\vbox{\centerline{\epsfysize=2cm\epsfbox{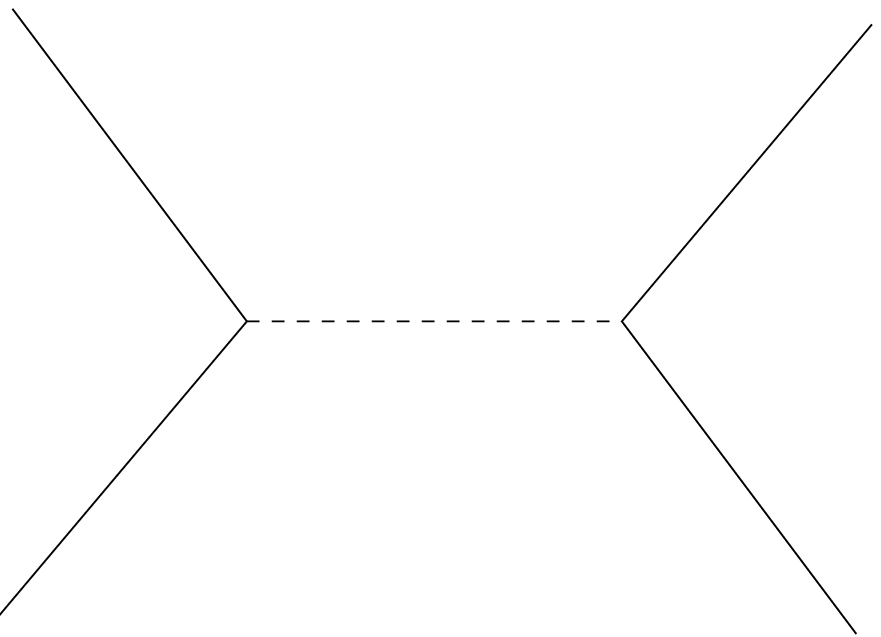}}
}$$
\else

\hglue2.5in \verb+\     /+

\hglue2.5in \verb+ ^^^^^+

\hglue2.5in \verb+/     \+

\noindent
\fi
which represents a transfer of momentum from one particle to the other via
the emission and absorption of a third particle carrying the force, and
higher order corrections involve diagrams with more complicated
topologies---loops are allowed, for example.
Such diagrams can be cut into simpler pieces, at the expense of performing
an integral over all possible intermediate states.  For example, the
interacting Feynman diagram illustrated above
can be decomposed
into two more primitive pieces,
\iffigs
$$\vbox{\centerline{\epsfysize=2cm\epsfbox{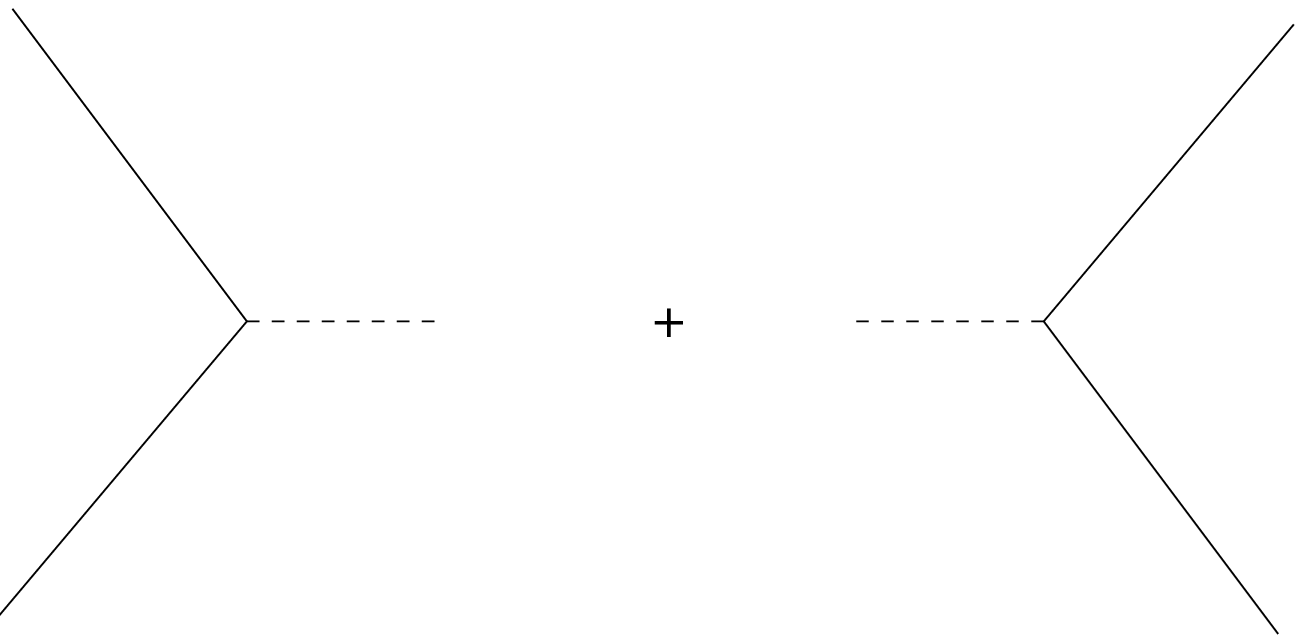}}
}$$
\else

\hglue2in \verb+\                 /+

\hglue2in \verb- ^^^     +     ^^^-

\hglue2in \verb+/                 \+

\noindent
\fi
each of which represents a fundamental ``interaction'' vertex.

In string theory, the paths are replaced by surfaces:  
the interacting 
diagram
might be represented as a sphere with four disks
removed (or perhaps as something with more complicated topology),
\iffigs
$$\vbox{\centerline{\epsfysize=2cm\epsfbox{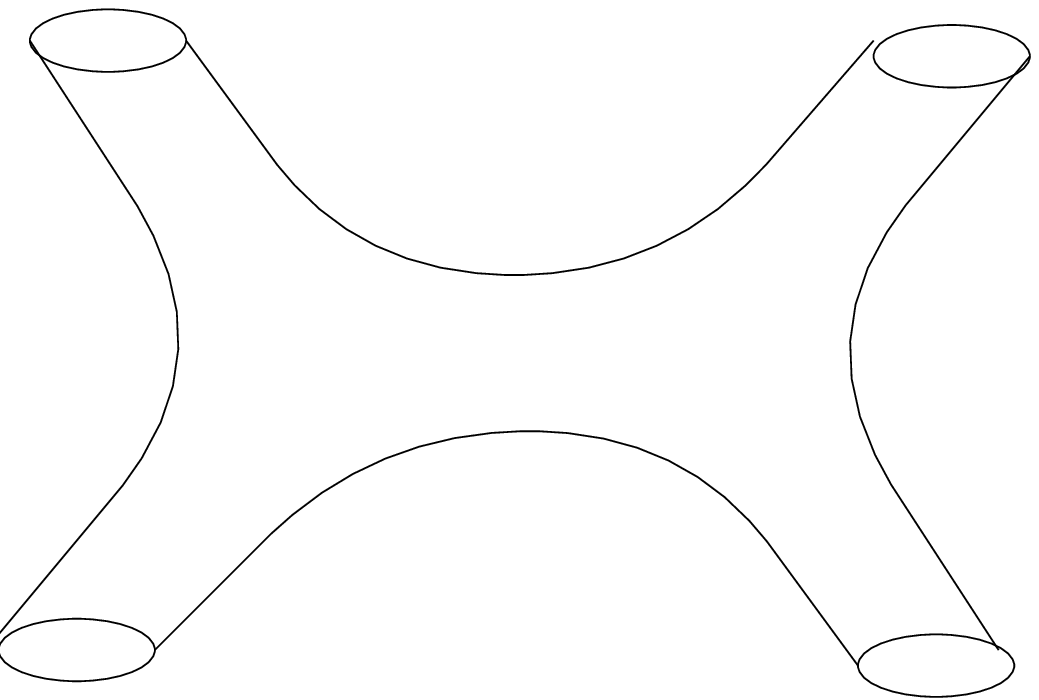}}
}$$
\fi
and this 
could also
be decomposed into more primitive pieces
\iffigs
$$\vbox{\centerline{\epsfysize=2cm\epsfbox{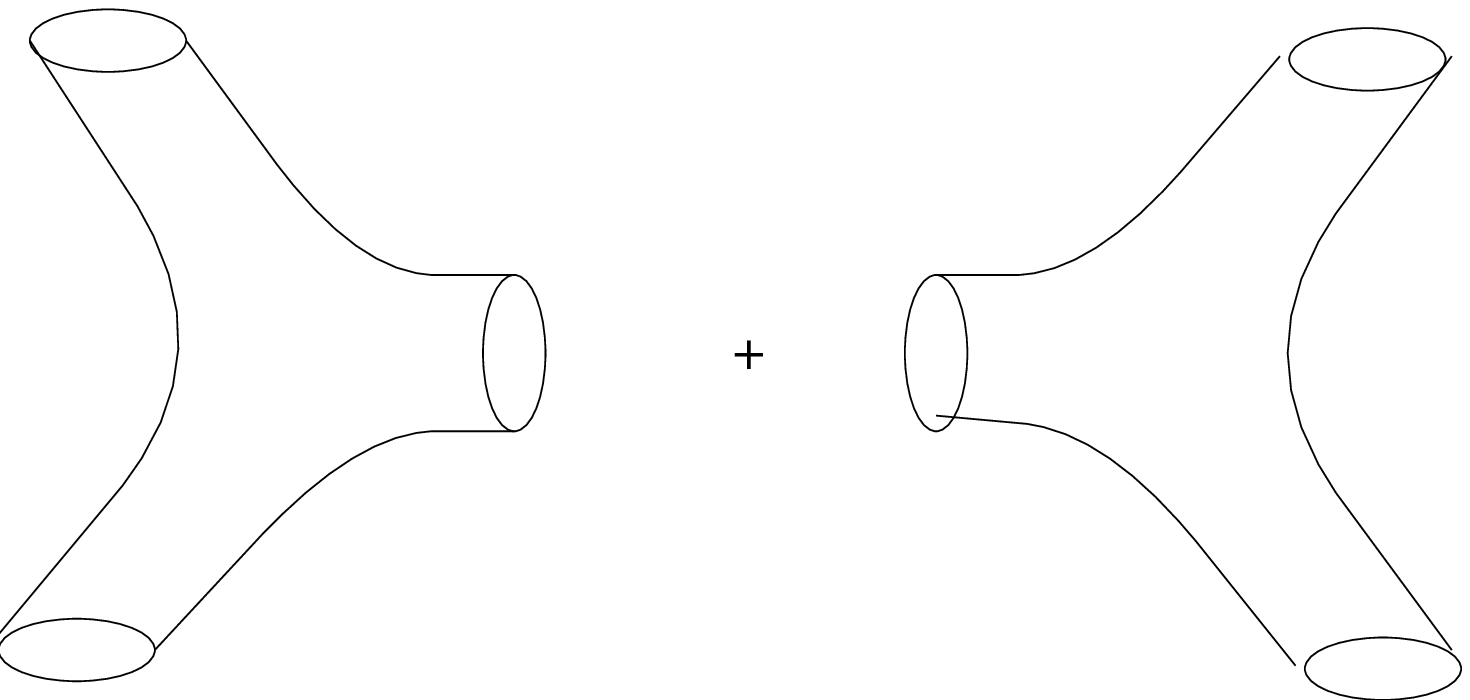}}
}$$
\fi
(called ``pair of pants'' surfaces).  One of the advantages of string theory is
that this fundamental piece, the ``pair of pants'' surface, is a smooth
surface, in contrast to the interaction vertex which introduced a
singularity into the worldline.

The methods of quantum field theory are applied to string theory in a
rather interesting way.  If we fix the spacetime and consider a string
propagating through it, the location of the worldsheet can be viewed as a
map from the worldsheet to the spacetime, and we can regard the coordinate
functions on the spacetime as functions on the worldsheet.  These spacetime
coordinate functions are then treated as the ``internal variables'' of a
two-dimensional quantum field theory---formulated on the worldsheet
itself---which captures many of the important
physical features of the string theory. 
(The functional integral in this theory involves an integration over all
possible metrics on the worldsheet as well as all possible maps from the
worldsheet to the spacetime.)  The two-dimensional quantum field theories
arising from string theory are
of a particular type known as a
{\em conformal field theory}; this means that a conformal change of metric
on the worldsheet will act as an automorphism of the theory (typically
acting linearly on various spaces of ``internal fields'' of the
theory).
The formulation in terms of conformal field theory
has turned out to be a very fruitful viewpoint for the study of string
theory.

\section{Correlation functions and pseudo-holomorphic curves}

The basic quantities which one needs to evaluate in any quantum field
theory are the {\em correlation functions}\/ which determine the
probabilities for a specified final state, given an initial state.
Specifying the initial and final states means not only specifying 
positions, but also the values of any ``internal variables'' which form a
part of the theory.  The possible initial or final
states in a conformal field theory can be
represented as
operators $\O_P$ on some fixed Hilbert space $\cH$, often referred to as
``vertex operators.''\footnote{Generally, in quantum field theories states
are represented as elements of a Hilbert space $\cH$ but in conformal field
theories there is also an operator interpretation.}
(The label $P$ indicates the position; we should in principle be
specifying initial conditions on an entire boundary circle, but in fact
it suffices to consider a
limit---within the conformal class of the given metric---in which
the circles have been shrunk to zero size and the vertex
operators are located at points.)
The conjugate transpose of an
initial state is a final state, so we often don't distinguish between
those in our notation; with these conventions, the correlation function of a
number of vertex operators is denoted by
\[\langle\O_{P_1}\O_{P_2}\dots\O_{P_k}\rangle.\]
(Note that the correlation functions are complex-valued and do not directly
calculate probabilities, but also include the phase of the
quantum-mechanical wavefunction.)

If we fix the topology of the worldsheet we must in general integrate over
the choice of metric on that worldsheet.  A conformal change of metric
leaves the correlation functions invariant, so we only need to integrate
over the set of conformal classes of metrics, i.e., over the
(finite-dimensional) moduli
space ${\MM}_{g,k}$ of $k$-punctured
 Riemann surfaces of genus $g$.

The two-dimensional quantum field theories which are related to mirror
symmetry have a subset of their correlation functions whose values do not
depend on the position of the points $P_j$ on the worldsheet; these are
called {\em topological}\/ correlation functions.
(We don't need to consider an
integral over ${\MM}_{g,k}$ in this case.)
They will be the primary
objects of interest for us.  In fact, due to the possibility of decomposing
the worldsheet into more primitive pieces, the main case to consider is the
case of three vertex operators
on surfaces of genus zero, i.e., we take $\Sigma$ to be
the ``pair of pants'' surface
$\CP^1-\{P_1,P_2,P_3\}$.
To evaluate a correlation function
\[\langle\O_{P_1}\,\O_{P_2}\,\O_{P_3}\rangle,\]
however, we must still integrate over the
infinite-dimensional space $\Maps(\Sigma,M)$ of maps from $\Sigma$
to the spacetime $M$.

To proceed further, we need to introduce the ``action'' functional on the
space of maps.  We fix Riemannian metrics\footnote{This is a ``Euclidean''
version of the theory, whose correlation functions are related by analytic
continuation to those of the
``Lorentzian'' version in which the worldsheet metric has signature
$(1,1)$.}  on both the worldsheet $\Sigma$
and the spacetime $M$ and define for any sufficiently smooth
$\varphi\in\Maps(\Sigma,M)$
\[\cS[\varphi]=\int_\Sigma\|d\varphi\|^2\,d\mu\]
using the metrics to define the norm.  (In practice, we take $M$ to be the
compact six-dimensional manifold rather than the entire space.)  The
properties of the action functional
are easier to analyze if we assume that $M$ has
some additional structure---the minimal structure needed is a symplectic
form $\omega$ and an almost-complex structure $J$ which is $\omega$-tamed.
(We will review the definitions of these
in lecture three.)  When these have been
chosen, there is a ``d-bar'' operator $\dbar_J$ on maps, and an alternate
formula for the action
\[\cS[\varphi]=\int_\Sigma\|\dbar_J\varphi\|^2\,d\mu
+\int_\Sigma\varphi^*(\omega).\]
A lower bound for $\cS[\varphi]$ in any homotopy class of maps is thus given by
$\int_\Sigma\varphi^*(\omega)$; this bound will be achieved by the
so-called pseudo-holomorphic maps---the ones for which
$\dbar_J\varphi\equiv0$. These have been extensively studied by Gromov
\cite{gromov} and others as a natural generalization of complex curves on
K\"ahler manifolds.

This action functional now appears in an integrand which is supposed to be
integrated over the infinite-dimensional space of all maps.  We will
outline the standard manipulations which are made with these functional
integrals
in physics in order to
express the correlation function as
an infinite sum of
finite-dimensional integrals.  We will subsequently use the outcome of those
manipulations to make
mathematical definitions of the corresponding quantities
in the form of a formal sum of these finite-dimensional integrals.

The topological correlation functions we are studying are to be evaluated
by a functional integral of the form
\begin{equation}\label{eq:path}
\begin{split}
\langle\O_{P_1}\dots\O_{P_k}\rangle&=
\int {\cD}\varphi\,\O_{P_1}\dots\O_{P_k}\,e^{-2\pi \,{\cS}[\varphi]}\\
&=e^{-2\pi\,\int_\Sigma\varphi^*(\omega)}
\int {\cD}\varphi\,\O_{P_1}\dots\O_{P_k}\,e^{-2\pi
\,\int_\Sigma\|\dbar_J\varphi\|^2\,d\mu}.
\end{split}
\end{equation}
(We are suppressing the ``fermionic'' part of this functional integral, which
is
actually very important, but explaining it would take us too far afield.)
The ``topological'' property of these correlation functions turns out to
imply \cite{tsm,witten:mirror} that if we introduce a parameter $t$ into
the exponent of the last functional integral in eq.~\eqref{eq:path} to produce
\[\int {\cD}\varphi\,\O_{P_1}\dots\O_{P_k}\,e^{-2\pi t
\,\int_\Sigma\|\dbar_J\varphi\|^2\,d\mu},
\]
then the resulting expression is independent of $t$ and can be evaluated in
a limit in which $t\to\infty$.  In such a limit, the only contributions to
the functional integral are the maps $\varphi$ for which
$\dbar_J\varphi\equiv0$,
i.e., the pseudo-holomorphic maps.
(This trick for reducing to a finite-dimensional integral
is known as the ``method of stationary phase.'')  The space of
pseudo-holomorphic maps in a given homotopy class is finite-dimensional, so
we have reduced the evaluation of our correlation function to an infinite
sum of finite-dimensional integrals, of the form
\begin{equation}\label{eq:reduced}
\langle\O_{P_1}\dots\O_{P_k}\rangle=\sum_{\text{homotopy classes}}
e^{-2\pi\,\int_\Sigma\varphi^*(\omega)}
\int_{\MM}\cD\varphi\,\O_{P_1}\dots\O_{P_k},
\end{equation}
where $\MM$ denotes the moduli space of pseudo-holomorphic maps in a fixed
homotopy class.
It is these finite-dimensional integrals on which we shall eventually base
our definitions.  The convergence of the infinite sum will remain an issue
in our approach, and will lead us to (in some cases) assign a provisional
interpretation to this formula as being
a formal power series only.  From the
physics one expects convergence whenever the volume of the corresponding
metric is sufficiently large.

\section{A glimpse of mirror symmetry}

If the target space $M$ for our maps is a Calabi--Yau manifold (equipped
with a Ricci-flat metric), all of
the vertex operators which participate in a given topological
correlation functions must be of one
of two distinct types.  Correlation functions involving vertex operators
of the first type are called
{\em $A$-model correlation functions}\/ while those involving vertex operators
of the second type are
known as {\em $B$-model correlation functions}\/
\cite{witten:mirror}.
(These are actually the correlation functions in two ``topological field
theories'' \cite{tsm} which are closely related to the original quantum
field theories.)
For each type, the vertex operators $\O_P$ in the quantum
field theory or topological field theory
have a geometric interpretation; we will treat the correlation
functions as functions of these geometric objects.

The $A$-model correlation functions can be defined in a much broader
context than Calabi--Yau manifolds: they can be defined for any semipositive
symplectic manifold $M$ (where semipositive roughly means that $-c_1(M)$ is
nonnegative---we will give the precise definition in lecture three).  The
vertex operators $\O_{P_j}$
in the topological field theory correspond to harmonic
differential forms $\alpha_j$
on $M$, and the correlation functions
$\langle\alpha_1\,\alpha_2\,\alpha_3\rangle$
take the form of an infinite series whose constant term---corresponding to
homotopically trivial maps from $\Sigma$ to $M$---is the familiar trilinear
function $\int_M\alpha_1\wedge\alpha_2\wedge\alpha_3$.

To evaluate the non-constant terms we need an integral over the moduli space
of pseudo-holomorphic two-spheres.  In the Calabi--Yau case, those
two-spheres are expected to be discrete (based on a formal dimension
count),
so there should be invariants which count the number of rational
curves in a given homology class.
(There are certain technical difficulties with this, as we shall see in
lecture two.)  More generally, the non-constant terms in the $A$-model
correlation
functions will be related to certain kinds of counting problems for
pseudo-holomorphic curves on a semipositive symplectic manifold.

The $B$-model correlation functions, on the other hand, require a choice of
nonvanishing holomorphic $n$-form $\Omega$
on $M$ for their definition, so they are
restricted to the Calabi--Yau case.  The vertex operators in the topological
correlation functions correspond to elements in
the space $H^q(\Lambda^p\Thol_M)$, where we use $\Thol$ to denote the
holomorphic tangent bundle of an almost-complex manifold.  (More precisely,
we should use Dolbeault cohomology to describe $H^q(\Lambda^p\Thol_M)$,
and take harmonic representatives to get the vertex 
operators in the topological
field theory.)  The ``first term'' in the correlation function
is then defined as a composition
of the standard map on cohomology groups
\[H^{q_1}(\Lambda^{p_1}\Thol_M)\times H^{q_2}(\Lambda^{p_2}\Thol_M)\times
H^{q_3}(\Lambda^{p_3}\Thol_M)  \to
H^{n}(\Lambda^{n}\Thol_M)\]
(for $p_1+p_2+p_3=q_1+q_2+q_3=n$) with some isomorphisms depending on the
choice of $\Omega^{\otimes2}$
\[H^n(\Lambda^n(\Thol_M))
\overset{\lhk\,\Omega}{\longrightarrow}
H^n(\O_M)\cong\left(H^0(K_M)\right)^*
\overset{\otimes\Omega}{\longrightarrow}
\C,\]
where the middle isomorphism is Serre duality.
(This can be written as an integral over $M$, and so can be thought of as
coming from integrating over the moduli space of homotopically trivial maps
from $\Sigma$ to $M$---this is why we identify it with the first term in an
expansion like eq.~\eqref{eq:reduced}.)
Remarkably, all of the other terms in the expansion \eqref{eq:reduced} of a
$B$-model 
correlation function are known to vanish on
physical grounds \cite{DG:exact,witten:mirror}, 
so we can calculate these correlation
functions exactly using geometry, and even use the geometric version of the
correlation function as a mathematical definition.

In brief, the idea of mirror symmetry is this.  There could be pairs of
complex manifolds $M$,
$W$ (each with trivial canonical bundle)
which produce identical physics when used for string compactification,
except that the r\^oles of the
$A$-model and $B$-model correlation functions are reversed.
In particular, this would imply the existence of isomorphisms
\[H^q(\Lambda^p(T_M^{(1,0)})^*)\cong
H^q(\Lambda^p(T_W^{(1,0)}))\]
(and vice versa), as well as formulas relating the $A$-model correlation
functions on $M$ (which count the number of rational curves) to the
$B$-model correlation functions on $W$ (which are related to period
integrals of $\Omega$).

\chapter*{}

\lecturename{Counting Rational Curves}
\lecture

\markboth{D. R. Morrison, Mathematical Aspects of Mirror Symmetry}{Lecture
2. Counting Rational Curves}

\noindent
In this lecture we begin the discussion of the problem of counting rational
curves
on a complex threefold with trivial canonical bundle (a ``Calabi--Yau
threefold'').  These curve-counting invariants will eventually
be used to formulate a mathematical version of the $A$-model correlation
functions.  In the present lecture, we focus on the problems one encounters
in formulating these invariants purely algebraically; we give a number of
examples.

Consider the deformation theory of holomorphic maps from
$\CP^1\to M$, where $M$ is a complex projective variety.
If we are given such a map $\varphi:\CP^1\to M$, then a first order variation
of that map can be described by specifying in which direction
(and at what rate) each point of the image moves.  That is,
we need to specify a holomorphic tangent vector of $M$ for every
point on $\CP^1$, or in other words, a section of
$H^0(\CP^1,\varphi^*(\Thol_M))$.  As might be expected from other
deformation problems, the obstruction
group for these deformations is
$H^1(\CP^1,\varphi^*(\Thol_M))$.  The moduli problem for such maps
will be best-behaved if the obstruction group vanishes, that is, if
$h^1(\CP^1,\varphi^*(\Thol_M))=0$.
When that is true, the moduli space will be a smooth complex manifold
of complex dimension $h^0(\CP^1,\varphi^*(\Thol_M))$.  More generally, the
Euler--Poincar\'e characteristic
\[\chi(\varphi^*(\Thol_M))=
h^0(\CP^1,\varphi^*(\Thol_M))-h^1(\CP^1,\varphi^*(\Thol_M))\]
can be regarded as the ``expected complex dimension'' of the moduli space.

Although the Euler--Poincar\'e characteristic can be easily computed from the
Riemann--Roch theorem for vector bundles,  we shall make a more
elementary calculation, based on a structure theorem for bundles on
$\CP^1$.

\begin{theorem}[Grothendieck]
Every vector bundle ${\cE}$ on $\CP^1$ can be written as a direct
sum of line bundles:
\[{\cE}\cong\O(a_1)\oplus\cdots\oplus\O(a_n).\]
\end{theorem}

Using a Grothendieck decomposition for $\varphi^*(\Thol_M)$, we can
calculate the cohomology directly.  For if
\[\varphi^*(\Thol_M)\cong\O(a_1)\oplus\cdots\oplus\O(a_n)\]
then using the fact the $h^0(\O(a))=1+a$ we find
\[h^0(\varphi^*(\Thol_M))=\sum_j
\begin{cases}
1+a_j&\text{if } a_j\ge-1\\
0&\text{if } a_j<-1
\end{cases}
\]
while since $H^1(\O(a_1)\oplus\cdots\oplus\O(a_m))\cong
H^0(\O(-2-a_1)\oplus\cdots\oplus\O(-2-a_m))^*$ we have
\[h^1(\varphi^*(\Thol_M))=\sum_j
\begin{cases}
-(1+a_j)&\text{if } -2-a_j\ge0\\
0&\text{if } -2-a_j<0
\end{cases}
\]
since $1+(-2-a_j)=-(1+a_j)$.

Taking the difference, we find
\[
\chi(\varphi^*(\Thol_M))
=\sum_j(1+a_j)=n+\sum_ja_j=\dim_\C M+\deg\varphi^*(-K_M).\]
Thus, the ``expected dimension'' is independent of the decomposition.
The same result can be obtained from Riemann--Roch.

But our calculation shows more---to
ensure vanishing of the obstruction group, we must have $a_j\ge-1$
for all $j$.
In addition to this condition on the $a_j$'s, they must also satisfy
$\max\{a_j{-}2\}\ge0$, which is seen as follows.  From the exact sequence
\[0\to\Thol_{\CP^1}\to\varphi^*(\Thol_M)\to N_\varphi\to0\]
(where $N_\varphi$ denotes the normal bundle) and the fact that
$\Thol_{\CP^1}\cong\O(2)$, we see that there must be a nontrivial
homomorphism
\[\O(2)\to\O(a_1)\oplus\cdots\oplus\O(a_m),\]
which implies that $\max\{a_j{-}2\}\ge0$ as claimed.  Without loss of
generality, we may therefore assume that $a_1\ge2.$

In the  case relevant to string theory
($K_M=0$, $\dim_\C M=3$) we then find that in order
to have vanishing obstruction group
we need
\[0=a_1+a_2+a_3\ge2-1-1=0\]
 and so $a_1=2$, $a_2=a_3=-1$.  In this
case, the moduli space of holomorphic maps will be smooth of dimension three;
if we mod out by the automorphism group $\PGL(2,\C)$, the moduli space
of unparameterized maps will be smooth of dimension $0$.  The points in
that space are what
we would like to ``count.''
We discuss some examples, drawn largely from \cite{katz:mirror}, to which
we refer the reader
for more details.

\begin{example} \label{exampleone}
{\it Lines on the Fermat quintic threefold.}
All of the lines on the Fermat quintic threefold
\[\{x_0^5+x_1^5+x_2^5+x_3^5+x_4^5=0\}\subset\CP^4\]
can be described as follows.\footnote{We use the Fermat quintic because
it is an easily-described nonsingular hypersurface, and because it will be
related to a mirror symmetry construction later on, {\em not}\/ because
Wiles announced a proof of Fermat's Last Theorem while
the 1993 Park City Institute was underway!}

\medskip

\noindent
{\it First type}\/ (375 lines):
The line described by $x_0+x_1=x_2+x_3=x_4=0$, and
others whose equations are obtained from these
by permutations and multiplication by fifth roots of unity.

\medskip

\noindent
{\it Second type}\/ (50 one-parameter families of lines):
The lines described parametrically by
\[(u,v)\mapsto(u,-u,av,bv,cv)\]
for fixed constants $a$, $b$, $c$ satisfying $a^5+b^5+c^5=0$,
and others whose parameterizations are obtained from these
by permutations and multiplication by fifth roots of unity.

\medskip

\noindent
So we see that the lines are not always finite in number, even for
smooth hypersurfaces.  (One might have suspected such a ``universal
finiteness for smooth hypersurfaces'' based on experience with cubic
surfaces---every smooth cubic surface in $\CP^3$ has precisely twenty-seven
lines.)
\end{example}

\begin{example}
{\it Lines on the general quintic threefold.}
However, if we deform from the Fermat
quintic threefold to a general one, it is possible
to show that the number of lines is finite.  The generic number of lines
can then be computed as follows.  Start from the Grassmannian
$\Gr(\CP^1,\CP^4)$ of lines in $\CP^4$.  Consider the universal bundle
$U$ whose fiber at a line $L$ is the two-dimensional
subspace $U_L\subset\C^5$ such that
$\P(U_L)=L$.  We define a bundle ${\cB}=\Sym^5(U^*)$ whose fibers
describe the quintic forms on the lines $L$.  Then every quintic
threefold $M$ determines a section $s_M\in\Gamma({\cB})$:  the
equation of $M$ is restricted to $L$ to give a homogeneous quintic there.
Clearly, the lines contained in $M$ are precisely those whose corresponding
points in the Grassmannian are zeros of the section $s_M$.

The Grassmannian $\Gr(\CP^1,\CP^4)$ has complex dimension six, and the
bundle ${\cB}$ has rank six; when things are generic, the section
$s_M$ will have finitely many zeros, which can be counted by calculating
\[\#\{L\suchthat s_M(L)=0\}=c_6({\cB})=2875.\]
\end{example}

\begin{exampleonebis}
Katz \cite{katz:mirror}
has found a way to  assign multiplicities to each of
the isolated lines, and one-parameter
families of lines, on the Fermat quintic threefold.
His multiplicity assignment for each of the 375 isolated lines is ``5,''
and that for each of the 50 one-parameter families is ``20.''
Thus, the total count is
\[5\cdot375+20\cdot50=2875.\]
Katz's methods of assigning multiplicities are not yet completely 
general,\footnote{See the ``Postscript: Recent Developments'' section for
the current status.}
but they do hold out the hope that a ``count'' of rational curves
might be made even in cases when the actual number of curves is not
finite.
\end{exampleonebis}

\begin{example}
{\it Conics on the general quintic threefold.}
We can make a similar calculation for conics on the general quintic
threefold.  The key observation is that every conic spans a $\CP^2$,
so the starting point for describing them is the Grassmannian
$\Gr(\CP^2,\CP^4)$.  We need the bundle over the Grassmannian whose
fiber is the set of conics in the $\CP^2$ in question:  this is described
by $\P(\Sym^2(U^*))$, where $U$ is the universal subbundle as before.

The space $\P(\Sym^2(U^*))$ contains degenerate conics (pairs of lines,
and double lines) as well as smooth conics.  However, if $M$ is sufficiently
general, then the actual locus of conics which lie in $M$ will be finite
in number, and contain only smooth conics.

The vector bundle which will get a section $s_M$ for every quintic $M$
is the bundle ${\cB}:=\Sym^5(U^*)/(\Sym^3(U^*)\oplus\O_\P(-1))$.
This describes
the effect of restricting the quintic equation to the conic:  one gets
a quintic equation on the $\CP^2$, but must mod out by those quintics
which can be written as the product of a cubic (the $\Sym^3(U^*)$ factor)
and the given conic.

We have $\dim_\C\P(\Sym^2(U^*))=\rank{\cB}=11$, so the computation
is made by calculating:
\[\#\{C\suchthat s_M(C)=0\}=c_{11}({\cB})=609250.\]
\end{example}

\begin{example}
{\it Twisted cubics on the general quintic threefold.}
The problem gets more difficult for twisted cubics.  Again, we can
look at the linear span (a $\CP^3$) and begin by considering a
Grassmannian $\Gr(\CP^3,\CP^4)$.  But this time we must use a bundle
${\cH}\to\Gr(\CP^3,\CP^4)$ whose fibers are isomorphic to
the Hilbert scheme
of twisted cubics in $\CP^3$.  That scheme contains limits which
are quite complicated.  (For example, there is a limit which is a
nodal plane curve with an embedded point at the node which points out
of the plane:
\iffigs
$$\vbox{\centerline{\epsfysize=2cm\epsfbox{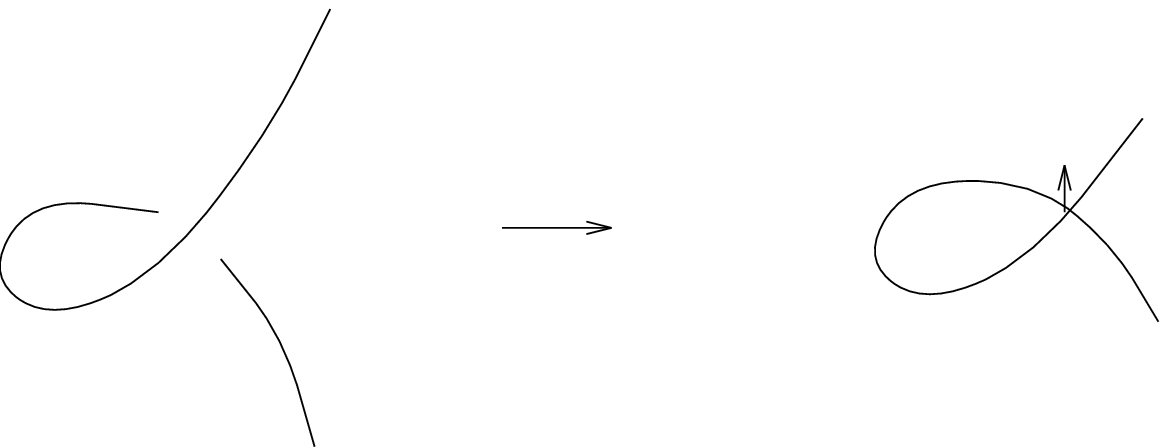}}
}$$
\fi
see Hartshorne \cite{Hartshorne}, pp.~259--260.)
Although the bundle ${\cB}$ and the
section $s_M$ can be
defined and understood at points representing
smooth twisted cubics, their extension to
the locus of degenerate cubics is by no means easy.

Ellingsrud and Str\o mme \cite{ES}
have, however, carried this out, and they find
that the number of twisted cubics on the general quintic threefold
is 317206375.
\end{example}

Clemens \cite{Clemens:AJ}
has conjectured that the general quintic threefold will have
only a finite number of rational curves of each degree, and that all
of them will satisfy $\varphi^*(\Thol_M)=\O(2)\oplus\O(-1)\oplus\O(-1)$.
This has been verified up through degree nine by Katz
\cite{katz:degree7}, Johnsen--Kleiman \cite{JK:nine} and
Nijsse \cite{nijsse}, and the prospects
are good for degrees as high as twenty-four \cite{JK:high}.  However, as we
have seen, making the calculation of the number becomes very difficult
past degree two.  In fact, for degree greater than three,
effective techniques for
calculating this number are not presently known.

We now turn to another example which demonstrates that we cannot always
expect finiteness, even for the generic deformation of a given threefold
with trivial canonical bundle.

\begin{example}
{\it Rational curves on double solids.}
We let $M$ be the double cover of $\CP^3$, branched along a general
surface $S$ of
degree eight in $\CP^3$; the double cover map is denoted by
$\pi:M\to\CP^3$.  We let $\pi^*(H)$ be the pullback of a hyperplane
$H$ from $\CP^3$; the {\em degree}\/ of a rational curve $C$ will mean
$\pi^*(H)\cdot C$.

To find ``lines'' on $M$, that is, curves $L$ with $\pi^*(H)\cdot L=1$
we consider their images $\pi(L)$.  Since $\pi^*(H)$ meets $L$ in a
single point $P$, $H$ meets $\pi(L)$ in the single point $\pi(P)$.
Thus, $\pi(L)$ must itself be a line.  But its inverse image on $M$
will necessarily have two components:  $\pi^{-1}(\pi(L))=L+L'$.
In order to have this splitting into two components, the line
$\pi(L)$ must be tangent to $S$ at every point of intersection with
$S$, i.e., it must be four-times tangent to $S$.  Now the Grassmannian
$\Gr(\CP^1,\CP^3)$ has dimension four, and it is one condition to be
tangent to a surface, so the dimension of the set of four-tangent lines
is nonnegative, and can be expected to be equal to zero.  
(In fact, it turns out to
equal zero as expected, when $S$ is general.)  The number of such lines in
the Grassmannian
can be calculated with the Schubert calculus; it turns out to be
14752.  The corresponding count of lines on $M$ is 29504.

Finding ``conics'' on $M$ is a different story, as has been observed
by Katz and by Koll\'ar.  Given a curve $C$ with
$\pi^*(H)\cdot L=2$, there are two possibilities for $\pi(C)$: it could
be a line, or it could be a conic.  In the latter case, the conic $\pi(C)$
must
be eight-times tangent to $S$.  But in the former case, in order to have
an irreducible double cover with a rational normalization, the line
$\pi(C)$ must be three-times tangent to $S$.  By our previous dimension
count, there is at least a one-parameter family of such lines for any
choice of $S$.
\end{example}

So we won't always have a finite number of things to ``count,'' even
if we perturb to a general member of a particular family.  And there
is an additional difficulty if we wish to count maps from $\CP^1$ to
$M$ when multiple covers are allowed, as the next example shows.

\begin{example}\label{ex:multiple}
{\it Multiple covers.}
Suppose that $\varphi:\CP^1\to M$ is generically one-to-one, but that
we consider a map $\varphi':=u\circ\varphi$, where $u:\CP^1\to\CP^1$
is a covering of degree $\mu$.
Even if $\varphi^*(\Thol_M)=\O(2)\oplus\O(-1)\oplus\O(-1)$,
we will get a bad splitting of the pullback via the new map:
\[\varphi'{}^*(\Thol_M)=\O(2\mu)\oplus\O(-\mu)\oplus\O(-\mu).\]
Furthermore, the dimension of the moduli space can be calculated:
the moduli space of maps $u:\CP^1\to\CP^1$ of degree $\mu$ has dimension
$2\mu{+}1$.  So we see that the dimension of the space of maps will go
up and up.
\end{example}

To handle cases such as multiple covers, ``virtual'' numbers of curves must
be introduced;  Katz's approach to this is to use excess intersection
theory \cite{fulton:intersection}.  However, this introduction of ``virtual''
numbers
leads to another complication, as our final example shows.

\begin{example}\label{ex:negative}{\it Negative numbers of curves}\/ (see
\cite{2param2}, section 8).
There are cases in which the ``virtual'' number of curves is negative.
In general, when the parameter space $B$ for a family of curves is smooth
of dimension $b$, the virtual number of curves should be the top Chern
class of the holomorphic
cotangent bundle $c_b((T^{(1,0)}_B)^*)$.  If $M$ is a complex threefold
with $K_M=0$
which contains $\CP^2$ as a submanifold (which can arise from
resolving a $\Z/3\Z$-quotient singularity, for example), then the lines
on $\CP^2$ are parameterized by $\CP^2$ and have virtual number
$c_2((T^{(1,0)}_{\CP^2})^*)=3$, but the conics on $\CP^2$, being parameterized
by $\CP^5$, have virtual number $c_5((T^{(1,0)}_{\CP^5})^*)=-6$.
This negative value actually agrees with the predictions of mirror symmetry
as shown in \cite{2param2}.
\end{example}

\chapter*{}

\lecturename{Gromov--Witten Invariants}
\lecture

\markboth{D. R. Morrison,
Mathematical Aspects of Mirror Symmetry}{Lecture
3. Gromov--Witten Invariants}

\section{Counting curves via symplectic geometry}

The difficulties we encountered in trying to count rational curves on a
Calabi--Yau threefold can be avoided by enlarging the category we are
considering, and using Gromov's theory of pseudo-holomorphic spheres in
symplectic manifolds \cite{gromov}.  This approach has the advantage that
the almost-complex structure can be slightly
perturbed to make the number of such
spheres finite, and the finite number so obtained is independent of the
choice of small perturbation.

Let $(M,\omega)$ be a compact
{\em symplectic manifold}\/ of  dimension $2n$.
This means that $M$ is a compact oriented differentiable
manifold of (real) dimension $2n$
and $\omega$ is a closed real two-form on $M$ which is nondegenerate in
the sense that its $n^{\text{th}}$ exterior power $\omega^{\wedge n}$
is nonzero at every point.

An {\it almost complex structure}\/ on a manifold $M$ is a map 
$J:T_M\to T_M$ whose square is $-1$.  If we complexify the
tangent spaces, we get $T_{M,p}\otimes\C=
T_{M,p}^{(1,0)}\oplus T_{M,p}^{(0,1)}$,
the decomposition into $+i$ and $-i$ eigenspaces for $J$.
If these subspaces are closed under Lie bracket, we say that the
almost-complex structure is {\em integrable}; in this case, these
subspaces give $M$
 the structure of a complex manifold.

If $(M,\omega)$ is a symplectic manifold,
an almost-complex structure $J$ on $M$ is said to be {\em $\omega$-tamed}\/ if
$\omega(\xi,J\xi)>0$ for all nonzero $\xi\in T_pM$.
If we have fixed an ($\omega$-tamed) almost-complex
structure $J$ on $M$, and $\varphi$ is a differentiable map from
$S^2$ to $M$, we define
\[\dbar_J\varphi=\frac12(d\varphi+J\,d\varphi\,J_0),\]
where $J_0$ is the standard almost-complex structure on $S^2$.

The main example we have in mind is this:
 $M$ is a compact K\"ahler manifold, $\omega$ is the K\"ahler form,
and $J$ is an $\omega$-tamed perturbation of
 the original complex structure on $M$.

\begin{definition}[McDuff \cite{McD:contact}]
$(M,\omega)$ is {\em semipositive}\/ if there is no map $\varphi:S^2\to M$
satisfying
\[\int_{S^2} \varphi^*(\omega)>0, \quad \text{and} \quad
3-n\le\int_{S^2}\varphi^*(-K_M)<0,\]
where we are writing $-K_M$ as in algebraic geometry to indicate the
first Chern class $c_1(M)$, which may be represented as a two-form.
\end{definition}

\begin{examples} Here are three ways of producing semipositive symplectic
manifolds.
\begin{enumerate}
\item If $K_M=0$ (the Calabi--Yau case)
then $(M,\omega)$ is  semipositive for any $\omega$.
\item If $M$ is a complex projective manifold with $|{-}K_M|$ ample
(a ``Fano variety''),
then we can take $\omega=-K_M$ to produce a semipositive $(M,\omega)$.
\item If $n\le3$ then $M$ is automatically semipositive.
\end{enumerate}
\end{examples}

Because it is sometimes difficult to check whether a homology
class $\eta$ is represented as the image of a map $\varphi:S^2\to M$
 we also introduce a variant of this property.

\begin{definition}
$(M,\omega)$ is {\em strongly semipositive}\/ if there is no class
 $\eta\in H_2(M,\Z)$
satisfying
\[\omega\cdot \eta>0, \quad \text{and} \quad 3-n\le(-K_M)\cdot \eta<0.\]
All three of our examples satisfy this stronger property as well.
\end{definition}

Fix a homology class $\eta\in H_2(M,\Z)$.  As we saw in example
\ref{ex:multiple}, there are technical problems caused by ``multiple-covered''
maps---maps whose degree onto the image is greater than one.
Let us call a map {\em simple}\/ if its degree onto its image is one.
We let $\Maps^*_\eta(S^2,M)\subset \Maps_\eta(S^2,M)$ denote the subset
of simple maps with fundamental class $\eta$.  We also let
 $\Maps^*_\eta(S^2,M)_{(p)}$ be the set of simple
differentiable maps $S^2\to M$
with fundamental class $\eta$
whose derivative lies in $L_p$.
Using an appropriate Sobolev norm, $\Maps^*_\eta(S^2,M)_{(p)}$ can be given the
structure of a Banach manifold.
We can then regard $\dbar_J$
as  a section of the bundle $\WW\to\Maps^*_\eta(S^2,M)_{(p)}$ whose fibers are
\[\WW_\varphi:=H^0_{(p)}(S^2,{\cA}^{(0,1)}_{S^2}
\otimes\varphi^*(\Thol_M)),\]
where the subscript $(p)$ denotes $L_p$-cohomology, and
${\cA}^{(0,1)}_{S^2}$ denotes the sheaf of $(0,1)$-forms on $S^2$
(with respect to the complex structure $J_0$).

The key technical properties we need are summarized in the following two
theorems.

\begin{theorem}[McDuff \cite{McD:examples}]
If $J$ is generic,
then
\[\MMhol\eta:=
\{\varphi\in\Maps^*_\eta(S^2,M)_{(p)}\suchthat\dbar_J\varphi=0\}\]
is a smooth manifold of dimension
\[\dim_\R\MMhol\eta=2\,\chi(\varphi^*(\Thol_M)).\]
(The dimension is calculated using the Atiyah--Singer index theorem,
which yields the same result as the Riemann--Roch theorem did in algebraic
geometry.)
\end{theorem}

\noindent
(This theorem would have failed if we had allowed multiple-covered maps
to be included.)

The next theorem is due to Gromov \cite{gromov}, based on some techniques
of Sacks--Uhlenbeck \cite{SacksUhlenbeck} and with further improvements by
several authors \cite{PW,wolfson,rye}.
(We refer the reader to those papers for a more precise statement.)

\begin{theorem}
$\MMhol\eta$ can be compactified by using limits of graphs of maps;
this compactification has good properties.
\end{theorem}

In the  case relevant to string theory
 in which $M$ is a projective manifold with $K_M=0$
of complex
dimension three, we find that for generic $J$, the (real) dimension of
$\MMhol\eta$ is six, and the dimension of
\[\MM^*_{(\eta,J)}:=\MMhol\eta/\PGL(2,\C)\]
is zero.  The space $\MM^*_{(\eta,J)}$ itself is already compact in this
case; the number of points in that space
is our desired invariant.  (These points may need to be counted with
multiplicity, or with signs.)
This invariant counts the number of rational curves
(of fixed topological type) on $M$ with respect to its original complex
structure, if that number is finite; it can be used as a substitute
for that count in the general case.\footnote{It has not yet been verified
that Katz's method of assigning multiplicities to positive-dimensional
components in the algebro-geometric context produces the same results
as this method from symplectic geometry. Because of the need to include
signs in certain circumstances, this invariant can even accommodate the
``negative virtual numbers'' which occurred in example \ref{ex:negative}.}

To describe the invariants in situations more general than complex
threefolds with trivial canonical bundle, we must introduce the
oriented bordism group $\Omega_*(M)$.  The elements of $\Omega_k(M)$
are equivalence classes of
pairs $(B^k,F)$ consisting of a compact oriented differentiable
manifold $B$ of dimension $k$ (but not necessarily connected), together
with a differentiable map $F:B^k\to M$.  We say that $(B^k,F)\sim0$
if there exists an {\em oriented bordism}\/ $(C^{k+1},H)$:  i.e.,
a differentiable
manifold $C$ of dimension $k{+}1$ and a differentiable map $H:C^{k+1}\to M$
with $\partial C^{k+1}=B^k$ and $H|_{B^k}=F$.  We add elements of
$\Omega_k(M)$ by means of disjoint union:
$(B^k_1,F_1)+(B^k_2,F_2)=(B^k_1\cup B^k_2,F_1\cup F_2)$; the additive
inverse is given by reversing orientation.

The oriented bordism group $\Omega_*(M)$ is a module over the Thom
bordism ring $\Omega$ (consisting of oriented manifolds modulo oriented
bordisms, with no maps to target spaces) via
\[N^j\cdot(B^k,F)=(N^j\times B^k,G)\]
where $G(x,y)=F(y)$.

\begin{theorem}[Thom \cite{Thom}, Conner--Floyd \cite{CF}]\quad

\begin{enumerate}
\item
$\Omega_*(M)\otimes\Q\cong H_*(M,\Q)\otimes\Omega$.

\item
If $H_*(M,\Z)$ is torsion-free, then
$\Omega_*(M)\cong H_*(M,\Z)\otimes\Omega$.
\end{enumerate}
\end{theorem}

To describe our basic invariants, we choose three classes $\alpha_1$,
$\alpha_2$, $\alpha_3$ in $\Omega_*(M)$ represented by elements
$(B^{k_1}_1,F_1)$, $(B^{k_2}_2,F_2)$, $(B^{k_3}_3,F_3)$,
and let $Z_j=\Image(F_j)$.  We call the  invariants defined below the
{\em Gromov--Witten invariants}, since it was Witten \cite{tsm} who pointed out
how Gromov's study of $\MMhol\eta$ could be used in principle to describe
invariants relevant in topological quantum field theory.  The detailed
construction of these invariants was recently carried out by Ruan
\cite{ruan}.  There are two cases to consider, with one being
more technically challenging than the other.\footnote{To simplify the
exposition, we have altered Ruan's description of the second case, ignored
the necessity of passing to the inhomogeneous $\dbar$ equation (introduced
already by Gromov \cite{gromov}), and built into our definition
the so-called ``multiple cover
formula'' expected from the physics \cite{CDGP,aspmor}.  (This latter step
is now justified thanks to a theorem of Voisin \cite{voisin:multiple};
there is also a related result of Manin \cite{Manin}.)  We
are also abusing notation somewhat by using $\Phi_\eta$ in both cases,
since the second case is actually related to Ruan's $\widetilde\Phi_\eta$
invariant.}

\begin{construction}[Ruan]
Let  $\eta\in H_2(M,\Z)$,
 let $\alpha_j=(B^{k_j},F_j)$ be a bordism
class, and let $Z_j=\Image(F_j)$, for $j=1,2,3$.
Suppose that $\sum_{j=1}^3(2n-k_j)=2n-2K_M\cdot \eta$, where $\eta$ is the
class of the
image of $\varphi$, and suppose that the almost-complex structure $J$ is
generic.
\begin{itemize}
\item[(a)]
If $-K_M\cdot\eta>0$, then
\[\{\varphi\in\MMhol\eta\suchthat \varphi(0)\in Z_1,
\varphi(1)\in Z_2, \varphi(\infty)\in Z_3\}\]
is a finite set.  Let $\Phi_{\eta}(\alpha_1,\alpha_2,\alpha_3)$
denote the signed number of points in this set,
with signs assigned according
to orientations at the specified points of intersection.
\item[(b)]
If $-K_M\cdot\eta=0$, then there exists an integer
$\Phi_{\eta}(\alpha_1,\alpha_2,\alpha_3)$ which agrees with
\[\#\{\text{generically injective}\ \varphi\in\MMhol\eta\suchthat
\varphi(0)\in Z_1,
\varphi(1)\in Z_2, \varphi(\infty)\in Z_3\}\]
(counted with signs) whenever the latter makes sense.
(The signs are all positive if the almost-complex structure is integrable.)
\end{itemize}
These invariants $\Phi_{\eta}(\alpha_1,\alpha_2,\alpha_3)$
depend only
on the bordism classes $\alpha_1$, $\alpha_2$, $\alpha_3$, and do
not change under small variation of $J$.
\end{construction}

\section{Simple properties of Gromov--Witten invariants}

In spite of the fact that we needed to pass to bordism to ensure that
the Gromov--Witten invariants are well-defined, their dependence
on bordism-related phenomena is minimal.  In fact, Ruan checks that
the invariants
are trivial with respect to the $\Omega$-module structure on
$\Omega_*(M)$, and so it follows from the theorem of Thom and
Conner--Floyd that we get a well-defined $\Q$-valued invariant on rational
homology $H_*(M,\Q)$.  If $M$ has no torsion in
homology, we even get an integer-valued invariant on integral homology.

The Gromov--Witten invariants $\Phi_\eta(\alpha_1,\alpha_2,\alpha_3)$
will vanish if  $\alpha_1$ corresponds to a class of real
codimension zero or one.  This is easy to see---if there are any
elements in the set
\[\{\varphi\in\MMhol\eta\suchthat \varphi(0)\in Z_1,
\varphi(1)\in Z_2, \varphi(\infty)\in Z_3\}\]
then the intersection of the image of $\varphi$
with the image $Z_1$ of $F_1$ has real dimension
two or one.  By varying the location of
$\varphi(0)$, we will produce a two-{} or one-parameter family of
maps.  This contradicts the set being finite; thus, the set must be
empty and the invariant vanishes.

Note what happens to the Gromov--Witten invariants in the case of
interest to string theory ($\dim_\C M=3$, $K_M=0$):  the only relevant
invariants are those with $k_1=k_2=k_3=4$.  (This is because
$k_j\le 4$ to get a nonzero invariant,
so that $6=\sum (6-k_j)\ge\sum_{j=1}^3 2=6$, which implies
that each $k_j$ is $4$.)
The possible location of $0$ under a generically injective
map is easy to spot:  the image
curve $\varphi(S^2)$ is some rational curve on $M$, and meets the four-manifold
$Z_1$ in precisely $\#(Z_1\cap \eta)=\alpha_1\cdot\eta$ points; we can choose
any of these
for the image of $0$.  Similar remarks about the images of $1$ and
$\infty$ lead to the calculation:
\[\Phi_\eta(Z_1,Z_2,Z_3)
=(\alpha_1\cdot\eta)(\alpha_1\cdot\eta)(\alpha_1\cdot\eta)\,
\#(\MM^*_{(\eta,J)}).\]

\section{The $A$-model correlation functions}

Although we have defined the Gromov--Witten invariants for oriented
bordism classes, we will now use them in cohomology instead.  As
previously remarked, thanks to the triviality of the invariants
under the $\Omega$-module structure, if we tensor with $\Q$ we can
move the invariants to homology (and then by Poincar\'e duality,
to cohomology).  This is at the expense of possibly allowing
them to become $\Q$-valued on integral classes.  One hopes that they
will remain integer valued on integer classes, but this has not
yet been established.  Therefore, we will give a presentation using
$\Q$-coefficients, but the reader should bear in mind that most of the formulas
are expected to be valid with integer coefficients if one uses integer
cohomology classes.

In brief, the bordism class of $\alpha=(B^k,F)$ gives rise to a homology class
$[Z]\in H_k(M,\Z)$ (using the image
 $Z$ of $F$ to represent the class),
and by duality to a cohomology class $\zeta=[Z]^\vee\in H^{2n-k}(M,\Q)$.
(Our retreat to $\Q$-coefficients will be in part because we do not know
that every integer cohomology class can be so represented.)  We extend
the definition of Gromov--Witten invariants to cohomology by defining
\[\Phi_\eta(\zeta_1,\zeta_2,\zeta_3):=\Phi_\eta(\alpha_1,\alpha_2,\alpha_3)\]
when $\alpha_j=(B^{k_j}_j,F_j)$ and $\zeta_j=[\Image(F_j)]^\vee$;  then
extend by linearity to all of $H^*(M,\Q)$.

Our ``$A$-model correlation functions'' are then built from the Gromov--Witten
invariants, following a calculation from the physics literature
\cite{strominger,DSWW,CDGP,aspmor}.
There is
a certain danger in using the {\em outcome}\/ of a physics
calculation as a {\em definition}---later, the physicists may become
interested in a slightly different problem, whose outcome is radically
different from the original one, and we mathematicians will find
that our definitions are inadequate.

Nevertheless, we will go ahead and
 define the $A$-model correlation functions.  These are trilinear
functions on the cohomology $H^*(M,\Q)$ defined by:
\begin{equation}\label{A:correlation}
\begin{split}
\langle\zeta_1\,\zeta_2\,\zeta_3\rangle:=
(\zeta_1\cup\zeta_2\cup\zeta_3)|_{[M]}\ \
&+\sum_{\substack{\eta\in H_2(M,\Z),\\{-}K_M\cdot\eta>0}}
  \Phi_\eta(\zeta_1,\zeta_2,\zeta_3)\,q^{\eta}\\
&+\sum_{\substack{0\ne\eta\in H_2(M,\Z),\\{-}K_M\cdot\eta=0}}
\Phi_\eta(\zeta_1,\zeta_2,\zeta_3)\,\sum_{m=1}^\infty q^{m\eta}
\end{split}\end{equation}
It is sometimes convenient to formally sum the geometric series in
the final term, and write $q^\eta/(1-q^\eta)$ in place of
$\sum_{m=1}^\infty q^{m\eta}$, in which case eq.~\eqref{A:correlation} becomes
\begin{equation}\label{A:correlationbis}
\begin{split}
\langle\zeta_1\,\zeta_2\,\zeta_3\rangle:=
(\zeta_1\cup\zeta_2\cup\zeta_3)|_{[M]}\ \
&+\sum_{\substack{\eta\in H_2(M,\Z),\\{-}K_M\cdot\eta>0}}
  \Phi_\eta(\zeta_1,\zeta_2,\zeta_3)\,q^{\eta}\\
&+\sum_{\substack{0\ne\eta\in H_2(M,\Z),\\{-}K_M\cdot\eta=0}}
\Phi_\eta(\zeta_1,\zeta_2,\zeta_3)\,
\frac{q^{\eta}}{1-q^{\eta}}
\end{split}\end{equation}

The terms with $K_M\cdot\eta=0$ have been separated out because they are
where the multiple-covered maps cause the greatest difficulty.
Heuristically, the coefficients in these functions (as we have defined
them) are expected to count the simple maps only.

The symbol $q^\eta$ which appears in these formulas has not yet been
defined.  In fact, there are two natural interpretations of
eq.~\eqref{A:correlation}, one algebraic and one geometric, and we consider
them in turn in the next two lectures.

\chapter*{}

\lecturename{The Quantum Cohomology Ring}
\lecture

\markboth{D. R. Morrison,
Mathematical Aspects of Mirror Symmetry}{Lecture
4. The Quantum Cohomology Ring}

\section{Coefficient rings}

There are several possible ways to interpret the ``$A$-model correlation
functions'' defined by eq.~\eqref{A:correlation}.  In this lecture, we
will focus on the algebraic interpretation, in which the symbol $q^\eta$
can be regarded as an element of a group ring or
semigroup ring.\footnote{I am grateful to A.~Givental for
pointing out the relevance of group rings.}
Recall that for any commutative
semigroup ${\cS}$ and any commutative ring $R$, the
{\em semigroup ring of ${\cS}$ with coefficients in $R$}\/ is
the ring
\[R[q;{\cS}]:=\left\{\sum_{\eta\in{\cS}}a_\eta q^\eta\ |\
a_\eta\in R\text{ and } \{\eta\ |\ a_\eta\ne0\}\text{ is finite}\right\}.\]
The symbol $q$ serves as a placeholder, translating the semigroup
operation (usually written additively) into a multiplicative structure
on a set of monomials.
If ${\cS}$ is a group, this coincides with the usual ``group ring''
construction.

In the case of a Fano variety, the sum in 
eq.~\eqref{A:correlation} defining the
$A$-model correlation function is
finite, and we can regard it as taking values
in the rational group ring\footnote{If we knew
that the Gromov--Witten invariants were integers, we could use the
integral group
ring $\Z[q;H_2(M,\Z)]$.  But when we passed from bordism to cohomology
we lost control of the integer structure.}\ \
$\Q[q;H_2(M,\Z)]$.
To be more concrete, if we assume for simplicity
that $H_2(M,\Z)$ has no torsion, and
choose a basis
$e_1$,\dots,$e_r$ of
$H_2(M,\Z)$,
then writing $\eta=\sum a^je_j$
we may associate to $\eta$ the
 rational monomial $q^{\eta}\in\Q(q_1,\dots,q_r)$ defined by
\[\log q^{\eta}=\sum a^j\log q_j .\]
(One can also write this multiplicatively:
\[q^{\eta}=\prod (q_j)^{(a^j)}\]
but then great care is required in distinguishing exponents
from superscripts.)

If we choose our basis so that
the coefficients $a^j$ are nonnegative for
all classes $\eta$
which
have nonvanishing Gromov--Witten invariants
$\Phi_\eta(\zeta_1,\zeta_2,\zeta_3)$ for some $\zeta_1$, $\zeta_2$,
$\zeta_3$, then each $q^{\eta}$ occurring in eq.~\eqref{A:correlation}
is a {\em regular}\/ monomial, i.e., $q^{\eta}$ belongs to
the polynomial ring
$\Q[q_1,\dots,q_r]$,
and we can calculate eq.~\eqref{A:correlation} in that ring.

In the Calabi--Yau case in which $K_M=0$, the sum in 
eq.~\eqref{A:correlation} is not
finite and we must work harder.
The simplest interpretation would be to simply allow infinite sums
$\sum a_\eta q^\eta$ as formal expressions.
However, in order to construct quantum
cohomology (which we shall do in the next section)
we need the values of the correlation function to lie in
a {\em ring}.
In the definition of semigroup rings one restricts to finite sums in order
to ensure that the partial sums which occur in the expansion of a product
will be finite.  That finiteness can still be guaranteed for products of
infinite sums if
the semigroup satisfies a special property, given below.

We say that a semigroup ${\cS}$ has the {\em finite partition property}\/
if for every $\eta\in\cS$ there are only finitely many pairs
$(\eta_1,\eta_2)\in{\cS}\times{\cS}$ such that $\eta=\eta_1+\eta_2$.
For such semigroups, any expression of the form
\[\sum_{\substack{(\eta_1,\eta_2)\text{ s.t.}\\ \eta_1+\eta_2=\eta}}
a_{\eta_1}a_{\eta_2}\]
(for fixed $\eta$) will be finite.  Thus, infinite sums can be multiplied.
So if ${\cS}$ is  a semigroup with the finite partition property,
we define the
{\em formal semigroup ring of ${\cS}$ with coefficients in $R$}\/
to be
\[R[[q;{\cS}]]:=\{\sum_{\eta\in{\cS}}a_\eta q^\eta\},\]
with the product defined by
\[(\sum_{\eta_1\in{\cS}}a_{\eta_1} q^{\eta_1})\cdot
(\sum_{\eta_2\in{\cS}}a_{\eta_2} q^{\eta_2})=
\sum_{\eta\in{\cS}}
(\sum_{\substack{(\eta_1,\eta_2)\text{ s.t.}\\ \eta_1+\eta_2=\eta}}
a_{\eta_1}a_{\eta_2}) q^\eta.\]

The semigroup $H_2(M,\Z)$ of interest to us
is actually a {\em group}\/ with a nontrivial free abelian
part, and so does not satisfy
the finite partition property.  However, in eq.~\eqref{A:correlation}
we are only required to sum over classes which
can be realized by pseudo-holomorphic curves---these generate a smaller
semigroup.  If we are using an
integrable almost-complex
structure $J$ on $M$ for which $M$ is a K\"ahler manifold,
this smaller semigroup is the {\em integral Mori semigroup}\/
defined  (in the case $h^{2,0}=0$, for simplicity) as
\[\NE(M,\Z):=\{\eta\in H_2(M,\Z)\ |\ (\omega,\eta)\ge0\ \forall\
\omega\in\overline{\cK}_J\},\]
where $\cK_J$ is the K\"ahler cone and $\overline{\cK}_J$ is its closure.
The Mori semigroup has the finite partition property (the free part lies in
a strongly convex cone, and the torsion part is finite), so we can form the
formal semigroup ring $R[[q;\NE(M,\Z)]]$.
 Presumably, by using
the symplectic version of the K\"ahler cone, we would find a similar
property for the analogous semigroup in the almost-complex case
and could form a similar ring in that case.

There is an important variant which we will have occasion to consider.
Let $\Aut_J(M)$ be the image in $\Aut H_2(M,\Z)$ of the group of
diffeomorphisms of $M$
compatible with the almost-complex structure $J$.  This group
acts on the pseudo-holomorphic curves and so permutes the Gromov--Witten
invariants.  The values of the $A$-model correlation function are preserved
by the group action, and can be
regarded as lying in the ring of invariants
\[R[[q;\NE(M,\Z)]]^{\Aut_J(M)}.\]

As in the Fano variety case, if we choose an appropriate basis (and assume
$H_2(M,\Z)$ is
torsion-free) then we can regard the correlation function defined in
eq.~\eqref{A:correlation} as taking values in a formal power
series ring $\Q[[q_1,\dots,q_r]]$.
Note that if we set all $q_j$'s to $0$, we simply recover the topological
trilinear function $(\zeta_1\cup\zeta_2\cup\zeta_3)|_{[M]}$.
But although the formal series in eq.~\eqref{A:correlation} is expected by the
physicists to converge near $q_j=0$, no convergence properties of the
series (as we have defined it) are known at present.

There is an alternative to using the semigroup rings: we could instead
use the Novikov rings \cite{Novikov} which have played a r\^ole
elsewhere in symplectic geometry \cite{HS}.
For each K\"ahler class $\omega$, the {\em Novikov ring}\/
$\Lambda_\omega$
consists of all formal power series
\[\sum_{\eta\in H_2(M,\Z)} a_\eta q^\eta\]
such that the set
\[\{\eta\suchthat a_\eta\ne0 \ \text{and}\ (\omega,\eta) <c\}\]
is finite for all $c\in\R$.  (If it is necessary to specify the ring $R$
in which the coefficients $a_\eta$ take their values, the notation
$\Lambda(\omega,R)$ is used.)  The product of two elements of
$\Lambda_\omega$ is well-defined, and also belongs to $\Lambda_\omega$.
In the case $H_2(M,\Z)=\Z^r$, $\Lambda_\omega$ is the
ring of {\em generalized Laurent series}\/
\[\{\sum a_{\vec{k}} q^{\vec{k}}\suchthat \ \text{there are only finitely
many terms with $\omega\cdot\vec{k}<c$ for any $c\in\R$}\}.\]

\section{A new algebra structure}

The correlation functions defined in the previous lecture can be
used to describe a new algebra structure
on the cohomology of $M$, in the following way.  Let $R$ be an integral domain
(usually we use $R=\Z$ or $R=\Q$), and choose a coefficient ring $\cR$ from
among
\begin{enumerate}
\item the group ring $R[q;H_2(M,\Z)]$ (in the case of a Fano variety),
\item the
formal semigroup ring with coefficients in $R$ for the Mori semigroup
$R[[q;\NE(M,\Z)]]$ (when this is well-defined, such as in the case of a
K\"ahler manifold),
\item the subring $R[[q;\NE(M,\Z)]]^{\Aut_J(M)}$ of $\Aut_J(M)$-invariants, or
\item one of the Novikov rings $\Lambda(\omega,R)$.
\end{enumerate}
We introduce a binary
operation $\zeta_1\star\zeta_2$
on $H^*(M,\cR)$ defined by the requirement
\[((\zeta_1\star\zeta_2) \cup \zeta_3)|_{[M]}
=\langle\zeta_1\,\zeta_2\,\zeta_3\rangle.\]
(This is well-defined since the cup product is a perfect pairing.)
The class $\one:=[M]^\vee\in H^0(M)$ which is dual to the
fundamental class $[M]\in H_{2n}(M)$ has the property that the
Gromov--Witten invariants $\Phi_\eta(\one,\zeta_2,\zeta_3)$ all vanish,
hence
\[\langle\one\,\zeta_1\,\zeta_2\rangle=(\zeta_2\cup\zeta_3)|_{[M]};\]
it follows that $\one$
serves as the identity element for
the binary operation $\star$.

This interpretation of the correlation function as a binary operation
also comes from physics \cite{MooreSeiberg,topgrav}.
Let us return to the picture we had of the ``pair of pants'' surface
\iffigs
$$\vbox{\centerline{\epsfysize=2cm\epsfbox{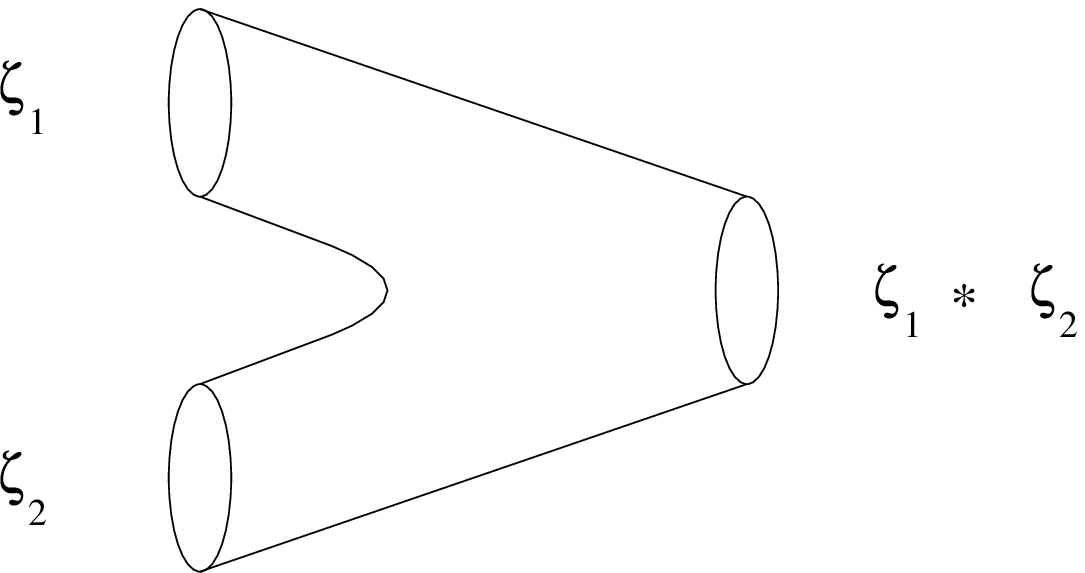}}
}$$
\else
\vglue2in\noindent
\fi
as describing a possible evolution between an initial state with two
``incoming'' vertex operators $\zeta_1$,
$\zeta_2$ on the left
and a final state with one ``outgoing'' vertex operator $\zeta_1\star\zeta_2$
 on the right.

This point of view leads
to the remarkable expectation
that the binary operation should
be associative!  A heuristic argument for this runs as follows:
the product $(\zeta_1\star\zeta_2)\star\zeta_3$
is evaluated by means of the surface
\iffigs
$$\vbox{\centerline{\epsfysize=3.5cm\epsfbox{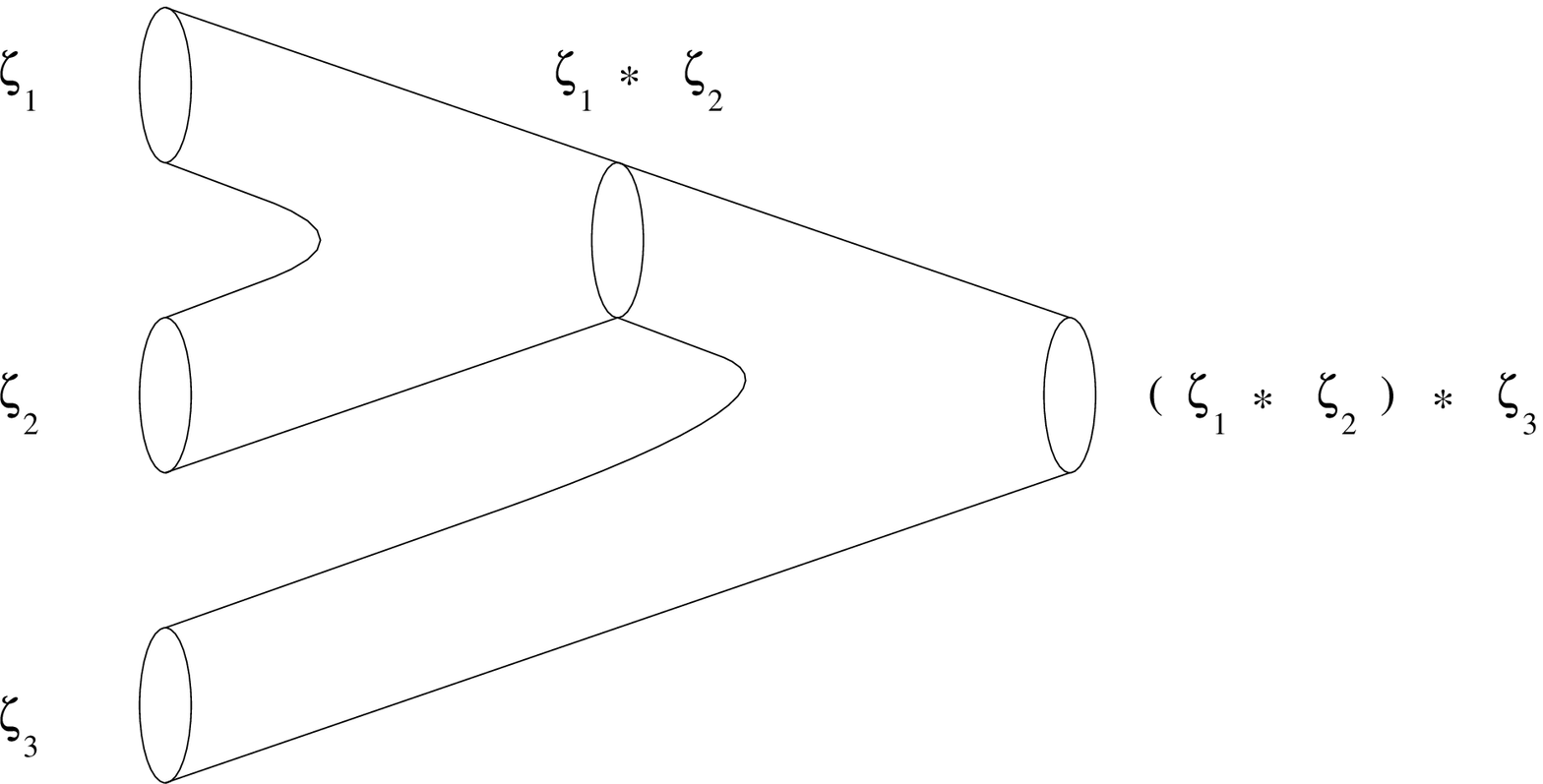}}
}$$
\else
\vglue2in\noindent
\fi
(with an outgoing vertex
operator of one piece attached to an incoming vertex operator of
the other)
while the product $\zeta_1\star(\zeta_2\star\zeta_3)$ is evaluated by
means of the surface
\iffigs
$$\vbox{\centerline{\epsfysize=3.5cm\epsfbox{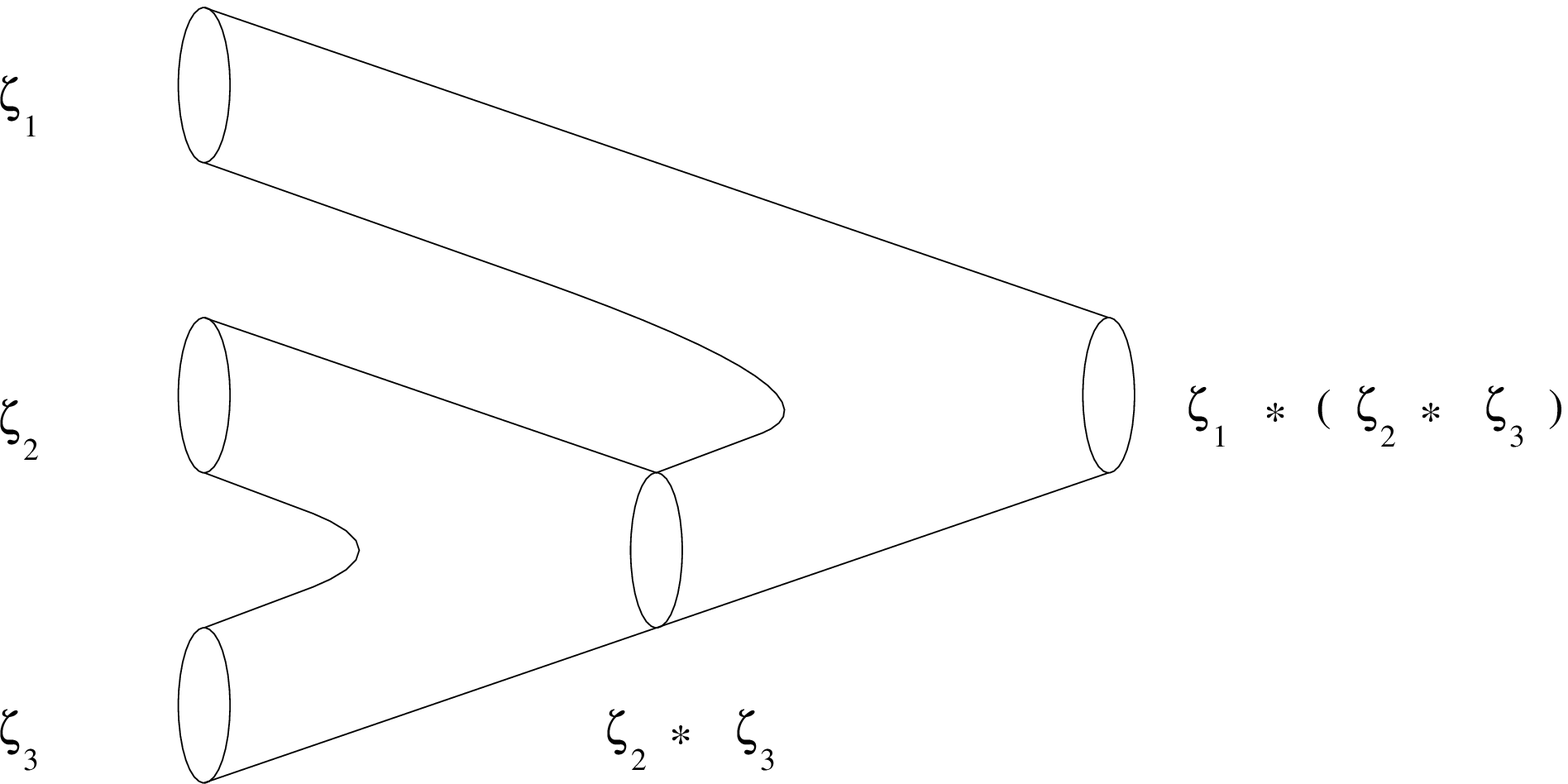}}
}$$
\else
\vglue2in\noindent
\fi
which is a deformation of the first one.  So long as the values of the
resulting quadrilinear function do not depend on the location of
the four points in $\CP^1$ used in defining it, these two products will
agree.  In fact, as we pointed out in the introduction, the correlation
functions we are studying are expected from the physics
to be precisely of this
``topological'' nature which makes them independent of the location of the
points \cite{tsm}.

This associativity property of the binary operation $\star$
can be rewritten as a set of relations
which must be satisfied among the Gromov--Witten invariants themselves.
This turns out
to be a very deep property, which had
not been proved at the time these lectures were delivered (although
proofs were given not too long thereafter \cite{RuanTian,Liu,MS}).
We have formulated the Gromov--Witten invariants and the binary operation
at this level of generality primarily because this associativity property
is such an interesting one.  However, as we will see in more detail below,
for the case of primary
interest in mirror symmetry---that of Calabi--Yau threefolds---the
associativity is automatic, and there is nothing to prove.  (Associativity
{\it does}\/ say 
something interesting for Calabi--Yau manifolds of higher dimension.)

The $\cR$-module $H^*(M,\cR)$ equipped with the binary operation $\star$
is called the {\em quantum cohomology ring}\/ of $M$, or the
{\em quantum cohomology algebra}\/ when we wish  to emphasize the
$\cR$-module structure.

We can give a more geometric description of the new binary operation,
by turning each Gromov--Witten invariant  itself into a kind of
binary operation.  Here is a heuristic description of what this
construction should look like.

We want a cohomology class $Q_\eta(\zeta_1,\zeta_2)$ with the property
that 
\[(Q_\eta(\zeta_1,\zeta_2)\cup \zeta_3)|_{[M]}=
\Phi_\eta(\zeta_1,\zeta_2,\zeta_3).\]
Consider the set of pseudo-holomorphic curves which satisfy the
conditions imposed by $\zeta_1$ and $\zeta_2$ only:
\[\MM_\eta(\zeta_1,\zeta_2):=
\{\varphi\in\MMhol\eta\suchthat\varphi(0)\in Z_1,\varphi(1)\in Z_2\},\]
where $\zeta_j=[Z_j]^\vee$.  To count the maps contributing to
$\Phi_\eta(\zeta_1,\zeta_2,\zeta_3)$, we must look for all maps in
$\MM_\eta(\zeta_1,\zeta_2)$ which also send $\infty$ into $Z_3$.  What
subset of $M$ has the property that
its intersection with $Z_3$ is in one-to-one correspondence with such
maps?  It is the subset consisting of {\em all possible}\/ points
$\varphi(\infty)$ which might be mapping to $Z_3$.  In other words,
we can write $Q_\eta(\zeta_1,\zeta_2)=[T_\eta(Z_1,Z_2)]^\vee$, where
$T_\eta(Z_1,Z_2)$ is the cycle defined by
\begin{align*}
T_\eta(Z_1,Z_2)&:=
\{P\in M\suchthat P=\varphi(\infty)
\text{ for some }\varphi\in\MM_\eta(\zeta_1,\zeta_2)\}
\\&=
\bigcup_{\varphi\in\MM_\eta(\zeta_1,\zeta_2)}\Image(\varphi)
.\end{align*}
Then $T_\eta(Z_1,Z_2)\cap Z_3$
will correspond to the maps
counted by $\Phi_\eta(\zeta_1,\zeta_2,\zeta_3)$, where
$\zeta_3=[Z_3]^\vee$.
Note that for this heuristic description to work, we need the set
$T_\eta(Z_1,Z_2)$ to be of the expected dimension.
A better formal definition of $Q_\eta(\zeta_1,\zeta_2)$ would be the
pushforward under evaluation at $\infty$ of the pullback of
$\MM_\eta(\zeta_1,\zeta_2)$ to the universal family of maps.

Expressed in these terms, then, the binary operation can be written:
\begin{equation}\label{A:binary}
\begin{split}
\zeta_1\star\zeta_2:=
\zeta_1\cup\zeta_2\ \
&+\sum_{\substack{\eta\in H_2(M,\Z),\\{-}K_M\cdot\eta>0}}
  q^{\eta}\,Q_\eta(\zeta_1,\zeta_2)\\
&+\sum_{\substack{0\ne\eta\in H_2(M,\Z),\\{-}K_M\cdot\eta=0}}
 \frac{q^{\eta}}{1-q^{\eta}}\,Q_\eta(\zeta_1,\zeta_2)
\end{split}\end{equation}

Recall that the Gromov--Witten invariant $\Phi_\eta(\zeta_1,\zeta_2,\zeta_3)$
with $\zeta_j\in H^{\ell_j}(M,\Q)$ is zero unless
\[\ell_1+\ell_2+\ell_3=2n+2({-}K_M\cdot\eta),
\quad\text{and\ \ }\ell_j\ge2.\]
It follows that if the cycle $Q_\eta(\zeta_1,\zeta_2)$ is nonzero,
we have
\[Q_\eta(\zeta_1,\zeta_2)\in
H^{2n-\ell_3}(M,\Q)=
H^{\ell_1+\ell_2-2({-}K_M\cdot\eta)}(M,\Q).\]
Thus, if $K_M=0$, then the binary operation $\star$ preserves the grading
on cohomology, while if $-K_M\cdot\eta>0$ the grading is shifted down by
$2({-}K_M\cdot\eta)$.  But note that in any case, the $\Z/2\Z$-grading
on cohomology is preserved.

Note also that $\ell_j\le2n$ implies $\ell_1+\ell_2+\ell_3\le6n$
and hence $-K_M\cdot\eta\le2n$.

\begin{exercise}  Show that
the semipositivity condition $3-n<-K_M\cdot \eta$ implies that
the grading cannot shift up, it can only shift down.
\end{exercise}

\begin{example}
(cf.\ \cite{example:pm,vafa})
We now compute an example of the quantum cohomology ring.  Let $M=\CP^n$
(with $\omega$ induced from the Fubini--Study metric, which will ensure
semipositivity).  The formal semigroup ring in this case can be written as
$\cR=\Q[[q]]$ (or we could use
$\cR=\Q[q]$ since we know the sums are finite, this being a
Fano variety).

If $C$ is any complex curve on $M$, then $-K_M\cdot C=d(n+1)$,
where $d$ is the degree of the curve.  Since $-K_M\cdot C\le2n$,
we must have  $d=1$.  So only lines (and constant maps) will contribute
to our correlation function.

Now the predicted real dimension of the space of maps $\CP^1\to M$
whose image $L$ has degree one is
\[2n+2(-K_M\cdot L)=2n+2(n+1)=4n+2\]
while the actual dimension is
\[\dim_\R\PGL(2,\C)+\dim_\R\Gr(\CP^1,\CP^n)
=6+2\cdot2(n-1)=4n+2\]
so we should be able to use the given complex structure to compute
the  invariants.

The Gromov--Witten invariants are evaluated as follows.  A basis for
$H^*(M,\Q)$ is given by the classes $\zeta^k\in H^{2k}(M,\Q)$
where $\zeta$ is the class of a hyperplane.  We choose $k_1$, $k_2$,
$k_3$, satisfying
\[2k_1+2k_2+2k_3=4n+2\]
and find that there is a {\em unique}\/ line in $\CP^m$ meeting three
fixed linear spaces of codimensions $k_1$, $k_2$ and
$k_3$.  And there is a unique map sending $0$, $1$, $\infty$ to
the intersection points with the three linear spaces.  Thus,
\[\Phi_L(\zeta^{k_1},\zeta^{k_2},\zeta^{k_3})=1.\]

Expressed in terms of the binary operation, we find that
\[\zeta^{k_1}\star\zeta^{k_2}=
\begin{cases}
\zeta^{k_1+k_2}&\text{if }k_1+k_2\le n\\
\zeta^{k_1+k_2-n-1}\,q&\text{if }k_1+k_2\ge n+1\\
\end{cases}.\]
It follows that the quantum cohomology ring can be described as:
\[\cR[\zeta]/(\zeta^{\star (n+1)}-q).\]
\end{example}

\begin{example}
If we consider the case relevant to string theory ($\dim_\C(M)=3$,
$K_M=0$), we find that the only products which differ from the
cup product are products $\zeta_1\star \zeta_2$, with $\zeta_1, \zeta_2\in
H^2(M)$,
and these are given by
\begin{equation*}
\zeta_1\star\zeta_2:=
\zeta_1\cup\zeta_2\ \
+ \sum_{0\ne\eta\in H_2(M,\Z)}
\left(
\zeta_1(\eta)\cdot\zeta_2(\eta)\cdot
\#(\MM^*_{(\eta,J)})\right)\,\frac{q^{\eta}}{1-q^{\eta}}\,\eta
\end{equation*}
Here, $\#(\MM^*_{(\eta,J)})$ denotes the number of curves in class
$\eta$ (counted with appropriate multiplicity).
\end{example}

\begin{remark}
Note that the associativity of the binary operation $\star$
is automatically satisfied by threefolds
with trivial canonical bundle, since only one of the products being
associated can be different from the cup product.
\end{remark}

\begin{example} \label{example43}
Let $\lambda\in H^2(M)$ be represented by $L$, a submanifold of real
codimension two.  If we define
\[\MM_{\eta}(\zeta):=
\{\varphi\in\MMhol\eta\suchthat\varphi(1)\in \zeta\},\]
and
\[\Gamma_\eta(\zeta):=[\{P\in M\suchthat P=\varphi(\infty)\text{ for some }
\varphi\in\MM_\eta(\zeta)\}]^\vee,\]
then we can expect that
\[Q_\eta(\lambda,\zeta)=\lambda(\eta)\cdot\Gamma_\eta(\zeta).\]
This is because the image of each $\varphi$ should meet
 $L$ in precisely $\lambda(\eta)$
points, any of which may be chosen as $\varphi(0)$.

The binary operation can then be written:
\begin{equation*}
\begin{split}
\lambda \star\zeta :=
\lambda \cup\zeta \ \
&+\sum_{\substack{\eta\in H_2(M,\Z),\\{-}K_M\cdot\eta>0}}
  \lambda(\eta)\,q^{\eta}\,\Gamma_\eta(\zeta)\\
&+\sum_{\substack{0\ne\eta\in H_2(M,\Z),\\{-}K_M\cdot\eta=0}}
 \lambda(\eta)\,
\frac{q^{\eta}}{1-q^{\eta}}\,\Gamma_\eta(\zeta)
\end{split}\end{equation*}
We regard $\Gamma_\eta$ as a map on cohomology, and call it the
{\em Gromov--Witten map}.
\end{example}

\section{Algebraic properties of the correlation functions}

Let $K$ be the field of fractions of our coefficient ring $\cR$; tensoring
the quantum cohomology ring with $K$  makes it into a $K$-algebra.
This quantum cohomology  algebra carries some additional
 structure which makes it into what is known as a
{\em Frobenius algebra}.\footnote{We
follow standard mathematical usage \cite{CR,Karp} and do not require
a Frobenius algebra to be commutative; our definition therefore differs
slightly from that in \cite{Dubrov}.
However, we will primarily be interested in the even part
$H^{ev}(M)$ of the cohomology of $M$, on which the quantum product will
in fact be commutative.}
By definition this is a $K$-algebra $A$ with a
multiplicative identity element $\one$, such that there exists
a linear functional
$\varepsilon:A\to K$ for which the induced bilinear
pairing $(x,y)\mapsto\varepsilon(x\star y)$ is nondegenerate.  There does
not seem to be a standard name for such a functional; we call it
an {\em expectation function}\/ (cf.~\cite{summing}).
If an expectation function exists at all,
then most linear functionals on $A$ can serve as expectation functions.
If $A$ is $\Z$-graded, we call $\varepsilon$ a {\em graded expectation
function}\/ when $\ker(\varepsilon)$ is a graded subalgebra of $A$ (and we
call $A$ a {\em graded Frobenius algebra}\/ when such a function exists). There
is much less freedom to choose graded expectation functions.

The cohomology of a compact manifold $M$ has the structure of a graded
Frobenius algebra, with multiplication given by cup product, $\one$ given
by the standard generator of $H^0(M)$, and a graded
expectation function given by ``evaluation on the fundamental class.''
The quantum cohomology algebra is a deformation of this algebra, with
the expectation function given by
\[\varepsilon(\zeta)= \langle\zeta\,\one\,\one\rangle,\]
which again can be interpreted as evaluation on the fundamental class.
The induced bilinear pairing
\[(\zeta_1,\zeta_1)\mapsto \varepsilon(\zeta_1\star\zeta_2)
=\langle\zeta_1\,\zeta_2\,\one\rangle\]
coincides with the usual cup product pairing.  Note that the correlation
function is also determined by $\varepsilon$ and $\star$, via
\[\langle\zeta_1\,\zeta_2\,\zeta_3\rangle=
\varepsilon(\zeta_1\star\zeta_2\star\zeta_3).\]
That is, rather than specifying the correlation function first and using
it to determine the quantum product, we can simply work with the quantum
product and the expectation function.

For most symplectic manifolds, the Frobenius algebra structure on
quantum cohomology is not graded; however, in the Calabi--Yau case we
get the structure of a graded Frobenius algebra.

Generally, given any associative
$K$-algebra $A$ with multiplicative identity,
and any linear functional $\varphi$ on $A$,
the kernel of the bilinear form $(x,y)\mapsto\varphi(x*y)$
is an ideal ${\cJ_\varphi}$, and the quotient ring
$A/{\cJ_\varphi}$ is a Frobenius algebra
with expectation function induced by $\varphi$.
If $A$ is itself a Frobenius algebra with an expectation function
$\varepsilon$,
then by a theorem of Nakayama \cite{Nakayama} (see \cite{Karp} for a modern
discussion),
$\varphi$ takes the form $\varphi(x)=\varepsilon(\alpha*x)$ for some fixed
element $\alpha\in A$, and $\cJ_\varphi$ coincides with
 the annihilator of $\alpha$.

Although the correlation functions determine the ring structure, the
opposite does not hold in general---there can be many
expectation functions on a given algebra. However,
if $A$ is a graded Frobenius algebra of finite length as a $K$-module
 and all elements of $A$ have nonnegative
degree, then the graded expectation functions
on $A$ are in one-to-one correspondence with degree $0$ elements of $A$
which are not zero-divisors.  (This is because they must all be of the
form $\varphi(x)=\varepsilon(\alpha*x)$ for some $\alpha$ which is not
a zero-divisor, but every element of degree ${}>0$ must be a zero-divisor.)
In particular, in the case of the quantum cohomology algebra of a Calabi--Yau
manifold $M$ (equipped with a symplectic structure),
we have a graded Frobenius algebra of finite length in which the
degree $0$ elements are just the one-dimensional vector space $H^0(M)$.
This means that the graded expectation function
is unique up to multiplication by an element of $K$, and
that the ring structure determines the
correlation functions up to this overall factor.
(It is not hard to see in the Calabi--Yau case that the graded
expectation function is nonzero precisely
on the top degree piece $H^{2n}(M)$, where $n=\dim_{\C}M$, and that
$H^{2n}(M)$ must also be one-dimensional.)
We will see this structure again when we study the $B$-model correlation
functions in lecture six.

\chapter*{}

\lecturename{Moduli Spaces of $\sigma$-Models}
\lecture

\markboth{D. R. Morrison,
Mathematical Aspects of Mirror Symmetry}{Lecture
5. Moduli Spaces of $\sigma$-Models}

\section{Calabi--Yau manifolds and nonlinear $\sigma$-models}\label{sec:51}

In this lecture, we wish to give a more geometric interpretation to the
$A$-model correlation functions as defined by eq.~\eqref{A:correlation}.
This geometric interpretation is motivated in part by a study of the moduli
spaces of the conformal field theories associated to Calabi--Yau manifolds,
so we begin with a description of those moduli spaces.

Let $M$ be  a K\"ahler manifold with
$K_M=0$.  Underlying $M$ is a differentiable manifold $X$ of
real dimension $2n$.  We can regard $M$ as consisting of
 $X$ together with a chosen integrable almost-complex structure $J$
and a K\"ahler
metric $g_{ij}$,
such that $K_M=0$.
(The complex manifold specified by $J$ will then be denoted $X_J$.)
If $\omega$ denotes the K\"ahler form of the metric,
then by a theorem of Calabi
\cite{calabi} there is at most one Ricci-flat metric whose K\"ahler form
is cohomologous to $\omega$; by a theorem of Yau \cite{yau}
such a Ricci-flat metric always exists.  The global holonomy of such
a metric is necessarily contained in $\SU(n)$.  (The metric being K\"ahler
implies that its holonomy is contained in $\U(n)\subset SO(2n)$;
the Ricci-flatness further restricts the holonomy to
$\SU(n)$, and also implies that the canonical bundle is trivial.
See \cite{beauville} for an account of these holonomy
conditions.)

We use the term
 {\em Calabi--Yau manifold}\/ to mean a compact connected orientable
 manifold $X$ of dimension $2n$
which admits Riemannian metrics whose (global) holonomy is contained in
$\SU(n)$. You should be aware that there are some places in the
literature (including papers of mine \cite{guide}) where
``Calabi--Yau manifold'' is used in the more restrictive sense
of a Riemannian manifold with holonomy precisely $\SU(n)$.  These
alternate definitions
will often also insist that a complex structure has been
chosen on $X$.

Given a Calabi--Yau manifold $X$ (in our sense) and a metric on it
whose holonomy lies in $\SU(n)$,
 there always exist complex structures on $X$ for
which the given metric is K\"ahler.  If $h^{2,0}=0$, then there are only
a finite number of such complex structures.
(If the universal cover
is a written as a product of indecomposable pieces, one may apply
conjugation on the various factors to obtain other complex structures.)
When $h^{2,0}>0$, however, the complex structures depend on parameters.
There are some very interesting cases with $h^{2,0}>0$, including the
famous K3 surfaces, but lack of time in these lectures
forces us to assume---with regret---that $h^{2,0}=0$ henceforth.

The physical model discussed in lecture one which
considers maps from surfaces to a six-dimensional target space is
a special case of a class of physical theories
 called ``nonlinear $\sigma$-models.''  One regards these as
quantum field theories
on the surfaces themselves, with various vertex operators and correlation
functions derived from the space of maps from the surface to the
target.  The target should be a fixed Riemannian manifold, usually
assumed to be compact.

When the Riemannian metric on the target
is (a particular perturbation of) one which
has holonomy in $\SU(n)$, the resulting ``nonlinear $\sigma$-model''
is believed to be invariant under conformal transformations of
the surface.  It thus is a type of ``conformal field
theory''---an even broader class of physical models
 which have a rich literature devoted to their study
(see \cite{ginsparg} for an introduction and further references).
Conformal field theories typically depend on finitely many
parameters, and in the case of a nonlinear
$\sigma$-model those parameters have a direct geometric interpretation.
In the Lagrangian formulation of the theory, one must specify
the metric $g_{ij}$ on the target $X$
together with an auxiliary harmonic two-form
$B$ on $X$ called the
``$B$-field.''  (To simplify matters,
we take our metrics to have holonomy in $\SU(n)$,
even though the true metrics of interest in physics will be
perturbations of those; we also assume that $H_2(X,\Z)$ has no
torsion.\footnote{The correct description of the moduli space will be
slightly different if torsion is included---see section \ref{sec:53}
below.})  The data consisting of the
pair $(g_{ij},B)$ accounts for
all local parameters in the conformal field theory moduli space,
so we get at least a good local description of moduli if we specify
such a pair.
More details about these moduli spaces can be found in \cite{ICM}.

Two  pairs
$(g_{ij},B)$ and $(g_{ij}',B')$ will
determine isomorphic conformal field theories if there is a
diffeomorphism $\varphi:X\to X$ such that $\varphi^*(g_{ij}')=g_{ij}$,
and $\varphi^*(B')-B\in H^2_{\DR}(X,\Z)$. (We use the notation
$H^k_{\DR}(X,\Z)$ to denote the image of integral cohomology
in de Rham cohomology.)
This second condition arises because the appearance of $B$ in the
Lagrangian is always in the form $\int_\Sigma B$, and the Lagrangian
is exponentiated (with an appropriate factor of $2\pi i$) in
every physically observable quantity.

We call the set of all isomorphism classes of such pairs  the
{\em  semiclassical nonlinear $\sigma$-model
moduli space}, or simply the
{\em $\sigma$-model moduli space}\/  (for short).
This may differ from the actual {\em conformal field theory moduli space}\/
for three reasons.
\begin{enumerate}
\item
It may happen that the physical theory does not exist for all values
of $g_{ij}$ and $B$.  Most of the study of these theories uses
perturbative methods, valid near a limit of ``large volume'' of the
metric, but it may be that the theory breaks down when the volume
(either of $X$, or of images of holomorphic maps into $X$) becomes
too small.

\item
On the other hand, there may be a sort of analytic continuation of
the theory beyond the region where the $\sigma$-model description is
valid.  (This was shown to occur in \cite{mmm,phases}.)
It was only claimed above that the specification of $(g_{ij},B)$
gave good {\em local}\/ parameters for the moduli.

\item
Furthermore, there could be subtle isomorphisms between conformal
field theories which do not show up in the $\sigma$-model interpretation.
This is known to happen in the K3 surface case \cite{AM:K3}, for
example (which we have no time to discuss here)---mirror symmetry
provides a new identification of conformal field theories.
\end{enumerate}
We will ignore these phenomena for the present, and concentrate
on the ``$\sigma$-model moduli space'' which parameterizes pairs
$(g_{ij},B)$ modulo equivalence.

To study this moduli space using the tools of algebraic geometry,
we must choose a complex structure on $X$.  In fact, if we consider
the set of triples $(g_{ij},B,J)$ modulo equivalence, with $J$
being an integrable almost-complex structure for which the metric
$g_{ij}$ is a Ricci-flat K\"ahler metric,
then the map from the set of equivalence classes
of triples to that of pairs is a finite map.  (It is a map of degree two if the
holonomy is precisely $SU(n)$.)

On the other hand, we can map the set of triples $(g_{ij},B,J)$
to the moduli space $\Mcx$ of complex structures on $X$.  That moduli
space is quite well-behaved, both locally and globally.  The local
structure is given by the theorem of Bogomolov--Tian--Todorov
\cite{bogomolov,tian,todorov},
which says that all first-order deformations are unobstructed.
(I recommend Bob Friedman's paper
\cite{Friedman:threefolds} for a very readable
account of this theorem.)  Thus, the moduli space $\Mcx$ will be smooth,
and the tangent space at $[J]$ can be canonically identified with
$H^1(\Thol_{X_J})$.
Globally, $\Mcx$ is known to be a quasi-projective variety (if one
specifies a ``polarization'') by a theorem of Viehweg \cite{viehweg}.
We will study the moduli space $\Mcx$ in more detail (using variations of
Hodge structure) in the next section.

The fibers of the map
\begin{equation}\label{fibrebundle}
\{(g_{ij},B,J)\}/{\sim}\ \to\ \Mcx\end{equation}
(from the set of equivalence classes of triples to the moduli space)
are spaces of the form $\cD/\Gamma$, with
\begin{align*}\cD&=H^2(X,\R)+i\,\cK_J\\
\Gamma&=H^2_{\DR}(X,\Z)\rtimes
\Aut_J(X).
\end{align*}
One hopes that the map \eqref{fibrebundle}
 is some kind of fiber bundle (at least generically);
this would require that both the family of K\"ahler cones and the
family of automorphism groups are generically locally constant.
This has been shown for the K\"ahler cones in the case of complex
dimension three by Wilson \cite{wilson}.

The tangent spaces to the fibers of the map \eqref{fibrebundle} can be
canonically identified with $H^1((\Thol_{X_J})^*)$.  Mirror symmetry
predicts that $X$ should have a mirror partner $Y$, such that the
moduli spaces of conformal field theories on $X$ and $Y$ should be
isomorphic, but with a reversal of r\^oles of $H^1(\Thol_{X_J})$
and $H^1((\Thol_{X_J})^*)$.  That is, under the isomorphism between
the conformal field theory moduli spaces, the part of the tangent
space corresponding to $H^1((\Thol_{X_J})^*)$ on $X$ should map to
the part corresponding to $H^1(\Thol_{Y_{J'}})$ on $Y$, and vice versa.

In particular, the r\^oles of base and fiber in \eqref{fibrebundle}
should be reversed.  This is at first sight a rather peculiar statement,
since the base and the fiber do not look much alike:  the base $\Mcx$
is a quasi-projective variety, whereas the fiber $\cD/\Gamma$
looks much more like a Zariski open subset of a bounded domain---a typical
model for the space is $(\Delta^*)^r$, where $\Delta^*$ is
the punctured disk.

This is in fact one of the indicators that the conformal field theory
moduli space must be analytically continued beyond the realm of
$\sigma$-models, as suggested in point 2 above.  We will see further
evidence of this at the end of lecture seven.

\section{Geometric interpretation of the $A$-model correlation functions}

We turn now to a geometric interpretation of the $A$-model correlation
functions, which in the case of Calabi--Yau manifolds will turn out
to be closely related to the spaces $\cD/\Gamma$ described above.

In the previous lecture, the symbols $q^\eta$ were treated purely
formally, which allowed us to discuss some algebraic aspects of quantum
cohomology.  Now, however, we
would like to make the new product more geometric by giving
specific values to the $q^\eta$'s, thereby making the quantum cohomology
ring into a deformation
of the usual cohomology ring.  Turning algebraic parameters
into geometric data is a familiar task for algebraic geometers; however
here, we only have formal parameters.  We will describe a natural
parameter space as a formal completion of a
certain geometric space---if some day someone
proves that the series \eqref{A:correlation} and
\eqref{A:binary} are convergent power series,
then the true parameter space
will be a neighborhood (in the classical topology)
of the completion point within
the geometric space which we will construct.

Let $\cR=\Q[[q;\NE(X_J,\Z)]]$ be the formal semigroup ring of
the integral Mori semigroup.  If $\NE(X_J,\Z)$ is finitely generated, then
we can take as the geometric space  $\Spec\C[q;\NE(X_J,\Z)]$ (the
spectrum of the semigroup ring),
and as its completion the formal scheme
$\Spf {\cR}_{\C}$,
where ${\cR}_{\C}$ denotes ${\cR}\otimes_{\Q}\C$ and
$\Spf$ is the formal spectrum.  More generally, if the ring of
$\Aut_J(X)$-invariants ${\cR}^{\Aut_J(X)}$ is the formal completion of
a ring of finite type over $\Q$,
we take
our completed parameter space to be $\Spf ({\cR}_{\C}^{\Aut_J(X)})$.

In the finitely generated case,
 this geometric
space $\Spec\C[q;\NE(X_J,\Z)]$ is in a natural way an affine toric
variety, and as such admits a rather concrete description:  the geometric
points are in one-to-one
correspondence with the set of semigroup homomorphisms
$\Hom_{\text{sg}}(\NE(X_J,\Z),\C)$, where $\C$ is given the structure of
a {\em multiplicative}\/ semigroup.  Any geometric point $\xi$
in the parameter space---regarded as
a semigroup homomorphism---specifies
compatible values $\xi(q^\eta)$ for the symbols $q^\eta$.
An important open problem is to decide for which $\xi$ the
series expressions \eqref{A:correlation} for the correlation functions
converge.  If convergent, the correlation functions would become actual
$\C$-valued functions on a parameter space (as expected by the
physicists), which would be an open subset of $\Spec\C[q;\NE(X_J,\Z)]$
in the classical topology.

To make this even more concrete, consider the case in which
the Mori semigroup is freely generated by
elements $e_1$, \dots, $e_r$ which also serve as a basis of
the lattice $H_2(X,\Z)$.  In this case, we can
define $q_j:=q^{e_j}$, and write the ring $\cR$ as a formal power series
ring
${\cR}=\Q[[q_1,\dots,q_r]]$.  The geometric space $\Spec\C[q_1,\dots,q_r]$
can then be identified
as $\C^r$ with coordinates $q_1,\dots q_r$.  One natural candidate
for the open set on which the correlation functions might converge
is
\[\{(q_1,\dots,q_r)\in\C^r\suchthat 0\le|q_j|<1\}.\]

More generally, still assuming that $H_2(X,\Z)$ is torsion-free,
suppose we choose a basis
$e_1$, \dots, $e_r$ whose span as a semigroup
{\em contains}\/ $\NE(X_J,\Z)$.  Then the corresponding
formal power series ring
$\Q[[q_1,\dots,q_r]]$ contains our coefficient ring $\cR$.
If we let $\sigma$ denote the open real cone generated by
the dual basis $e^1$, \dots, $e^r$, then that formal power series
ring can be more canonically described as the formal semigroup
ring $\cR_\sigma:=\Q[[q;\check\sigma\cap H_2(X,\Z)]]$.  The same cone
$\sigma$ can be used to give a canonical description of the
 open set specified by $0<|q_j|<1$ in the form
\[(H^2(X,\R)+i\sigma)/H^2(X,\Z).\]
(To see this, write a general element of $H^2(X,\C)$ modulo $H^2(X,\Z)$
in the form
\[\frac1{2\pi i}\sum(\log q_j)e^j,\]
and note that the condition $0<|q_j|<1$ is equivalent to
$\Im(\frac1{2\pi i}\log q_j)>0$.)
The Mori semigroup $\NE(X_J,\Z)$ will be contained in the semigroup spanned by
$\{e_j\}$ precisely when the cone $\sigma$ is contained in
the K\"ahler cone of $X_J$.

We will treat such a choice of cone $\sigma$ as specifying a coordinate
chart on the geometric space we are trying to construct.  For any such
cone, we define
\[\cD_\sigma=H^2(X,\R)+i\,\sigma\subset H^2(X,\C)\]
In terms of local coordinates, as pointed out above we have
\[\cD_\sigma/H^2(X,\Z)=\{(q_1,\dots,q_r)\suchthat 0<|q_j|<1\}.\]
The open subset of our desired geometric space will be a partial
compactification of this, defined by
\[(\cD_\sigma/H^2(X,\Z))^-=\{(q_1,\dots,q_r)\suchthat 0\le|q_j|<1\}.\]
We call the origin $0\in(\cD_\sigma/H^2(X,\Z))^-$ the {\em distinguished limit
point}\/ in
this space.

It is hoped that the expressions for the $A$-model correlation functions, or
for the binary operation $\zeta_1\star\zeta_2$, will converge in a
neighborhood of the distinguished limit point
$0$ in $(\cD_\sigma/H^2(X,\Z))^-$.  The different possible
choices of $\sigma$ will correspond to operations---such as blowing
up the boundary---which change the
compactification without changing the underlying space.

Intrinsically, we can describe $\cR_\sigma\otimes\C$ as the formal
completion of the local ring of
$(\cD_\sigma/H^2(X,\Z))^-$ at its distinguished
limit point $0$.

The geometric space which is emerging from this discussion is very
closely related to the space $\cD/\Gamma$ which formed part of the
nonlinear $\sigma$-model moduli space in the case of a Calabi--Yau
manifold with $h^{2,0}=0$.
In fact, if $\cK_J$ is the K\"ahler cone of such a Calabi--Yau
manifold which can be partitioned into cones $\sigma_\alpha$ which
are spanned by various bases of $H^2(X,\Z)$, then $\cD/H^2(X,\Z)$
is the interior of the closure of
 the union of the sets $\cD_{\sigma_\alpha}/H^2(X,\Z)$.
Ideally, one could make such a partition in an $\Aut_J(X)$-equivariant way.
This would be guaranteed by the following conjecture.

\begin{ConeConjecture}
Let $X$ be a Calabi--Yau manifold
on which a complex
structure $J$ has been chosen, and suppose that
$h^{2,0}(X)=0$.
Let $\cK_J$ be the K\"ahler
cone of $X$, let $(\cK_J)_+$
be the convex hull of $\overline{\cK}_J\cap H^2(X,\Q)$,
and let $\Aut_J(X)$ be the group of holomorphic automorphisms of $X$.
Then there exists a rational polyhedral cone $\Pi\subset(\cK_J)_+$
such that $\Aut_J(X).\Pi=(\cK_J)_+$.
\end{ConeConjecture}

A nontrivial case of this conjecture---Calabi--Yau threefolds which are
fiber products of generic rational elliptic surfaces with section (as
studied by Schoen \cite{schoen})---has been checked by Grassi and the author
\cite{GM}.  There are some other pieces of supporting evidence in examples
worked out by Borcea \cite{borcea} and
Oguiso \cite{oguiso}.

When this conjecture holds, there is a partial compactification of
$\cD/\Gamma$ constructed in \cite{compact} by gluing together the spaces
$(\cD_{\sigma_\alpha}/H^2(X,\Z))^-$ for an $\Aut_J(X)$-equivariant
partitioning
of $\cK_J$, and modding out by $\Aut_J(X)$.  This produces
a ``semi-toric'' partial compactification of the type introduced by Looijenga
\cite{Looijenga}.
Because it is covered by explicit coordinate charts, this is a convenient
type of compactification for making comparisons of correlation functions.

There is also a ``minimal'' semi-tori compactification determined from the
same data, which partially
compactifies $\cD/\Gamma$ more directly, adding several
new strata but only a single stratum
of maximal codimension (the analogue of the
``distinguished limit points'').
When the cone conjecture holds, the ring of invariants
$\cR^{\Aut_J(X)}$ is the formal completion of a ring of finite type
over $\Q$, and
the completion of the local ring of the minimal semi-toric compactification
at its distinguished point $P$ coincides with
$\Spf(\cR_{\C}^{\Aut_J(X)})$.
On such a compactification, we will expect
\begin{equation}\label{eq:limA}
\lim_{Q\to P}\langle\zeta_1\,\zeta_2\,\zeta_3\rangle_Q
=(\zeta_1\cup\zeta_2\cup\zeta_3)|_{[X]}
\end{equation}
(the ``$q_j=0$ values'' in coordinate charts).
Such a point is called a ``semiclassical
limit'' in the physics literature \cite{AL}.

\section{The r\^ole of torsion in the moduli space}\label{sec:53}

Up to this point, we have not considered the effects of possible torsion in
$H_2(X,\Z)$ and in fact we have explicitly assumed at several points that
there was no torsion.  If torsion is present, we can define the formal
semigroup ring $\cR=\Q[[q;\NE(X_J,\Z)]]$ as before, and it will have a
torsion part ${\cR}_{\text{torsion}}$ whose
spectrum is a finite set of geometric points.
This can be identified with the set of connected components of
our parameter space.  It can also be seen in the following description
of the $\sigma$-model moduli space.

The complete description of the $\sigma$-model moduli space (with the
torsion included) considers
the quantity
$e^{2\pi i(B+i\omega)}$ to lie in $\Hom(H_2(X,\Z),\C^*)$.  This can be
thought of concretely as having a torsion part, together with a
free part which lies in the space
\[\Hom(H_2(X,\Z)/\text{torsion},\C^*)\cong
H^2(X,\C^*)\cong H^2(X,\C)/H^2_{\text{DR}}(X,\Z)\]
where (as in section \ref{sec:51}) $H^2_{\text{DR}}(X,\Z)$ is the
image of $H^2(X,\Z)$ in de Rham cohomology, isomorphic to
$H^2(X,\Z)/\text{torsion}$.  A representative of the free part
can be written as $B_{\text{free}}+i\omega\in H^2(X,\C)$, where $\omega$ is
the K\"ahler form and $B_{\text{free}}$ is the
real two-form which appeared in section \ref{sec:51}.
The torsion part of $e^{2\pi i(B+i\omega)}\in\Hom(H_2(X,\Z),\C^*)$
can be identified with the
torsion part of our coefficient ring ${\cR}_{\text{torsion}}$
from the algebraic interpretation.
One way to interpret this ``$B$-field with torsion included'' is to regard
it as an element of $H^2(X,\R/\Z)$.

\chapter*{}

\lecturename{Variations of Hodge Structure}
\lecture

\markboth{D. R. Morrison,
Mathematical Aspects of Mirror Symmetry}{Lecture
6. Variations of Hodge Structure}

\section{The $B$-model correlation functions}

Our goal in this lecture is to describe the $B$-model correlation functions
and how they are related to variations of Hodge structure.  We work
with Calabi--Yau manifolds on which complex
structures have been chosen.  That is, we let $W$ be a complex manifold
with $K_W=0$. The assumption of trivial canonical bundle
is needed in order to define the $B$-model correlation functions.

Let us define
\[H^{-p,q}(W):=H^q(\Lambda^p(\Thol_W)),\]
and consider all of these groups together:
\[H^{-*}(W):=\bigoplus_{p,q}H^{-p,q}(W).\]
There is a natural ring structure on $H^{-*}(W)$ which can be thought of as
a sheaf cohomology version of the
cup product pairing:
\[H^q(\Lambda^p(\Thol_W))\otimes H^{q'}(\Lambda^{p'}(\Thol_W))
\to H^{q+q'}(\Lambda^{p+p'}(\Thol_W)).\]
Note that since these are sheaf cohomology groups, this ring structure is
not ``topological'' in nature; in fact, it depends heavily on the choice of
complex structure on $W$.

Recall that in the case of the $A$-model correlation functions
on a symplectic manifold $M$, the expectation
function which determined the Frobenius algebra structure
was a very familiar object, given by evaluating a cohomology class
on the fundamental class of $M$ (which determines a canonical map
$H^{n,n}(M)\to\C$).  By contrast, the ring structure on quantum cohomology
was unusual.  In this new ``$B$-model'' case, however,
the ring structure is straightforward but the expectation function is more
elusive. To define it, we must choose a
nonvanishing global section $\Omega^{\otimes2}$ of
$(K_W)^{\otimes2}$.
This is then used in two
steps to specify the expectation function:
\[H^{-n,n}(W)=H^n(\Lambda^n(\Thol_W))
\overset{\lhk\,\Omega}{\longrightarrow}
H^n(\O_W)\cong\left(H^0(K_W)\right)^*
\overset{\otimes\Omega}{\longrightarrow}
\C,\]
where the middle isomorphism is Serre duality.

Using this expectation function and the ``sheaf cup product'' binary operation,
we define the $B$-model correlation functions (in the standard way from
the Frobenius algebra structure):
\[\langle\beta_1\,\beta_2\,\beta_3\rangle=
((\beta_1\cup\beta_2\cup\beta_3)\lhk\,\Omega)\otimes\Omega.\]
(Once again we have a definition which is inspired by the outcome of a
calculation in the physics literature
\cite{strwit}.)
This gives a map
\[H^{-p,q}(W)\times H^{-p',q'}(W)\times H^{-(n-p-p'),n-q-q'}(W)\to\C.\]
Note that as in the $A$-model case, we actually have a graded Frobenius
algebra of finite length, so the expectation function is uniquely defined
up to a scalar multiple (which can be absorbed in the choice of
$\Omega^{\otimes 2}$.)

In order to relate this correlation function to a more familiar mathematical
object, we can proceed as follows:  first use the two $\Omega$'s to
transform two of the arguments, and then use the cup product:
\[\langle\beta_1\,\beta_2\,\beta_3\rangle=
((\beta_1\lhk\,\Omega)\cup\beta_2\cup(\beta_3\lhk\,\Omega)).\]
This variant of the correlation function can be regarded as a map
\[H^{n-p,q}(W)\times H^{-p',q'}(W)\times H^{p+p',n-q-q'}(W)\to\C,\]
or, if we treat it as a modified ``binary operation,'' as a map
\[H^{n-p,q}(W)\times H^{-p',q'}(W)\to H^{n-p-p',q+q'}(W).\]
This version of the
``binary operation'' expresses the cohomology $H^*(W)$ as
a module over the ring $H^{-*}(W)$.  As we shall see, this variant
has the pleasant property that it can be directly interpreted in
terms of variations of Hodge structure and the differential of
the period map.  Of course, the original version of the correlation
function can be recovered from this, once we have specified
$\Omega^{\otimes 2}$.

\section{Variations of Hodge structure}

We now briefly review the theory of variations of Hodge structure,
in order to explain the mathematical origin of the $B$-model
correlation functions.
Variations of Hodge structure were introduced as a tool for measuring
how the complex structure on a differentiable manifold can vary.
Good general references for this are
Griffiths et al.~\cite{transcendental}, and Schmid \cite{schmid}.

There are two primary ways one can view deformations of complex structure.
In the first viewpoint, we
fix a compact
differentiable manifold $Y$, and consider various integrable
almost-complex structures $J$ on $Y$.  Then the set of such, modulo
diffeomorphism, is known to be a finite-dimensional space.

In the second viewpoint,
we consider proper holomorphic maps $\pi:\WW\to S$ with $W_s=\pi^{-1}(s)$
diffeomorphic to $Y$.  Each fiber $W_s$ has an induced structure of a complex
manifold.
If $S$ is contractable, then $\pi$ can be trivialized in the $C^\infty$
category, and we can regard $\pi$ as specifying a family of complex
structures.  One wants to represent the functor
\[S\mapsto\{\pi:\WW\to S\}/(\text{isomorphism}),\]
by maps to a moduli space which has a ``universal family.''
This is generally too much to hope for, but there are often ``coarse
moduli spaces'' whose points are in one-to-one correspondence with the
possible complex structures.
(The appendices in \cite{MumfordFogarty}
 provide good background for moduli problems in general.)

We will study complex structures on $Y$ by studying the Hodge decomposition
induced on cohomology by each choice of complex structure.
In general, if $W_s$ is a  K\"ahler manifold there is a {\em
Hodge decomposition}\/ of the cohomology:
\begin{equation}\label{eq:hodge}
H^k(W_s,\C)\cong \bigoplus_{p+q=k}H^{p,q}(W_s).
\end{equation}
Now in a family over a contractable base, the bundle of $H^k(W_s,\C)$'s
may be canonically
trivialized.  
Over more general bases $S$ (assumed to be connected),
it is convenient to consider $R^k\pi_*\C_\WW$,
which is simply the sheaf whose local sections are topologically constant
families of cohomology classes.  This sheaf has the structure of a
{\em local system}:  it can be characterized by its fiber $H^k(W_s,\C)$
at a particular point $s\in S$ together with a representation of the
fundamental group
\[\rho:\pi_1(S,s)\to\Aut(H^k(W_s,\C))\]
which specifies what happens when the locally constant sections are
followed around loops.  There is useful dictionary \cite{rsp} between
local systems and pairs $({\cH},\nabla)$ consisting of a holomorphic
vector bundle ${\cH}$ on $S$ and a flat holomorphic connection
\[\nabla:\cH\to(\Thol_S)^*\otimes\cH.\]
The way the dictionary works is this:  given a local system $\BH$,
define $\cH=\O_S\otimes\BH$, and
$\nabla(\sum \varphi_j h^j)=\sum d\varphi_j\otimes h^j$ for $\{h^j\}$ a local
basis of sections of $\BH$.  Conversely, given $(\cH,\nabla)$,
define $\Gamma(U,\BH)=\{h\in\Gamma(U,\cH)\suchthat\nabla(h)=0\}$ for every
open set $U$.

In the case of the cohomology local system $R^k\pi_*{\C}_{\WW}$, the
associated
connection $\nabla$ on $\cH^k$ is called the {\em Gauss--Manin connection}.
An explicit version of this Gauss--Manin connection goes like this:
if we choose a local basis $\alpha^1,\dots,\alpha^r$ for the space of sections
$\Gamma(U,R^k\pi_*{\C}_{\WW})$, then any
$\beta(s)\in\Gamma(U,{\cH}^k)$ can be written
$\beta(s)=\sum f_j(s)\alpha^j$ for some coefficient functions
$f_j\in\Gamma(U,{\O}_S)$.  Then
\[\nabla(\beta)=\sum df_j\otimes \alpha^j
\in\Gamma(U,(\Thol_S)^*\otimes{\cH^k}).\]
This can be given an interpretation in terms of classical
``period integrals''  as follows.  The basis
$\alpha^1,\dots,\alpha^r$ is dual to some basis
$\gamma_1,\dots,\gamma_r\in H_k(W_{s_0},{\C})$.
Then the coefficient functions are the period integrals
 $f_j(s)=\int_{\gamma_j}\beta(s)$.
(We use integration to denote
the pairing between homology and cohomology.)

The great advantage of expressing everything in terms of the Gauss--Manin
connection is that the Gauss--Manin connection can be computed algebraically,
without knowing the topological cycles in advance.

Although the sheaf $R^k\pi_*\C_{\WW}$
of cohomology groups can be locally trivialized over the
base $S$, the Hodge decomposition \eqref{eq:hodge} will vary as we vary the
complex structure.  The properties of this variation are more conveniently
expressed using the {\em Hodge filtration}:
\[F^p(W_s):=\bigoplus_{p'\ge p}H^{p',k-p'}(W_s)\subset H^k(W,\C)\]
rather than the Hodge groups $H^{p,q}(W_s)$ directly.  The 
spaces $F^p(W_s)$ in the Hodge filtration vary
holomorphically with parameters, fitting 
together to form a holomorphic subbundle
$\cF^p\subset\cH^k$.
One might also try to construct a bundle of $H^{p,q}$'s by the simple procedure
\[{\cH}^{p,q}_{C^\infty}:=\bigcup_{s\in S}
H^{p,q}(W_s)\subset{\cH^k}.\]
As the notation indicates,
this defines a $C^\infty$ bundle, but it is not in general holomorphic.
There is a holomorphic bundle ${\cH}^{p,k-p}$ defined
by the exact sequence
\begin{equation}\label{nonsplit}
0\to{\cF}^{p+1}\to{\cF}^p\to {\cH}^{p,k-p}\to0,
\end{equation}
but this exact sequence {\em has no canonical splitting}, and $\cH^k$ cannot
in general be written as
a direct sum of these holomorphic $\cH^{p,k-p}$ bundles.

The key property satisfied by the Hodge bundles is known as {\em Griffiths
transversality}:
when we differentiate with respect to parameters by using the
Gauss--Manin connection, the Hodge filtration only shifts by one,
i.e.,
\[\nabla(\cF^p)\subset(\Thol_W)^*\otimes \cF^{p-1}.\]

To study the totality of complex structures on $W$, we can map the moduli
space, or any parameter space $S$ for a family, to the classifying space
for Hodge structures.  Each Hodge structure on a fixed vector space $H^k$
determines a point in a flag variety
\[\Flags_{(f_j)}:=\{\{0\}\subset F^k\subset\cdots\subset F^0= H^k\suchthat
\dim F^j=f_j\},\]
with the $f_j$'s specifying the dimensions of the spaces making up the
filtration. 
The group $\GL(f_0,\C)$ acts transitively on such flags, and if we fix a
reference flag $F_0^\dt$, then the flag variety can be described as
$\GL(f_0,\C)/\Stab(F_0^\dt)$.  (The stabilizer $\Stab(F_0^\dt)$
is the group of block lower triangular matrices.)
There are some additional conditions which should be imposed to get a good
Hodge structure (cf.~\cite{transcendental,schmid}); these restrict us to an
open subset $\cU$ of a subvariety\footnote{We must pass to a subvariety to
restrict to the so-called {\it polarized}\/ Hodge structures---see 
\cite{transcendental} or \cite{deligne} for an explanation of this.}
of the flag variety on which a discrete group $\Gamma$
acts, and the desired classifying space for Hodge structures is $\cU/\Gamma$.
The classifying map $S\to\cU/\Gamma$ for a family is often referred to as
the {\em period map}.

The tangent space to the flag variety can be described as
\[\bigoplus_j \Hom(F^j/F^{j+1},H^k/F^j).\]
So another way of stating Griffiths transversality is to say that the
differential of the period map $S\to\cU/\Gamma$ sends $\Thol_S$ to the
subspace 
\[\bigoplus_j \Hom(F^j/F^{j+1},F^{j-1}/F^j)
=\bigoplus_j \Hom(H^{j,k-j}(W_s),H^{j-1,k+1}(W_s))\]
of the tangent space.

The differential of the map $S\to\Flags_{(f_j)}$ factors through a map
$\Thol_S\to H^1(\Thol_W)$ which describes the first-order deformations
represented by $S$ at $[W]$.  The  map which then induces the differential
is the map
\begin{equation}\label{eq:differential}
H^1(\Thol_W)\to
\bigoplus_j \Hom(H^{j,k-j}(W),H^{j-1,k+1}(W))
\end{equation}
given by sheaf cup product.  

The success of this approach to studying the moduli of complex structures
derives from the {\em local Torelli theorem}\/ for Calabi--Yau manifolds,
which states that the map \eqref{eq:differential} is injective.  This means
that at least locally, the moduli space can be accurately described by
using variations of Hodge structure.  However, that same map can now be
given a new interpretation, as a $B$-model correlation function.
That is, {\em the $B$-model correlation
function
\[H^1(\Thol_W)\times H^{j,k-j}(W)\to H^{j-1,k+1}(W)\]
coincides with the differential of the period map!}

We now restrict our attention to the middle-dimensional cohomology
$H^n(W,\C)$.
Stated in terms of the Gauss--Manin connection, we find the following
``bundle version'' of our correlation function \cite{guide}:
given a vector field $\theta$ on the moduli space
and sections  $\alpha\in\cF^j$, $\beta\in\cF^{j-1}$,
 the correlation function is
\[\langle\theta\,\alpha\,\beta\rangle=\int_W\nabla_\theta(\alpha)\wedge\beta\]
(where $\nabla_\theta=\theta\lhk\nabla$
denotes the directional derivative in direction
$\theta$).

However, as used in physics the correlation function is a specific
function rather than a map between bundles.  To find this interpretation,
we will need to choose specific sections of these bundles on which
to evaluate the map.  It is this issue to which we now turn.

\section{Splitting the Hodge filtration}

Our method for specifying sections of the Hodge bundles will be given
in terms of a choice of splitting for the Hodge filtration on the
middle-dimensional cohomology $H^n(W,\C)$, i.e., a set of
splittings of the exact sequences
(\ref{nonsplit}) (but defined only locally in the parameter space).
We determine such a splitting by means of a filtration on {\em homology},
which we think of as specifying ``which periods to calculate.''

Let
${\BS}_\dt$
be a filtration of the homology local system
$\Hom(R^n\pi_*{\C}_{\WW},{\C}_S)$ by sub-local systems, and let
\[{\BS}^\ell:=\Ann({\BS}_{\ell-1}):=
\{\alpha\in \cH^n\suchthat 
\int_\gamma\alpha=0\ \forall\ \gamma\in{\BS}_{\ell-1}\}.\]
be the associated filtration of annihilators of $\BS_\dt$ in cohomology.
We say that $\BS_\dt$ is a {\em splitting filtration for $\cF^\dt$}\/ if
$({\cH}^n)_s \cong ({\cF}^p)_s\oplus({\BS}^{n-p+1})_s$ for every $s\in S$ and
for every $0\le p\le n$.  (In this case,
$\BS^\dt$ and $\cF^\dt$ are called
{\em opposite filtrations of weight $n$}\/ \cite{deligne}.)

One way of producing examples of splitting filtrations is as follows:
fix a point $s\in S$, and consider the conjugate of the Hodge
filtration at $s_0$, namely, $\overline{F^q}_{s_0}$.
The ``opposite'' property for these filtrations is easy to check:
by definition
\begin{align*}({\cF}^p)_s&=H^{n,0}(W_s)\oplus\cdots\oplus H^{p,n-p}(W_s)\\
\intertext{and so}
(\overline{{\cF}^{n-p+1}})_s&=
\overline{H^{n,0}(W_s)\oplus\cdots\oplus H^{n-p+1,p-1}(W_s)}\\
&=H^{0,n}(W_s)\oplus\cdots\oplus H^{p-1,n-p+1}(W_s),
\end{align*}
where we have used the fact that $\overline{H^{p,q}(W_s)}=H^{q,p}(W_s)$.
The Gauss--Manin connection can be used to extend this from a filtration at
one point to a filtration of the local system.  Although this filtration
only coincides with the
conjugate of the Hodge filtration at one point in the parameter space, it
remains opposite to the Hodge filtration at all points nearby.

Given a splitting filtration ${\BS}_\dt$, we define
\[{\cH}^{p,q}_{\BS}:={\cF}^p\cap\Ann({\BS}_{q-1}),\]
on any open set on which $\BS_\dt$ is single-valued.  Then
\[{\cH}=\bigoplus_{p=0}^n {\cH}^{p,q}_{\BS} \quad \text{and} \quad
\cF^p=\bigoplus_{p'\ge p}\cH^{p',n-p'}_\BS.\]
(This is the promised splitting of the Hodge filtration.)
More concretely, this space can be described in terms of
 conditions on the periods as follows.
The sections of ${\cH}^{p,q}_{\BS}$ over $U$ are
\[\Gamma(U,\cH^{p,q}_{\BS}):= \{\beta\in\Gamma(U,{\cF}^p)\ |\
\int_\gamma\beta=0\ \forall\ \gamma\in{\BS}_{q-1}\}.\]
We also define a space of {\em distinguished sections}\/
of ${\cH}^{p,q}_{\BS}$ by
\[\Gamma(U,{\cH}^{p,q}_{\BS})_{\text{dist}}:=
\{\beta\in\Gamma(U,{\cH}^{p,q}_{\BS})\ |\
d\left(\int_\gamma\beta\right) =0\ \forall\ \gamma\in{\BS}_{q}\}.\]
(That is, the period integrals $\int_\gamma\beta$ are constant for
all $\gamma\in{\BS}_{q}$, and vanish for all $\gamma\in\BS_{q-1}$.)

For each ${\BS}_\dt$, then, we can define specific
$B$-model correlation functions,
using the $\Omega$ coming from the distinguished section of
${\cH}^{n,0}_{\BS}$  (which is well defined up to
a complex scalar multiple).
This has the advantage that the correlation functions have been
turned into actual functions on a parameter space (in accord with
the physicists' interpretation) rather than sections of a bundle.
The disadvantage is that further parameters---in the form of a choice
of splitting---have been introduced.  However, the necessity of considering
further parameters such as these, on which the correlation functions
will depend anti-holomorphically rather than holomorphically, was
recently realized in the physics literature \cite{t:tstar}.

In addition to the distinguished $n$-form $\Omega$, our choice of
splitting determines  a family of
distinguished vector fields which when contracted with $\Omega$ yield the
distinguished sections of $\cH^{n-1,1}_{\BS}$.
These vector fields
can be integrated into {\em canonical coordinates}, well-defined
up to a $\GL(r,{\C})$ transformation.  (The flexibility of that final
$\GL(r,\C)$
choice comes from the constants of integration, which must also be
specified in order to completely determine a set of canonical coordinates.)

A bit more explicitly, if $\gamma_0$ spans ${\BS}_0$
and $\gamma_0, \gamma_1,\dots,\gamma_r$ span ${\BS}_1$, then the distinguished
$\Omega$ satisfies $\int_{\gamma_0}\Omega=\text{constant}$,
and the coordinates are given by
\[\int_{\gamma_1}\Omega,\dots,\int_{\gamma_r}\Omega.\]
If we start with an arbitrary $n$-form $\widetilde\Omega$, we can write the
distinguished $n$-form as
\[\Omega:=\frac{\widetilde\Omega}{\int_{\gamma_0}\widetilde\Omega}\]
and the canonical coordinates as
\[\frac{\int_{\gamma_1}\widetilde\Omega}{\int_{\gamma_0}\widetilde\Omega},
\dots,
\frac{\int_{\gamma_r}\widetilde\Omega}{\int_{\gamma_0}\widetilde\Omega}.\]
This is the most general possible form for canonical coordinates (and
a distinguished $n$-form) needed for the physical theory, according
to recent work in physics
\cite{BCOV:KS}.

Let us fix a splitting filtration ${\BS_\dt}$.  Consider a basis
$\{\beta^i\}$ of
${\cH^n}$ consisting of distinguished sections of the bundles
${\cH}_{\BS}^{p,q}$ (ordered so that the basis is also
 adapted to the Hodge filtration ${\cF}^{\dt}$),
and a multi-valued basis $\{\gamma_j\}$ of the homology local system
$\Hom(R^n\pi_*{\C}_{\WW},{\C}_S)$, adapted to the splitting
filtration ${\BS_\dt}$.  Then the period matrix $(\int_{\gamma_j}\beta^i)$
(which has multi-valued entries) will take a block upper triangular
form with constant diagonal blocks.  And if we calculate the connection
matrix in the basis $\{\beta^i\}$, it takes the special form
\[\begin{pmatrix}
0&A^1_0&0&&\cdots&0\\
&0&A^1_1&0&\cdots&0\\
&&\ddots&\ddots&&\vdots\\
\vdots&\vdots&&&0&A^1_{n-1}\\
0&0&\cdots&&&0
\end{pmatrix}\]
in which the only nonzero entries are in the first block superdiagonal
of the matrix.  The entries $A^1_j$ precisely contain the data for the
$B$-model correlation functions,  calculated in our distinguished
basis.

\chapter*{}

\lecturename{The $A$-Variation of Hodge Structure}
\lecture

\markboth{D. R. Morrison, Mathematical Aspects of Mirror Symmetry}{Lecture
7. The $A$-Variation of Hodge Structure}

\section{Variations of Hodge structure near the boundary of moduli}

In this lecture, we begin by reviewing the asymptotic behavior of a
variation of Hodge structure near the boundary of moduli space, and the
behavior of the $B$-model correlation functions there.  Comparing to the
$A$-model correlation functions will
reveal some similarities---this is one of the
hints of mirror symmetry.  We make the similarities even more apparent by
using the $A$-model correlation functions to construct a new variation of
Hodge structure, which we call the $A$-variation of Hodge structure.

Let $S=(\Delta^*)^r\subset \overline{S}=\Delta^r$, and suppose we are given
a family $\pi:\WW\to S$ of complex manifolds.  We will assume
that
there is a way to complete this to a family
$\bar\pi:\overline{\WW}\to\overline{S}$ in which $\bar\pi$ is still
proper (but no longer smooth).  
Thus, $0\in\overline{S}$ is a boundary point in the parameter space.
Pick a basepoint $s\in S$; then
the fundamental group $\pi_1(S,s)$ is generated by loops
$\gamma^{(1)}$,\dots,$\gamma^{(r)}$ with $\gamma^{(j)}$
homotopic to the standard
generator of $\pi_1(\Delta^*_j)$, where $\Delta^*_j$ is the $j^{\text{th}}$
factor in $(\Delta^*)^r$.

\begin{MonodromyTheorem}[Landman \cite{monodromy}]
The action of each generator $\gamma^{(j)}$ gives a quasi-unipotent
automorphism $T^{(j)}$ of $H^k(W_s,\Q)$, i.e.,
$(((T^{(j)})^{b_j}-I)^{r_j}=0$.  (This is called {\em unipotent}\/ if $b_j=1$.)
\end{MonodromyTheorem}

We will restrict attention to the unipotent case.  This is partially
for technical convenience, but in fact, in the examples which have
been calculated for mirror symmetry purposes, only unipotent monodromy
transformations have played a r\^ole.

When $T^{(j)}$ is unipotent, its logarithm can be defined by the following
sum (which is finite).
\[N^{(j)}:=\log T^{(j)}:= (T^{(j)}-I) - \frac12\,(T^{(j)}-I)^2+\cdots.\]
(Note that the $T^{(j)}$'s and $N^{(j)}$'s all commute.)
Let $z_1,\dots,z_r$ be coordinates on $\overline{S}$, with $z_j$
a coordinate on the $j^{\text{th}}$ disk.
Consider the operator
\begin{multline*}
\cN:=\exp\left(-\frac1{2\pi i}\sum\,\log z_j\,N^{(j)}\right)=\\
I + \left(-\frac1{2\pi i}\sum\,
\log z_j\,N^{(j)}\right)
+ \frac1{2!}\left(-\frac1{2\pi i}\sum\,\log z_j\,N^{(j)}\right)^2+\cdots
\end{multline*}
(also a finite sum).  For any section $e$ of the local system
$R^k\pi_*(\C_\WW)$, a simple calculation shows that
\begin{equation}\label{eq:GMext}
\nabla(\cN (e)) = -\frac1{2\pi i}\sum\frac{dz_j}{z_j}\,N^{(j)}(e).
\end{equation}
The key facts about the asymptotic behavior are as follows.

\begin{NilpotentOrbitTheorem}[Schmid \cite{schmid}]
Assume that each monodromy transformation $T^{(j)}$  is unipotent.
Let $e_1(s),\dots,e_r(s)$ be a multi-valued basis of $R^k\pi_*(\C_\WW)$, and
let $\eta_\ell:=\cN(e_\ell)$.  Then each $\eta_\ell$
is a single-valued section of $\cH^k$ on $S$, and together they can be used to
generate an extension $\overline{\cH}^k$ of $\cH^k$ to $\overline{S}$.
By eq.~\eqref{eq:GMext}, the Gauss--Manin connection
extends to a connection on $\overline{\cH}^k$ (again
denoted by $\nabla$) with
{\em regular singular points}, i.e., the extended connection is a map
\[\nabla:\overline{\cH}^k\to
(\Thol_{\overline{S}})^*(\log B)\otimes \overline{\cH}^k\]
where $(\Thol_{\overline{S}})^*(\log B)$ is the free
$\O_{\overline{S}}$-module generated by $\frac{dz_j}{z_j}$, $j=1,\dots,r$.
Moreover, the Hodge bundles $\cF^p$
have locally free extensions to subbundles
$\overline{\cF}^p\subset\overline{\cH}^k$ such that
\[\nabla(\overline{\cF}^p)\subset
(\Thol_{\overline{S}})^*(\log B)\otimes \overline{\cF}^{p-1}.\]
\end{NilpotentOrbitTheorem}

The asymptotic behavior as $z_j\to0$ of the $B$-model correlation functions
\[\langle\theta\,\alpha\,\beta\rangle=\int_W\nabla_\theta(\alpha)\wedge\beta\]
can be deduced from this theorem.
If we let $\theta_j=2\pi i\,z_j\,\frac{d}{dz_j}$
(chosen to remove poles in the asymptotic expression for the
correlation function) then the leading term in
$\langle\theta_j\,\eta_\ell\,\beta\rangle$ is given by the monodromy:
\begin{equation}\label{eq:limB}
\lim_{z_j\to0}\langle\theta_j\,\eta_\ell\,\beta\rangle
=-\int_W N^{(j)}(e_\ell)\wedge\beta.
\end{equation}

The essential properties of the monodromy are captured by the
{\em monodromy weight filtration}\/ ${\BW}_{\dt}$
on the cohomology, which has the properties that
$N^{(j)}{\BW}_\ell\subset {\BW}_{\ell-2}$, and that for any positive
real numbers $a_1$, \dots, $a_r$, the operator $N:=\sum a_j N^{(j)}$
induces isomorphisms $N^\ell:\Gr^{\BW}_{n+\ell}\to\Gr^{\BW}_{n-\ell}$.
Any splitting filtration which we use to make calculations of $B$-model
correlation functions must be somehow compatible with this monodromy weight
filtration, if those calculations are to make sense near the boundary.

If mirror symmetry is going to hold, there must be a correspondence between
the limiting behaviors described in eqs.~\eqref{eq:limA} and \eqref{eq:limB}.
In fact, the first thing to notice is that the natural flat coordinates on
the $A$-model moduli space are multiple-valued, with the ambiguity
precisely specified by $H_{\DR}^2(M,\Z)$.  So there must be some part of
the monodromy weight filtration which matches that behavior.  This
motivated the following definition, first given in \cite{guide,compact}
(cf.~also \cite{deligne}).

We say that a boundary point is {\em maximally unipotent}\/ if
\[{\cH}_s=({\cF}^n)_s\oplus({\BW}_{2n-2})_s\]
and
\[{\cH}_s=({\cF}^{n-1})_s\oplus({\BW}_{2n-4})_s\]
for all $s$ near the point.  With this definition, the distinguished
holomorphic $n$-form and the canonical coordinates can be
defined as in lecture six.

There is an alternate version of this ``maximally unipotent monodromy''
condition, which agrees with the original one for
Calabi--Yau threefolds, but is more restrictive
in higher dimension.
We say that a boundary point is {\em strongly
maximally unipotent}\/ if the weight filtration ${\BW}_{\dt}$
has nontrivial graded pieces in even degree only, and if
the induced filtration on homology defined by
\[\BS_\ell:=\Ann(\BW_{2n-2\ell+2})\]
is a splitting filtration.  (Note that the corresponding filtration on
cohomology is then
\[\BS^\ell:=\Ann(\BS_{\ell-1})=\BW_{2n-2\ell};\]
this is the filtration which should be opposite to the Hodge filtration.)
In this case, we will be able to use distinguished sections to
calculate $B$-model correlation functions, as explained earlier.

At the moment, only the original version of the definition
has  been justified to the
satisfaction of physicists as an appropriate characterization of points
which should be useful for
mirror symmetry.  To completely carry out a mirror symmetry type
calculation, though, the second version would seem to be necessary.
And as we shall see, that version has been extremely successful
in examples.

Actually, even just at the level of the monodromy action, the parallels
between the structure of the Lefschetz operators on the cohomology and the
action of monodromy are rather striking, as was first observed by Cattani,
Kaplan and Schmid 
\cite{CKS}.
The
operators $\ad(e^j)$ describe Lefschetz decompositions of the cohomology
of $M$, which have many structural parallels 
to the monodromy weight filtration
at a maximally unipotent point.

\section{Reinterpreting the $A$-model correlation functions}

Let $M$ be a Calabi--Yau manifold on which a complex structure and K\"ahler
metric have been fixed.
Inspired by some of the similarities between the two different types of
correlation functions, we wish to improve the analogy by translating the
$A$-model correlation functions into data describing a variation of Hodge
structure.
Consider the moduli space $\cD/\Gamma$ for $A$-model correlation functions,
and a coordinate chart specified by a cone $\sigma$:
\[\begin{array}{ccccc}
\cD/\Gamma&\leftarrow&\cD_\sigma/H^2(M,\Z)&\cong&(\Delta^*)^r\\
&&\cap\raise1pt\hbox{$\scriptstyle|$}&&\cap\raise1pt\hbox{$\scriptstyle|$}\\
&&(\cD_\sigma/H^2(M,\Z))^-&\cong&\Delta^r
\end{array}\]
We assume that the cone $\sigma$ (which we call a {\em framing}\/) is generated
by a basis $e^1$, \dots, $e^r$ of $H^2(M,\Z)$.  Let $t_1$, \dots, $t_r$
be coordinates on $H^2(M,\C)$ dual to this basis (so that elements of
$H^2(M,\C)$ take the form $\sum t_j e^j$).  The natural vector fields
for making calculations of correlation functions which involve
a term from the tangent space $H^2(M,\C)$ are the vector fields
$\partial/\partial t_j$.  These are the analogues of the distinguished
vector fields which we had on the $B$-model side.

On the other hand, natural coordinates on $\cD_\sigma/H^2(M,\Z)$
are furnished by $q_j=\exp(2\pi i\,t_j)$.  Then
\[\frac{\partial}{\partial t_j}=2\pi i\,q_j\,\frac{\partial}{\partial q_j},\]
from which we conclude that those correlation functions should
naturally be evaluated on the basis $2\pi i\,q_j\,\partial/\partial q_j$ of
the sheaf of logarithmic vector fields on
the space $(\cD_\sigma/H^2(M,\Z))^-$.

We identify $\partial/\partial t_j$ with the operation of taking the
quantum product with the basis element $e^j\in H^2(M,\Q)$.  
The resulting map is
determined by the correlation functions of the form
$\langle e^j\, \alpha\,\beta\rangle$.
We had a particularly simple form for
these correlation functions, given in example \ref{example43}, in terms of the
Gromov--Witten maps $\Gamma_\eta$.  We now wish
to reinterpret that formula in the following way.  We will describe
a holomorphic bundle\footnote{There are a few variants to this
construction, in which one uses slightly different bundles.
Essentially, one can restrict to any subbundle of $\bigoplus H^{\ell,\ell}(M)$
which is preserved by cup products with the part of $H^{1,1}(M)$ which it
contains.}
$\cE:=\left(\bigoplus H^{\ell,\ell}(M)\right)\otimes
\O_{(\cD_\sigma/H^2(M,\Z))^-}$
with a
 connection\footnote{I am indebted to P. Deligne for advice
which led to this form of the formula (cf.~\cite{deligne}).}
(with regular singular points)
\[\nabla:=\frac1{2\pi i}\,\left(
\sum  \dlog q_j\otimes\ad(e^j)+
\sum_{0\ne\eta\in H_2(M,\Z)}
\dlog\left(\frac1{1-q^{\eta}}\right)\otimes\Gamma_\eta
\right)\]
which was
derived from the formulas for $e^j{\star}$, where $\ad(e^j):H^k(M)\to
H^{k+2}(M)$ is defined by $\ad(e^j)(A)=e^j\cup A$.  We also define a ``Hodge
filtration''
\[\cE^p:=\left(\bigoplus_{0\le\ell\le m-p}
H^{\ell,\ell}(M)\right)\otimes\O_{(\cD_\sigma/H^2(M,\Z))^-}.\]
This describes a structure we call the {\em framed $A$-variation of
Hodge structure with framing $\sigma$}.  To be a bit more precise,
we should study ``formally degenerating variations of Hodge structure,''
since the series used to define $\nabla$ is only formal.
(We won't formulate that theory in detail here.)

The connection $\nabla$ which we defined from the Gromov--Witten
invariants is in fact a {\em flat}\/ holomorphic connection \cite{topgrav}.
The
flatness follows from the associativity (and commutativity)
of the binary operation.
In fact, since the directional derivatives with respect to $\nabla$
corresponded to binary products
$e^j\star\zeta$ (where $e^j$ describes the direction
of the derivative), iterated directional derivatives have the form
$e^k\star(e^j\star\zeta)$.  We would simply need to know
that reversing the order of $j$ and $k$ produces the same result,
and this is guaranteed by the commutativity and associativity.

In particular, the flatness is automatic when $\dim_\C M=3$, a case in
which 
there is no issue of associativity.  The recent proofs of associativity of
quantum cohomology \cite{RuanTian,Liu,MS} guarantee that this connection is
flat in arbitrary dimension.

As in the geometric case, there is an additional structure associated to
this variation of Hodge structure: a local system.  The local system on
homology takes the simple form
\[\BS_\ell:=H_{0,0}\oplus H_{1,1}\oplus\dots\oplus H_{\ell,\ell},\]
and the corresponding local
system on cohomology then becomes
\[\BS^\ell=H^{\ell,\ell}\oplus H^{\ell+1,\ell+1}\oplus\dots\oplus H^{n,n}.\]
The logarithms of the
monodromy actions which define these local systems are specified
by the topological pairings, and coincide with the cup-product maps
\[H^2(M,\Z)\otimes \BS^\ell\to \BS^{\ell+1}.\]

In the next lecture, we will formulate a precise conjecture which equates
this $A$-variation of Hodge structure with the geometric variation of Hodge
structure on a mirror partner.

\section{Beyond the K\"ahler cone}\label{sec73}

We indicated in lecture five that the conformal field theory moduli space
is actually {\em larger}\/ than the nonlinear $\sigma$-model moduli space.
We can now explain how this comes about---it is due to an analysis of the
effect of flops on the conformal field theory.

\begin{figure}
\iffigs
$$\vbox{\centerline{\epsfysize=3cm\epsfbox{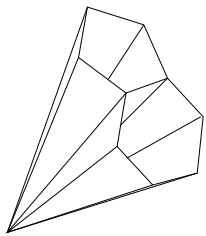}}
}$$
\else
\vglue2in\noindent
\fi
\caption{Adjacent K\"ahler cones} \label{fig0}
\end{figure}

Flops are birational transformations among Calabi--Yau threefolds which
have been studied extensively as part of the minimal model program
(see for example \cite{CKM}).
The effect of flops on
the K\"ahler cone of a Calabi--Yau threefold is as follows.  Given a
Calabi--Yau threefold $X$ with a complex structure $J$,
and a linear system $|L|$ inducing a flopping contraction
from $X_J$ to $\widehat X_{\widehat J}$, the K\"ahler cones $\cK_J$ and
$\widehat{\cK}_{\widehat J}$
share a common wall, which contains the class of $|L|$, as depicted
in figure \ref{fig0}.
The K\"ahler cone has already occurred in our discussion
of the moduli spaces of $\sigma$-models.  The natural question arises:
suppose we attempt to ``attach'' the moduli spaces $\cD/\Gamma$
and $\widehat{\cD}/\widehat{\Gamma}$ along (the images of) their common wall?
In fact, it now appears likely
that the conformal field theory moduli spaces of
$X$ and $\widehat X$ are analytic continuations of each other, and that this
``attached space'' is a part of the full conformal field theory moduli
space \cite{mmm,phases}.  
(This at least seems to happen in examples---the arguments
for this rely on mirror symmetry, and involve finding regions in the mirror's
moduli space  which correspond to the $X_J$ and
$\widehat X_{\widehat J}$ theories, respectively.)  One of the consequences
of this would be an analytic continuation of correlation functions
from $\cD/\Gamma$
to $\widehat{\cD}/\widehat{\Gamma}$.

Here is a formal calculation from \cite{phases,small}
which supports this analytic continuation idea (see also
\cite{beyond} for a more mathematical treatment).
The union of all of the K\"ahler cones of birational models of $X_J$ is
known as the {\em movable cone}\/ $\Mov{X_J}$ \cite{kawamata}.
We compute in the
formal semigroup ring
$\Q[q;\Mov(X_J)^\vee]$ (which we identify canonically with
the same ring for $\widehat X_{\widehat J}$),
and so the computation is purely formal.

Consider the simplest
flop: the flop based on a collection of disjoint holomorphic rational
curves $\Gamma_i\subset X_J$ (in a common homology class $[\Gamma]$) such that
the normal bundle is $N_{\Gamma_i/X_J}=\O(-1)\oplus \O(-1)$.
(These curves must be flopped
simultaneously in order to ensure
that the flopped variety is K\"ahler.)
A reasonable genericity assumption about the {\em other}\/ rational
curves on $X_J$ is this:
all (pseudo-)holomorphic curves in classes $\eta\not\in\R_{>0}[\Gamma]$
are disjoint from the $\Gamma_i$'s.
Since there is a proper transform map on divisors,
the Gromov--Witten invariants (which in this case are determined entirely
by intersection properties of $\eta$ and the number of elements in
$\MM^*_{(\eta,J)}$) do not change when passing
from $X_J$ to $\widehat X_{\widehat J}$,
except for the invariants $\Phi_{[\Gamma]}$ themselves.  The cup product
can also change.

The $A$-model correlation functions on $X_J$ can be written in the form
\begin{align*}
\langle A\,B\,C\rangle=A\cdot B\cdot C
&+\frac{q^{[\Gamma]}}{1-q^{[\Gamma]}}\,(A\cdot\Gamma)(B\cdot\Gamma)
(C\cdot\Gamma)\,n_\Gamma \\
&+\sum_{\substack{\eta\in H_2(X,\Z)\\ \eta\ne\lambda\Gamma}}
\frac{q^\eta}{1-q^\eta}\,\Phi_\eta(A,B,C).
\end{align*}
Only the first terms change when passing to $\widehat X_{\widehat J}$
 and in fact
we claim that
\begin{align*}
A\cdot B\cdot C
&+\frac{q^{[\Gamma]}}{1-q^{[\Gamma]}}\,(A\cdot\Gamma)(B\cdot\Gamma)
(C\cdot\Gamma)\,n_\Gamma \\
= & {\widehat A}\cdot {\widehat B}\cdot {\widehat C}
+\frac{q^{[\widehat \Gamma]}}{1-q^{[\widehat \Gamma]}}\,
(\widehat A\cdot\widehat \Gamma)(\widehat B\cdot\widehat \Gamma)
(\widehat C\cdot\widehat \Gamma)\,n_{\widehat \Gamma},
\end{align*}
where $\widehat A$, $\widehat B$, and $\widehat C$ are the proper transforms
of $A$, $B$, and $C$.
(In other words, the change in the topological term is precisely
compensated for by the change in the $q^{[\Gamma]}$ term.)

We will check this formula in the case
in which $A$ and $B$ meet one of the curves $\Gamma$ transversally
at $a$ and $b$ points, respectively,
and $\widehat C$ meets $\widehat\Gamma$ transversally at $c$ points.
(The general case can
be deduced from this one.)
Then $C$ must
contain $\Gamma$ with multiplicity $c$,
and the configuration of divisors is as in figure \ref{fig1}
(which illustrates the case $a=b=c=1$ for simplicity).
$A$ and $B$ have no intersection points along $\Gamma$, but both
$\widehat A$ and $\widehat B$ contain $\widehat \Gamma$, and they
meet $\widehat C$.  The total number of intersection points
of $\widehat A$, $\widehat B$ and $\widehat C$ (counted
with multiplicity) which lie in $\widehat \Gamma$ is thus $abc$.

\refstepcounter{figure}\label{fig1}
\begin{figure}
\iffigs
$$
\matrix\epsfxsize=2in\epsfbox{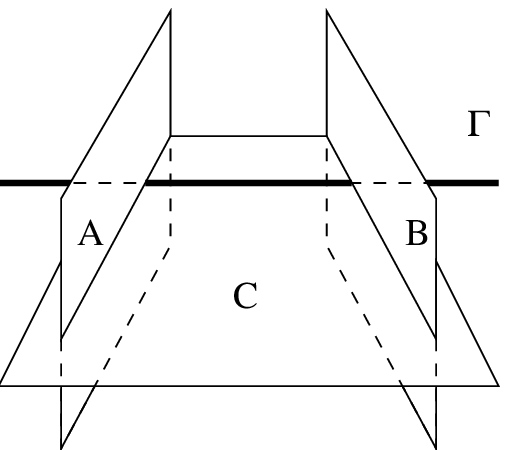} & \qquad &
\epsfxsize=2in\epsfbox{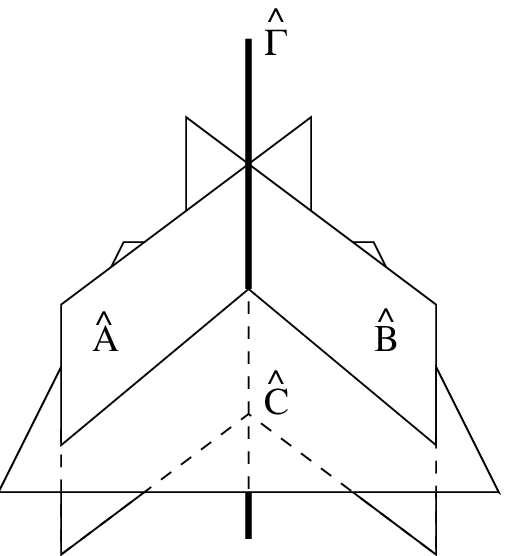} \cr
\quad & & \cr
\hbox{{\footnotesize\bfseries Figure
\ref{fig1}a}{.\footnotesize\mdseries\upshape\enspace Before the flop.}} & &
\hbox{{\footnotesize\bfseries Figure
\ref{fig1}b}{.\footnotesize\mdseries\upshape\enspace After the flop.}}
\cr\endmatrix
$$
\else
\vglue3in\noindent
\fi
\end{figure}

Since a similar thing happens for each curve $\Gamma_i$
in the numerical equivalence
class, we see that
\begin{equation}
\widehat A\cdot\widehat B\cdot\widehat C
-A\cdot B\cdot C = abc\,n_\Gamma
=-(A\cdot\Gamma)(B\cdot\Gamma)(C\cdot\Gamma)\,n_\Gamma
\label{eq:three}
\end{equation}
(using $A\cdot\Gamma=a$, $B\cdot\Gamma=b$, $C\cdot\Gamma=-c$).
On the other hand, since $[\widehat\Gamma]=-[\Gamma]$ and
$n_{\widehat\Gamma}=n_\Gamma$, we can compute:
\begin{equation}\label{eq:four}\begin{split}
\frac{q^{[\Gamma]}}{1-q^{[\Gamma]}}\,&(A\cdot\Gamma)(B\cdot\Gamma)
(C\cdot\Gamma)\,n_\Gamma
- \frac{q^{[\widehat \Gamma]}}{1-q^{[\widehat \Gamma]}}\,
(\widehat A\cdot\widehat \Gamma)(\widehat B\cdot\widehat \Gamma)
(\widehat C\cdot\widehat \Gamma)\,n_{\widehat \Gamma}\\
=&\frac{q^{[\Gamma]}}{1-q^{[\Gamma]}}\,
(A\cdot\Gamma)(B\cdot\Gamma)
(C\cdot\Gamma)\,n_\Gamma
+\frac{q^{-[\Gamma]}}{1-q^{-[\Gamma]}}
(A\cdot\Gamma)(B\cdot\Gamma)
(C\cdot\Gamma)\,n_\Gamma \\
=&\left(\frac{q^{[\Gamma]}}{1-q^{[\Gamma]}}
+\frac{1}{q^{[\Gamma]}-1}\right)
(A\cdot\Gamma)(B\cdot\Gamma)
(C\cdot\Gamma)\,n_\Gamma                    \\
=&-(A\cdot\Gamma)(B\cdot\Gamma)
(C\cdot\Gamma)\,n_\Gamma.
\end{split}\end{equation}
Adding 
eqs.~\eqref{eq:three} and \eqref{eq:four}
proves the desired formula.

The conclusion from all of this should be that the mirror symmetry
phenomenon is really about birational equivalence classes.
For if there is any analytic continuation of the correlation function
from the region associated to $\cK_J$
out into the next cone
$\widehat{\cK}_{\widehat J}$,
the calculation above shows that this analytic continuation
must in fact reproduce the correlation function of the flopped model
$\widehat X_{\widehat J}$.

It is tempting to think that if we combined the $\sigma$-model moduli
spaces for all birational models of $X$ we would fill out the
entire conformal field theory moduli space.  However, some examples
that have been worked out by Witten \cite{phases} and by Aspinwall,
Greene and the author \cite{catp}
show that this is not the case.
In those examples, there are other
regions in the moduli space which correspond to rather different
kinds of physical model, including some
called {\em Landau--Ginzburg theories}\/ which will play a r\^ole again
in the next lecture.

\chapter*{}
\lecturename{Mirror Symmetry}
\lecture

\markboth{D. R. Morrison, Mathematical Aspects of Mirror Symmetry}{Lecture
8. Mirror Symmetry}

\section{Mirror manifold constructions}

The original speculations about mirror symmetry were based on the
appearance of arbitrariness of a choice that was made in identifying
certain constituents of the conformal field theory associated to a
Calabi--Yau manifold with geometric objects on the manifold.
The distinction between vertex operators which appear in the $A$-model and
$B$-model correlation functions is simply a difference in sign of a certain
quantum number; if that sign is changed,
the geometric interpretation is altered dramatically.  This led
Dixon \cite{Dixon}
and Lerche--Vafa--Warner \cite{LVW} to propose that there might be a second
Calabi--Yau manifold
producing essentially the same physical theory as the first, but
implementing this change of sign.

Some time later,\footnote{At about the same time, another important
piece of evidence for mirror symmetry was given by Candelas, Lynker,
and Schimmrigk \cite{CLS}, who found an almost perfect symmetry under the
exchange $h^{1,1}\leftrightarrow h^{2,1}$ on the set of  Hodge
numbers coming from
 Calabi--Yau threefolds which can be realized as weighted
projective hypersurfaces.} an explicit construction was made by Greene and
Plesser \cite{GreenePlesser}
which showed that this phenomenon does indeed occur in physics.
The construction rests on a chain of equivalences which are believed to
hold among
different physical models, as follows.

\begin{enumerate}
\item
Certain $\sigma$-models on Calabi--Yau manifolds are believed to correspond
to so-called Landau--Ginzburg theories \cite{GVW}.  (It has recently been
recognized \cite{phases} that this
correspondence is not direct, but involves analytic continuation
on the moduli space.)  Roughly speaking, the class of Calabi--Yau
manifolds for which this correspondence can be made is the class
of ample anti-canonical hypersurfaces in toric varieties.
Such a hypersurface will have an equation of the form
$\Phi(x_1,\dots,x_{n+1})=0$ (in some appropriate coordinates on the
torus), and this same polynomial is used as a ``superpotential''
in constructing the Landau--Ginzburg theory.

\item
Certain Landau--Ginzburg theories---quotients of
the ones for which the superpotential
is of ``Fermat type''
\[\Phi(x_1,\dots,x_{n+1})=x_1^{d_1}+\dots+x_{n+1}^{d_{n+1}}\]
by certain finite groups $\Gamma$---are
believed to correspond to yet another type of conformal field
theory.  This other theory is described in terms of discrete
series representations $V^{(k)}$ of the ``$N{=}2$ superconformal
algebra,'' and it takes the form
\[\left(\bigotimes_j V^{(d_j+2)}\right)/G\]
where $G$ is a slight enlargement of the group $\Gamma$.  (Note that
the case of $\Gamma$ being trivial is allowed, but then $G$ is not trivial.)

The representation theory of the $N{=}2$ superconformal
algebra is related to these things by analyzing the conformal field theory
 on an infinite
cylinder.
(The superconformal algebra can be described in terms of automorphisms of the
cylinder.)

\item
By studying the representation theory, Greene and Plesser
 find a kind of duality
among the finite groups $G$:  there is a dual group $\widehat{G}$
and an isomorphism
\[\left(\bigotimes_j V^{(d_j+2)}\right)/G\cong
\left(\bigotimes_j V^{(d_j+2)}\right)/\widehat{G}\]
which has the ``sign-reversing property'' of mirror symmetry.

\item
The duality can be extended to the groups $\Gamma$, and the mirror
Landau--Ginzburg theory of $\Phi/\Gamma$ is $\Phi/\widehat{\Gamma}$.
This looks a bit asymmetric, since for example the case $\Gamma$ trivial
leads to a rather large group $\widehat{\Gamma}$.  But the group
$\Gamma$ continues to act as a group of ``quantum symmetries''
on the quotient theory, in a way that restores symmetry to this
construction.

\item
Finally, the Calabi--Yau which is the quotient of the Fermat hypersurface
by $\Gamma$ should have as its mirror the one which is the quotient by
$\widehat{\Gamma}$.

\end{enumerate}
This is called the {\em Greene--Plesser orbifolding construction}.

\medskip

There is a conjectural generalization of this construction, which
as of yet has no basis in conformal field theory---it is simply
a mathematician's guess.  This generalization would work for an
arbitrary
family of Calabi--Yau hypersurfaces in toric varieties.
The construction is due to
V.~Batyrev \cite{batyrev1}.\footnote{We restrict ourselves
to the hypersurface case here; further generalizations---to
complete intersections---were subsequently given by Borisov and
Batyrev--Borisov \cite{borisov,BB:dual}.}

Take an ample anticanonical hypersurface $M$ in a toric variety
$V$, and let $\{M_t\}$ be the family of such.  This family is
determined by the Newton polygon of the corresponding equations---that is
a polygon $P\subset L_\R:=L\otimes \R$, where
$L$ is the {\em monomial lattice}\/ of the torus $T$ (of which $V$ is
a compactification).

Batyrev shows that the Calabi--Yau condition admits a particularly
simple characterization in terms of $P$:  the polyhedron $P$
is {\em reflexive}, which means that each hyperplane $H$ which
supports a face of codimension one of $P$ can be written in the
form
\[H=\{y\in L_\R\suchthat (\ell,y)=-1\}\]
for some appropriate vector $v\in \Hom(L,\Z)$.  (The key property
here is the {\em integrality}\/ of the vector $v$---there would
always be some $v\in \Hom(L,\R)$ to define $H$.)

\begin{lemma}[Batyrev]
If $P$ is reflexive, then the {\em polar polyhedron}\/
\[P^o:=\{x\in \Hom(L,\R)\suchthat (x,y)\ge-1\text{ for all }y\in P\}\]
is also reflexive.
\end{lemma}

The conjectured generalization is that the mirror of the family $\{M_t\}$
of hypersurfaces determined by $P$ should be the family $\{W_s\}$
of hypersurfaces (in a compactification of the dual torus of $T$)
determined by the polar polyhedron $P^o$.

One of the pieces of evidence for this conjecture is

\begin{theorem}[Batyrev]\label{batthm}
\[\dim H^{\pm1,1}(\widehat{M})=\dim H^{\mp1,1}(\widehat{W}),\]
where $\widehat{M}$ and $\widehat{W}$ are $\Q$-factorial terminalizations
of $M$ and $W$ respectively.
\end{theorem}

\noindent
Batyrev and collaborators have also explored the Hodge structures of these
hypersurfaces in considerable detail \cite{Bat:vmhs,BvS,BatCox}.

A refinement of Batyrev's theorem 
called the {\em monomial-divisor
mirror map}\/ was introduced in \cite{mondiv}.  This map
gives an explicit combinatorial correspondence between (appropriate
subspaces of) $H^{\pm1,1}(\widehat{M})$ and $H^{\mp1,1}(\widehat{W})$,
and is expected to correctly determine the derivative of the mirror
map near the large radius limit point.  That derivative data is
precisely what one needs in order to evaluate the ``constants of
integration'' in finding the canonical coordinates $q_j$.

\medskip

There is another mirror manifold construction for a class of threefolds
which has been proposed by
Voisin \cite{Voisin:K3} and Borcea \cite{Borcea:K3}.  Let $S$ be a K3
surface with an involution $\iota$ such that $\iota^*(\Omega)=-\Omega$ for
any holomorphic two-form $\Omega$ on $S$, and let $E$ be an elliptic curve.
The quotient 
$\overline{M}=(S\times E)/(\iota\times(-1))$ has singularities along the
fixed curves of the involution $\iota\times(-1)$, but they can be resolved
by a simple blowing up to produce a Calabi--Yau threefold $M$.

Involutions of this type on K3 surfaces have been classified by Nikulin
\cite{Nikulin:involutions}, who found that they
fall into a pattern with a remarkable
symmetry; when the Hodge numbers of the associated Calabi--Yau threefold
are calculated, this symmetry becomes the expected mirror relation among
Hodge numbers.  The detailed knowledge which is available concerning the
variations of Hodge structure on K3 surfaces can be used to study the
correlation functions in detail for these models \cite{Voisin:K3}, which
provides further evidence that mirror partners have been correctly
identified.  In fact, there is also a physics argument explaining why these
pairs of conformal field
theories are actually mirror to each other \cite{AM:K3}, based on
the physics of mirror symmetry for K3 surfaces.

\section{Hodge-theoretic mirror conjectures}

We can now formulate the main conjecture in the mathematical study of
mirror symmetry.

\begin{HTmirrorconjecture}
\quad Given a boundary point $P\in\overline{\MM}_W$
with maximally unipotent monodromy (or perhaps with strongly
maximally unipotent monodromy), there should exist
 a mirror partner $M$ of $W$,
a framing $\sigma$ of $M$, a neighborhood $U$ of $P$ in $\MM_W$,
and a ``mirror map''
\[\mu:U \to (\cD_\sigma/L)^-\]
which is determined up to constants of integration by the property that
\[\mu^*(\dlog q_j)=
d\left(\frac{\int_{\gamma_j}\Omega}{\int_{\gamma_0}\Omega}\right),\]
such that $\mu$ induces an isomorphism between appropriate sub-variations
of Hodge structure of
\begin{enumerate}
\item the formal completion of the geometric variation of Hodge structure
at $P$, and
\item the framed $A$-variation of Hodge structure with framing $\sigma$.
\end{enumerate}
(The sub-variations of Hodge structure should contain the entire first two
terms of the Hodge filtration on both sides.)
\end{HTmirrorconjecture}

There are additional conjectures one wants to make about the relationship
between $M$ and $W$:  there should also be isomorphisms
\[H^{p,q}(W)\cong H^{-p,q}(M)\quad p\ge0,\]
and these should preserve all correlation functions.  (In particular, the
``reverse'' mirror isomorphism should hold, and there should also
be isomorphisms between correlation functions which do not come
from variations of Hodge structure.)  Of course, such isomorphisms
only make sense if we have specified the constants of integration.
In fact, one wants to conjecture that the entire conformal field theory
moduli spaces are isomorphic, but this is a difficult conjecture to
make precisely at present since we do not have a complete mathematical
understanding of conformal field theory moduli spaces.

If we start with the $A$-variation of Hodge structure, there is another
conjecture we can make.

\begin{converse}
Conversely, given $(M,\sigma)$, the corresponding $A$-variation of Hodge
structure comes from
geometry, in the sense that there is a family $\cZ\to\overline{S}$ of
varieties degenerating at $0\in\overline{S}$ such that the framed $A$-variation
of Hodge structure
is isomorphic to the formal completion at $0$ of a (Tate-twisted)
sub-variation of Hodge structures of the variation of Hodge structures
on some cohomology of $Z_s$.
\end{converse}

Due to the phenomenon of rigid Calabi--Yau manifolds, we
can't assume any stronger properties about $Z_s$: Calabi--Yau threefolds
with $h^{2,1}=0$ cannot have mirror partners in the usual sense,
since such a mirror partner would satisfy $h^{1,1}=0$, which is absurd.
However,
there is an example in the physics literature of
a rigid Calabi--Yau manifold, known as the ``$Z$-orbifold,''
which has a mirror physical theory that was
  worked out recently by Candelas, Derrick and Parkes \cite{CDP}
(see also \cite{AspGr}).
In this example, the variation of Hodge structure
associated to the mirror theory
can be described by the family of cubic sevenfolds in $\P^8$
(with a suitable Tate twist).

\section{Some computations}

We explain some of the evidence in favor of the mirror symmetry
conjectures which has been accumulated through specific
computations.\footnote{The computations presented here are taken from
the original paper of Candelas, de la Ossa, Green and Parkes \cite{CDGP} on
the quintic threefold, and a paper of Greene, Plesser and the
present author \cite{GMP}
on higher dimensional mirror manifolds.  A survey of other calculations
of this type (and the methods for making them) can be found in
\cite{predictions}.}
We will compute with Calabi--Yau hypersurfaces of dimension $n\ge3$
in ordinary projective
space $\CP^{n+1}$; the degree of the hypersurface must be $n+2$.
The family of  such hypersurface includes a Fermat hypersurface,
which is part of the ``Dwork pencil'' with defining equation:
\[x_0^{n+2}+\cdots+x_{n+1}^{n+2}-(n+2)\psi\,x_0{\cdots}x_{n+1}=0,\]
where $\psi$ is a parameter.
The group
\[\Gamma:=\{(\alpha_0,\dots,\alpha_{n+1})\suchthat
\alpha_j\in\mmu_{n+2}, \prod\alpha_j=1\}/\{(\alpha,\dots,\alpha)\}\]
acts on the fibers of this family by componentwise
multiplication on the coordinates.

Using either the Greene--Plesser orbifolding construction, or Batyrev's
polar polyhedron construction, one sees that the family $\{M\}$
of hypersurfaces of degree $n+2$ in $\CP^{n+1}$ has as its predicted mirror
the family $\{W\}$ described as the Dwork pencil modulo $\Gamma$
(living in the quotient space $\CP^{n+1}/\Gamma$).
In fact, we can describe the moduli space of this mirrored family
in terms of the parameter $\psi^{n+2}$---the reason for passing
to a power is the existence of an additional automorphism, acting
on the family as a whole, generated by componentwise multiplication
in the $x$'s by $(\alpha,1,\dots,1)$ while simultaneously multiplying
$\psi$ by $\alpha^{-1}$.

It is not difficult to compute where this family becomes singular.
The partial derivatives of the defining equation are all of the form
\[(n+2)\left(x_j^{n+1}-x_j^{-1}\psi\,x_0{\cdots}x_{n+1}\right)\]
and for these to vanish simultaneously we must have $\psi^{n+2}=1$.
Moreover, the additional automorphism of the family fixed the fiber
$\psi=0$, and so causes additional singularities there.
Thus, we can describe the moduli space as $\CP^1-\{0,1,\infty\}$,
with its natural compactification being $\CP^1$.

What is the monodromy behavior at the boundary points?  (We label the
monodromy transformations according to the point.)  At $\psi^{n+2}=0$,
we find that the monodromy has finite order, at $\psi^{n+2}=1$ it is unipotent
but $(T_1-I)^2=0$ so the order is not maximal (since $n\ne1$),
and at $\psi^{n+2}=\infty$ we find maximal order of unipotency
 $(T_\infty -I)^n\ne0$.  In fact, this point is maximally unipotent,
and even strongly maximally unipotent, in the terminology established
earlier.

To compute canonical coordinates and correlation functions near
$\psi^{n+2}=\infty$
we need to know the period functions there.  These can be found by
studying the differential equations which they satisfy.  In this case
of toric hypersurfaces, we have a special method available---the
representation of cohomology by means of residues of differential
forms on the ambient space with poles along the hypersurface.
A basis for the primitive cohomology can be written (in the affine
chart $x_0=1$, say) as
\[\beta_j:=\Res\left(
\frac{\psi^{j+1}\,(x_1{\cdots}x_{n+1})^j\,dx_1\wedge\cdots\wedge dx_{n+1}}
{\left(1+x_1^{n+2}+\cdots+x_{n+1}^{n+2}
-(n+2)\psi\,x_1{\cdots}x_{n+1}\right)^{j+1}}
\right)\]
The connection matrix in this basis can then be found using Griffiths'
``reduction of pole order'' lemma \cite{Griffiths}
to calculate coefficients $\theta_{ij}$
such that
\[\nabla(\beta_i)=\sum\theta_{ij}\beta_j.\]
To find the period matrix from the connection matrix, one must solve
some differential equations.  For if $\{e_k\}$ is a basis for the
local system and we write $e_k=\sum\eta_{ki}\beta_i$ then
\[0=\nabla(e_k)=\sum d\eta_{ki}\, \beta_i + \sum\eta_{ki}\theta_{ij}\beta_j\]
gives differential equations for the unknown coefficient functions
$\eta_{ki}$:
\[d\eta_{ki}=-\sum\eta_{k\ell}\theta_{\ell i}.\]
The flatness of $\nabla$ is equivalent to the integrability
of these equations, which can therefore be solved.

\begin{table}
\begin{center}
\begin{tabular}{|l|l|} \hline
$n$&$n$-point function\\ \hline
$3$&
$5+2875\,q+4876875\,q^2+8564575000\,q^3+15517926796875\,q^4
$\\ &$\phantom{5}
+28663236110956000\,
q^5+53621944306062201000\,q^6
$\\ &$\phantom{5}
+101216230345800061125625\,q^7+
192323666400003538944396875\,q^8
$\\ &$\phantom{5}
+367299732093982242625847031250\,q^9
$\\ &$\phantom{5}
+704288164978454714776724365580000\,q^{10}
$\\ &$\phantom{5}
+1354842473951260627644461070753075500\,q^{11}
$\\ &$\phantom{5}
+2613295702542192770504516764304958585000\,q^{12}
$\\ &$\phantom{5}
+5051976384195377826370376750184667397150000\,q^{13}
$\\ &$\phantom{5}
+9784992122065556293839548184561593434114765625\,q^{14}
$\\ &$\phantom{5}
+18983216783256131050355758292004110332155634496875\,q^{15}
$\\ &$\phantom{5}
+36880398908911843175757970052077286676680907186572875\,q^{16}
$\\ &$\phantom{5}
+71739993072775923425756947313710004388338109828244718125\,q^{17}
$\\ &$\phantom{5}
+139702324572802672116486725324237666156179096139345867681250\,q^{18}
$\\ &$\phantom{5}
+\dots$\\[6pt]
$4$&
$6 + 120960 \,q \!+\! 4136832000 \,q^2  \!+\! 148146924602880 \,q^3  \!+
    5420219848911544320 \,q^4
$\\ &$\phantom{6}
+ 200623934537137119778560 \,q^5
+    7478994517395643259712737280 \,q^6
$\\ &$\phantom{6}
+ 280135301818357004749298146851840 \,q^7
$\\ &$\phantom{6}
+    10528167289356385699173014219946393600 \,q^8
$\\ &$\phantom{6}
+    396658819202496234945300681212382224722560 \,q^9
$\\ &$\phantom{6}
+    14972930462574202465673643937107499992165427200 \,q^{10}
$\\ &$\phantom{6}
+    566037069767251121484562070892662863943365345190400 \,q^{11}
$\\ &$\phantom{6}
+    21424151141341932048068067497996096856724987411324108800 \,q^{12}
\!+\! \dots$\\[6pt]
$5$&
$7 + 3727381 \,q + 2637885990187 \,q^2  + 1927092954108108787 \,q^3
$\\ &$\phantom{7}
+1425153551321014327663291 \,q^4  + 1060347883438857662557634869906 \,q^5
$\\ &$\phantom{7}
+    791661306374088776109692880989252173 \,q^6
$\\ &$\phantom{7}
+    592348256908461616176898022359492565546566 \,q^7
$\\ &$\phantom{7}
+    443865568545713063761643598030194801299861575595 \,q^8
$\\ &$\phantom{7}
+    332947403131697202086626568381790256001850741509664373 \,q^9
+\dots$\\[6pt]
$6$&
$8 + 106975232 \,q + 1672023727001600 \,q^2  + 26611692333081695092736 \,q^3
$\\ &$\phantom{8}
+    426129121674687823674948571136 \,q^4
$\\ &$\phantom{8}
+    6842148599241293047857339542861643776 \,q^5
$\\ &$\phantom{8}
+    110018992594692024449889564415904439556898816 \,q^6
$\\ &$\phantom{8}
+    1770551943055574073245974844490813198478975912902656 \,q^7
$\\ &$\phantom{8}
+    28508925683951911989843155602330000507452539542539447947264 \,q^8
$\\ &$\phantom{8}
+\cdots$\\
\hline
\end{tabular}
\end{center}

\medskip

\caption{$n$-point functions in dimension $n$}

\end{table}

If we work in a local coordinate $z=\psi^{-n-2}$ near $\psi^{n+2}=\infty$,
we find that a basis $e_0(z)$, \dots, $e_{n+1}(z)$ of local solutions
can be found such that $e_0(z)$ is single-valued near $z=0$, and
\[e_{j+1}(z) = (\log z)\, e_j(z) + \text{single-valued function}.\]
(This is a consequence of the maximally unipotent monodromy.)
The vectors $e_j(z)$ form the columns of the period matrix.

One can then use row operations to put the period matrix in upper
triangular form, with constant diagonal elements.  (Let us choose the
diagonal elements to all be $n+2$.)  This implements the change of basis
to a basis consisting of distinguished sections of $\cH^{p,q}_\BS$.
The nonzero entries $A_j^1$ in the
connection matrix are then calculated by differentiating rows of the
period matrix, and writing the result as a multiple of a subsequent
row.  Each such entry takes the form
\[A^1_j=Y^1_j\,\frac{dq}q,\]
and the functions $Y^1_j$ represent correlation functions $\langle
(\partial/\partial t)\,\beta_j\,\beta_{n-j-1}\rangle$.

This can all be done very explicitly, using power series expansions
of the unknown single-valued functions, in these examples.
(I advise using {\sc maple} or {\sc mathematica} if
you would like to try it for yourself.)
We show two kinds of
calculations in the tables.  For the first, only the ``maximally
unipotent'' assumption is required, since the calculation requires
only the distinguished $n$-form and the canonical coordinates.
What is computed in table 1 is
 the ``$n$-point function,''
which iterates the differential of the period map $n$ times.
(This was introduced some years ago  in the variation of Hodge structures
context by Carlson, Green, Griffiths
and Harris \cite{CGGH}.)

\begin{table}
\begin{center}
\begin{tabular}{|l|} \hline
$Y_1^1=
5+2875\,\cuone3+609250\,\cu23+317206375\,\cu33+242467530000\,\cu43
$\\$\phantom{Y_1^1=5}
+229305888887625\,\cu53+
248249742118022000\,\cu63
$\\$\phantom{Y_1^1=5}
+295091050570845659250\,\cu73+375632160937476603550000\,\cu83
$\\$\phantom{Y_1^1=5}
+503840510416985243645106250\,\cu93
$\\$\phantom{Y_1^1=5}
+704288164978454686113488249750\,\cu{10}3
$\\$\phantom{Y_1^1=5}
+1017913203569692432490203659468875\,\cu{11}3
$\\$\phantom{Y_1^1=5}
+1512323901934139334751675234074638000
\,\cu{12}3
$\\$\phantom{Y_1^1=5}
+2299488568136266648325160104772265542625\,\cu{13}3
$\\$\phantom{Y_1^1=5}
+3565959228158001564810294084668822024070250\,\cu{14}3
$\\$\phantom{Y_1^1=5}
+5624656824668483274179483938371579753751395250\,\cu{15}3
$\\$\phantom{Y_1^1=5}
+9004003639871055462831535610291411200360685606000\,\cu{16}3+\dots
$\\
\hline
\end{tabular}
\end{center}

\medskip

\caption{Three-point function in dimension three}

\end{table}

\begin{table}
\begin{center}
\begin{tabular}{|l|} \hline
$Y_1^1=6+60480\,\cuone2+440884080\,\cu22+6255156277440\,\cu32$\\
$\phantom{Y_1^1=6}+117715791990353760\,\cu42
+2591176156368821985600\cdot5^2\,\cu52
+\dots$\\
\hline
\end{tabular}
\end{center}

\medskip

\caption{Three-point function in dimension four}

\end{table}

\begin{table}
\begin{center}
\begin{tabular}{|l|} \hline
$Y_1^1=7+1009792\,\cuone2+122239786088\,\cu22
+30528671745480104\,\cu32$\\
$\phantom{Y_1^1=7}+10378199509395886153216\,\cu42
+\dots$\\[6pt]
$Y_2^1=7+1707797\,\cuo1+510787745643\,\cuo2
+222548537108926490\,\cuo3$\\
$\phantom{Y_2^1=7}+113635631482486991647224\,\cuo4
+\dots$\\
\hline
\end{tabular}
\end{center}

\medskip

\caption{Three-point functions in dimension five}

\end{table}

\begin{table}
\begin{center}
\begin{tabular}{|l|} \hline
$Y_1^1=8+15984640\,\cuone2+33397159706624\,\cu22
+154090254047541417984\,\cu32
$\\$\phantom{Y_1^1=8}
+1000674891265872131899670528\,\cu42+\dots$\\[6pt]
$Y_2^1=8+\!37502976\,\cuo1\!+\!224340704157696\,\cuo2
\!+\!2000750410187341381632\,\cuo3$\\
$\phantom{Y_2^1=8}
+21122119007324663457380794368\,\cuo4+\dots$\\[6pt]
$Y_2^2=8+\!59021312\,\q\!+\!821654025830400\,\qq
\!+\!\!12197109744970010814464\,\qqq$\\
$\phantom{Y_2^2=8}
+186083410628492378226388631552\,\qqqq+\dots$\\
\hline
\end{tabular}
\end{center}

\medskip

\caption{Three-point functions in dimension six}

\end{table}

The other computations, displayed in tables 2--5, are of
 three-point functions $Y^a_b$, read off of the connection matrix in
a distinguished basis.
(There is a symmetry $Y^a_b=Y^a_{n-a-b}=Y^b_{n-a-b}$ so we only
show some of these.)  The coefficients in the series expansions are the
predicted values of the Gromov--Witten invariants.
The three-point function $Y^1_0$ has
the value $n+2$ (a constant, due to the definition of canonical
coordinates)  and is not shown in the tables.  The other
functions $Y^1_j$ come directly from the connection matrix.
In dimension six, there is also a ``secondary'' function, which (by the
$B$-model version of the associativity, which is simply the associativity
of the ``sheaf cup product'' pairing)
can be calculated as $Y^2_2 =(Y^1_2)^2/Y^1_1$.

There is a relation between the computations in table 1, and those
in tables 2--5, which can be explicitly verified from these tables:  it is
\[\text{$n$-point function } =
\frac{Y^1_0\cdot Y^1_1\cdot {} \dotsm {} \cdot Y^1_{n-1}}{(n+2)^{n}}.\]

The functions $Y^a_b$
are predicted to agree with quantum products
on the mirror manifolds
\[\zeta^a\star\zeta^b\star\zeta^{n-a-b},\]
where $\zeta^j$ is the class of a linear space (in $\CP^{n+1}$)
of complex codimension $j$.  In fact, we have displayed things in tables
2--5 with this in mind, writing series in terms of $q^k/(1-q^k)$.

Also in tables 2--5, we have pulled out some factors of the degree of
the rational curve.  If there are $\ell$ occurrences of ``1'' among
$\{a,b,n-a-b\}$, then there will be $\ell$ of the linear spaces of
codimension one, and each meets a given rational curve $\Gamma$
in $\deg(\Gamma)$ points, giving rise to a factor of $(\deg(\Gamma))^{\ell}$
in the Gromov--Witten invariants.
Pulling out those factors makes the comparison with ``counting'' problems
more transparent.

All of the predicted Gromov--Witten invariants in degrees one and two
in these tables have
been verified by  Katz \cite{katz:verifying};
most of the invariants in degree three have
been verified by Ellingsrud and Str{\o}mme \cite{ES,ESii}.

\chapter*{}

\lecturename{Postscript: Recent Developments}
\lectureoptionstar{POSTSCRIPT:}{Recent Developments}

\markboth{D. R. Morrison,
Mathematical Aspects of Mirror Symmetry}{Postscript: Recent Developments}

As mentioned in the introduction, the subject of mirror symmetry is a
rapidly developing one, and much has happened since the lectures on which
these notes are based were delivered.  We will briefly sketch some of these
developments in this postscript.

The Gromov--Witten invariants and their generalizations have been studied
particularly intensively.  The definition of Ruan \cite{ruan} which we
presented in the lectures has been supplanted by other definitions drawn
from symplectic geometry (cf.~\cite{MS,RuanTian}) which work directly in
cohomology (avoiding the bordism technicalities)
and are also more general.  In full generality these extended
Gromov--Witten invariants are not only
associated to curves of genus zero with three
vertex operators, but also to curves of arbitrary genus $g$ with $k$ vertex
operators (provided that $2g-2+k>0$) and even to some non-topological
correlation functions.\footnote{There have also been investigations
into the physical interpretation of these
higher genus invariants, and how they should transform under mirror
symmetry (in the case of Calabi--Yau threefolds) \cite{BCOV:anom,BCOV:KS}.
At one time, it had been expected that for Calabi--Yau threefolds
the genus zero topological correlation functions would completely determine the
conformal field theory, but now it is known that higher genus invariants
are needed as well \cite{chiral}.}
There are at least three proofs of the associativity relations for these
symplectic Gromov--Witten invariants \cite{RuanTian,Liu,MS},
including proofs of a stronger form of associativity known as the
Witten--Dijkgraaf--Verlinde--Verlinde (WDVV) equations
\cite{topgrav,DVV,Wit:twoDgrav,Dubrov} which are relevant in the
case of higher genus.  As in the genus zero case, these higher genus
invariants can be used to encode a kind of quantum cohomology ring
(somewhat larger than the one we studied here);
it is also possible to interpret
the WDVV associativity relation as the flatness of a certain connection
\cite{Dubrov}.
A very accessible exposition of this circle of ideas has been written by
McDuff and Salamon \cite{MS}.

Parallel to this development, Gromov--Witten invariants have also been defined
purely within algebraic geometry.  The methods of Katz described in the
lectures were developed further (see \cite{katz:GW} and the appendix to
\cite{BCOV:anom}),
and similar methods based
on the construction of a ``virtual moduli cycle'' were developed
independently by Li and Tian \cite{LiTian}.  The foundations for an
algebraic theory of Gromov--Witten invariants were carefully laid by
Kontsevich and Manin \cite{KM} (again, the higher genus invariants and the
WDVV equations play an important r\^ole), and the program they initiated was
ultimately carried out \cite{BehMan,BehFant,Beh}, producing a definition of
Gromov--Witten invariants based on stable maps.  (The work of Li--Tian
mentioned above \cite{LiTian} is also closely related to this program.)
Even before this program
was complete, Kontsevich had applied it to obtain some spectacular results
in enumerative geometry, including a verification of the predicted number
$242467530000$ of rational quartics on the general quintic threefold
\cite{Kontsevich}.  The stable map theory is nicely explained, with further
references, in \cite{FulPan}.

Kontsevich has also formulated a ``homological'' version of the
mirror conjecture \cite{kont:icm}
involving what are 
known as $A^\infty$-categories (cf.~\cite{stasheff}), 
which is related to the
``extended moduli space'' introduced by Witten \cite{witten:mirror}.  
By a construction of Fukaya \cite{fukaya}, to 
every compact symplectic manifold $(Y,\omega)$
with vanishing first Chern class, one
can associate an $A^\infty$-category whose objects are essentially the
Lagrangian submanifolds of $Y$, and whose morphisms are determined by the
intersections of pairs of submanifolds.  Kontsevich's conjecture relates
the bounded derived category of the Fukaya category of $Y$ (playing the
r\^ole of the $B$-model) to the bounded derived category of the category of
coherent sheaves on a mirror partner $X$ (playing the r\^ole of the
$A$-model). I must refer the reader to \cite{kont:icm} for further details
concerning this fascinating conjecture.

The art of making predictions about enumerative geometry from calculations
with the variation of Hodge structure on a candidate mirror partner has
been considerably refined: see \cite{predictions} for a survey and
references to the literature.  The era of numerical experiments in mirror
symmetry seems to be largely over, and has been supplanted by a more
analytical period.  Witten's analysis of the physics related to Calabi--Yau
manifolds which are hypersurfaces in toric varieties \cite{phases}
was further developed
in \cite{summing}, where techniques were found---somewhat related to methods
introduced by Batyrev \cite{Bat:qcoho} for the study of quantum cohomology
of toric varieties---for precisely calculating a
variant of the quantum cohomology ring
of the Calabi--Yau manifold.  (The variant is derived from enumerative
problems on the ambient space rather than directly on the Calabi--Yau
manifold.)
There is a physics argument, but not a complete
mathematics argument, which explains why this variant should coincide with
the usual quantum cohomology ring after a change of coordinates in the
coefficient ring.
This variant {\em can}\/ be rigorously shown to
agree with the correlation functions of the
mirror Calabi--Yau manifold, again calculated in the ``wrong'' coordinates.
In this way, the results of \cite{summing} provided
the first analytical proof that some kind of enumerative
problem on one side of the mirror could be related to a variation of Hodge
structure calculation on the other side.  Further development of these
ideas in \cite{towards-duality} led to a preliminary argument to the effect
that the
physical theories associated to a Batyrev--Borisov pair should actually be
mirror to each other.

In a striking recent development, Givental has proved
\cite{Givental:homological,Givental:ICM,Givental:equivariant}
that for Calabi--Yau
complete intersections in projective spaces, the ``predicted'' enumerative
formulas which one calculates by using a Batyrev--Borisov candidate
mirror partner
are in fact correct evaluations of
the Gromov--Witten invariants.  This establishes, for example,
the accuracy of {\em
all}\/ of the predictions about the general quintic threefold made by
Candelas et al. \cite{CDGP} (and which we listed
in table 2).  Givental's remarkable proof actually has very little to
do with mirror symmetry {\em per se}: in studying an equivariant version of
quantum cohomology, he finds enough structure to enable a calculation which
is formally similar to (and certainly inspired by) the variation of Hodge
structure calculations on the candidate mirror partner.

The last several years have also been a period of dramatic developments in
string theory.  There are new techniques which go by the names of
``duality'' and
``nonperturbative methods,'' and a number of the recent results
 have been closely related
to Calabi--Yau manifolds and mirror symmetry.  One of the earliest
nonperturbative results \cite{Str:,bhole}
was the discovery\footnote{This  had been
anticipated some time earlier in the physics literature
\cite{CDLS,GreenHubsch,CGH,Cd:con} based on the discovery of and
speculations about conifold transitions in the mathematics literature
\cite{Clemens:double,Friedman:simult,Hirzebruch:examples,tianyau,%
reid,Friedman:threefolds}, but an understanding of the physical mechanism
behind the attachment of the moduli spaces
was lacking.} that the string theory moduli spaces associated to
Calabi--Yau manifolds should be attached along loci corresponding to
``conifold transitions''---a process in which a collection of rational
curves is contracted to ordinary double points and the resulting space is
then smoothed to produce another Calabi--Yau manifold.  This new attaching
procedure supplements, but is rather different from, the gluing of
K\"ahler cones which we discussed in section \ref{sec73}.  In the new
procedure, a moduli space of a different dimension (corresponding to a
Calabi--Yau manifold with different Hodge numbers than the original) is
cemented on at the same point where the two
like-dimensional pieces (K\"ahler cones differing by a
flop) have been glued together.  The ``cement'' which holds these two
spaces
 together (i.e., the physical process responsible) is a phase transition
between charged black holes on one component of the moduli space and
elementary particles on the other.

The string theory moduli spaces mentioned above are actually somewhat
larger than the conformal field theory moduli spaces which were one of the
primary subjects
of these lectures.  There are two variants of string theory
which are relevant, called type IIA and type IIB string theories, and the
additional parameters which must be added to the conformal field theory
moduli space differs between the two.  In the case of type IIA, the extra
parameters are a choice of holomorphic $3$-form and the choice of an
element in the intermediate Jacobian of the Calabi--Yau threefold.
(Some of the mathematical structure of these spaces related to the
intermediate Jacobians was anticipated in work of Donagi and Markman
\cite{DonMark}.)
In the case of type IIB, the new parameters are similar, but related to the
even cohomology of the manifold.
These two types of parameters should
be mapped to each other under mirror symmetry \cite{udual,mirrorII}.
In fact, a large number of other related structures called ``D-brane moduli
spaces'' should also correspond under mirror symmetry---the precise
implications of this correspondence (which appears to be connected to
Kontsevich's homological mirror symmetry conjecture) are
still being worked out.

Finally, in a very exciting recent development, a completely new geometric
aspect of mirror symmetry has been discovered by Strominger, Yau and Zaslow
\cite{SYZ}.  A Calabi--Yau manifold $X$ of real dimension $2n$ on which a
complex structure $J$ and K\"ahler form $\omega$
have been fixed has a natural class of
$n$-dimensional submanifolds $M$ defined by the property that
$\omega|_M\equiv0$ and $\Im(\Omega)|_M\equiv0$ for some choice of
holomorphic $n$-form $\Omega$.  These {\em special Lagrangian
submanifolds}\/ were introduced by Harvey and Lawson \cite{HL} as a natural
class of volume-minimizing submanifolds; they have many other interesting
properties, including an exceptionally well-behaved deformation theory
\cite{mclean}.  Strominger, Yau and Zaslow argue on physical grounds (using
the correspondence of D-brane moduli spaces mentioned above) that
whenever $X$ has a mirror partner, then $X$ must admit a map $\rho:X^{2n}\to
B^n$ whose generic fiber is a special Lagrangian $n$-torus, and which has a
section $\sigma:B\to X$ whose image is itself a special Lagrangian submanifold.
Given this structure, the
mirror partner of $X$ is then predicted to be
 a compactification of the family of dual
tori of the fibers of $\rho$.  (The section specifies a point $p_b:=\sigma(b)$
on each torus $T_b:=\rho^{-1}(b)$; the dual torus is then
$\Hom(\pi_1(T_b,p_b),\U(1))$.)
There is also an argument---quite similar in nature
to \cite{towards-duality}---that such a structure should suffice for
producing a mirror isomorphism between the corresponding physical theories.
A mathematical account of this construction can be found in
\cite{underlying}, which attempts to make the mathematical implications of
this story precise: given a ``special Lagrangian $m$-torus fibration,'' all
of the structure we have seen relating the quantum cohomology and the
variation of Hodge structure should (conjecturally) follow as a
consequence.  For the Voisin--Borcea threefolds, the structure of these
special Lagrangian torus fibrations
(using a mildly degenerate metric)
has been worked out in complete detail by Gross and Wilson
\cite{GrossWilson}, who find compatibility
with the previously observed mirror phenomena in a beautiful
geometric form.

\backmatter
\chapter*{}

\lecturename{Bibliography}
\lecturestar{BIBLIOGRAPHY}

\markboth{D. R. Morrison,
Mathematical Aspects of Mirror Symmetry}{Bibliography}

\bibliographystyle{amsplain}

\ifx\undefined\bysame
\newcommand{\bysame}{\leavevmode\hbox to3em{\hrulefill}\,}
\fi

\addvspace\linespacing
\noindent{\large\bfseries Books}\par
\addvspace{.5\linespacing}

\addvspace\linespacing
\noindent{\large\bfseries Mathematics Articles}\par
\addvspace{.5\linespacing}

\addvspace\linespacing
\noindent{\large\bfseries Physics Articles}\par
\addvspace{.5\linespacing}

\vfill

\noindent
{\sc Department of Mathematics, Duke University, Durham NC 27708-0320}

\end{document}